\documentclass{mcom-l}

\newtheorem{theorem}{Theorem}[section]
\newtheorem{lemma}[theorem]{Lemma}

\theoremstyle{definition}

\newtheorem{assumption}[theorem]{Assumption}
\newtheorem{example}[theorem]{Example}

\theoremstyle{remark}
\newtheorem{remark}[theorem]{Remark}

\numberwithin{equation}{section}

\usepackage{changes}
\usepackage[numbers, square, comma]{natbib}
\usepackage{mathrsfs}
\usepackage{algorithm}
\usepackage{algorithmicx}
\usepackage{algpseudocode}
\usepackage{appendix}
\usepackage{url}

\begin{document}

\title[PDFB splitting for Wasserstein-like Gradient Flows]{Efficient Primal-dual Forward-backward Splitting Method for Wasserstein-like Gradient Flows with General Nonlinear Mobilities}


\author{Yunhong Deng}
\address{School of Mathematical Sciences, University of Electronic Science and Technology of China, Chengdu, Sichuan 611731, China; and School of Mathematics, University of Minnesota, Twin Cities, Minnesota 55455, United States}
\email{deng0335@umn.edu}

\author{Chaozhen Wei}
\address{School of Mathematical Sciences, University of Electronic Science and Technology of China, Chengdu, Sichuan 611731, China}
\email{cwei4@uestc.edu.cn}
\thanks{The second author is supported by the National Natural Science Foundation of China under grants 12371392, 12622128 and 12431015.}
\thanks{The second author is the corresponding author for the publication}


\date{January 1, 2001 and, in revised form, June 22, 2001.}


\begin{abstract}
We construct an efficient primal-dual forward-backward (PDFB) splitting method for solving a class of minimizing movement schemes with nonlinear mobility transport distances, and apply it to compute Wasserstein-like gradient flows. Our approach introduces a new saddle point formulation of the minimizing movement schemes, derived using a support function characterization from the Benamou-Brenier dynamical formulation of optimal transport. This framework enables the flexible computation of Wasserstein-like gradient flows and naturally extends to general nonlinear mobilities by solving the corresponding saddle point problem at the fully discrete level. We also present a detailed convergence analysis of the PDFB splitting method, along with practical remarks on its implementation and application. The effectiveness of the PDFB splitting method is demonstrated through several challenging numerical examples.
\end{abstract}

\maketitle
\section{Introduction}\label{sec1}
In this work, we consider Wasserstein-like gradient flows given as follows,
\begin{equation}\label{1-1}
    \begin{aligned}
        \partial_{t}\rho - \nabla\cdot\bigg(M(\rho) \nabla \frac{\delta\mathcal{E}(\rho)}{\delta \rho}\bigg) &= 0 \text{ in } \Omega \times (0, T)\\
        M(\rho) \nabla \frac{\delta\mathcal{E}(\rho)}{\delta \rho} \cdot \boldsymbol{\mathrm{n}} &= 0 \text{ on } \partial\Omega \times (0, T)
    \end{aligned}
\end{equation}
where $\Omega \subseteq \mathbb{R}^{d}$ is a regular bounded domain, and $\textbf{n}$ is the unit outer normal to the domain $\Omega$, and $T > 0$ is a time horizon where the equation is well-posed. In many practical scenarios of physical modeling which we will describe below, $\rho$ represents the density or concentration of some quantity of interests and $M : \mathbb{R} \to \mathbb{R}$ is the concentration-dependent mobility function. To ensure the physical meaning of the equation, we start the initial value $\rho_{0}$ within the region where $M$ is positive in the sense that $M(\rho_{0}) \ge 0$. $\mathcal{E}(\rho)$ is the free energy functional that drives the dynamics of the system. A standard choice is
\begin{equation}\label{1-12}
    \mathcal{E}(\rho) = \underbrace{\int_{\Omega} U(\rho(x))\,\text{d}x}_{\text{internal energy}} + \underbrace{\int_{\Omega} V(x)\rho(x) \,\text{d}x}_{\text{potential energy}} + \underbrace{\frac{\varepsilon^{2}}{2} \int_{\Omega} |\nabla \rho|^{2} \,\text{d}x}_{\text{Dirichlet energy}},
\end{equation}
and $\delta\mathcal{E}(\rho)/\delta \rho$ is the first variation of $\mathcal{E}(\rho)$, defined in a distributional sense:
\begin{equation*}
    \begin{aligned}
        \int_{\Omega} \frac{\delta\mathcal{E}(\rho)}{\delta \rho}(x) \varphi(x) \,\text{d}x &= \lim_{t \to 0} \frac{\mathcal{E}(\rho + t \varphi) - \mathcal{E}(\rho)}{t},\quad \forall \varphi \in C^{\infty}_{0}(\Omega).
    \end{aligned}
\end{equation*}
Our method has the ability to cope with a general energy functional $\mathcal{E}(\rho)$, and we will include other interesting forms of energy functional from the literature in our numerical examples.

Equation~\eqref{1-1} has been employed for modeling in various fields, including biology, physics and materials science, with different mobility functions $M$ considered for different applications. For example, the typical Wasserstein gradient flows correspond to the special case where $M(\rho)=\rho$ \cite{Jordan1996, Carrillo2003, Ambrosio2005}. In chemotaxis models that account for overcrowding prevention, a nonlinear degenerate mobility $M(\rho)=\rho(\sigma - \rho)$ is employed, where $\sigma$ refers to the saturation level of the population density \cite{Burger2006}. In Cahn-Hilliard equations for phase separation, a mobility $M(\rho)=1 - \rho^{2}$ is commonly utilized for modeling the dynamics of surface diffusion \cite{Cahn1996, Elliott1996}. In lubrication models for thin films, the mobility appears in a power form given by $M(\rho)=\rho^{3}$, where $\rho$ represents the thickness of the thin film, and the exponent is derived from the no-slip boundary condition on the fluid-solid interface \cite{Bertozzi1998, Bertozzi2002}.

One major difficulty in solving~\eqref{1-1} numerically is the stability issue, which arises not only from the nonlinear energy functional but also the possible nonlinearity and degeneracy of the mobility $M(\rho)$. To overcome this issue, it is often desirable to design numerical methods that preserve, at the fully discrete level, the underlying energy dissipation structure of the gradient flows \cite{Dolbeault2008, Carrillo2009, Lisini2012}:
\begin{equation}\label{EDI}
    \frac{\text{d}}{\text{d}t}\mathcal{E}(\rho(\cdot, t)) = -\int_{\Omega} M(\rho(x, t))|\textbf{v}(x, t)|^{2} \,\text{d}x \le 0,\quad \textbf{v}(\cdot, t) = -\nabla \frac{\delta \mathcal{E}(\rho(\cdot, t))}{\delta \rho},
\end{equation}
which provides an energy stability property to numerical solutions. Additionally, it is crucial to designing numerical methods to preserve fundamental properties of the physical solutions, such as bound preservation (i.e., the discrete maximal bound principle) for degenerate mobilities \cite{beretta1995nonnegative,Elliott1996}, and mass conservation.

Recently, many structure-preserving numerical methods have been proposed to maintain specific properties for certain types of gradient flows. For example, significant work has focused on energy-stable methods, such as the convex splitting method \cite{Eyre1998}, the stabilized semi-implicit method \cite{Xu2006}, the invariant energy quadratization method \cite{Yang2017}, and the scalar auxiliary variable method \cite{Shen2018}. Bound-preserving methods include, but are not limited to, the cutoff method \cite{Lu2012}, the flux-limiting discontinuous Galerkin method \cite{Frank2020}, the exponential time-differencing method \cite{Du2021Siamreview} and the Lagrangian multiplier method \cite{Cheng20212}. Intricate structure-preserving methods for Wasserstein gradient flows have also been considered in literature. For example, in \cite{Carrillo2014, Bailo2021}, finite-volume methods are proposed, which consider appropriate flux limiters to ensure bound preservation and upwind approximations of velocity to ensure energy dissipation; in \cite{Sun2017, Sun2018}, discontinuous Galerkin methods are proposed, which admit an entropy inequality on the discrete level, and positivity-preserving limiter to obtain nonnegative solutions. However, these methods often require intricate constructions to preserve energy stability, bound preservation and mass conservation simultaneously.

On the other hand, based on its form of energy dissipation given in~\eqref{EDI}, equation~\eqref{1-1} can be seen as a gradient flow with respect to a transport metric induced by the mobility $M(\rho)$ in the following admissible set \cite{Dolbeault2008, Lisini2012}:
\begin{equation*}
    \mathcal{P}_{M}(\Omega) = \bigg\{\rho : \Omega \to \mathbb{R} : \int_{\Omega} \rho(x) \,\text{d}x = m \text{ and } M(\rho) \ge 0\bigg\},
\end{equation*}
where $m$ is a fixed constant referred to as the total mass. Following the variational characterization of gradient flows as minimizing movements \cite{DeGiorgi, Jordan1996}, one can derive a semi-discrete energy-dissipation scheme for solving equation~\eqref{1-1} as a sequence of minimization problems \cite{Carrillo2023, Li2023}. Specifically,  given the approximation $\rho^{n}$ at $t = t^{n}$, the numerical solution at the next time step $t = t^{n+1}$ can be computed by
\begin{equation}\label{JKO}
    \rho^{n + 1} = \arg\min_{\rho} \Delta t \mathcal{E}(\rho) + \frac{1}{2}W_{M}(\rho^{n}, \rho)^{2} \text{ over } \rho \in \mathcal{P}_{M}(\Omega),
\end{equation}
where $\Delta t$ is the time step size and $W_{M}$ is a transport metric induced by the mobility function \cite{Dolbeault2008}, defined in a dynamical formulation as
\begin{equation}\label{1-5}
    \begin{aligned}
        W_{M}(\rho_{0}, \rho_{1}) := \min_{\rho, \textbf{m}} &\bigg\{\int_{0}^{1}\int_{\Omega} 2f\big(M(\rho(x, t)), \textbf{m}(x, t)\big) \,\text{d}x \text{d}t: \\
        &(\rho, \textbf{m}) \in \mathfrak{D}_{\text{ce}}(0, 1),\,
        \rho(\cdot, 0) = \rho_{0} \text{ and } \rho(\cdot, 1) = \rho_{1}\bigg\}^{\frac{1}{2}}\,.
    \end{aligned}
\end{equation}
where $\textbf{m}$ is the momentum and $f$ is a lower semi-continuous, convex, and one-homogeneous function:
\begin{equation}\label{1-4}
    f\big(M, \textbf{m}\big) :=
    \begin{cases}
        \displaystyle\frac{|\textbf{m}|^{2}}{2M} & M > 0\,;\\
        0 & \big(M, \textbf{m}\big) = (0, \boldsymbol{0})\,;\\
        +\infty & \text{ otherwise }\,.\\
    \end{cases}
\end{equation}
which will be referred to as the action function in the following, and $\mathfrak{D}_{\text{ce}}(0, 1)$ is the admissible set given by the constraints of continuity equation and no-flux boundary conditions given as follows:
\begin{equation*}
    \begin{aligned}
        \mathfrak{D}_{\text{ce}}(0, 1) := \Big\{(\rho, \textbf{m}) : \Omega \times (0, 1) &\to \mathbb{R} \times \mathbb{R}^{d} : \partial_{t}\rho(x, t) + \nabla\cdot \textbf{m}(x, t) = 0 \text{ in } \Omega \times (0, 1)\\
        &\text{ and } \textbf{m}(x, t) \cdot \textbf{n}(x) = 0 \text{ on } \partial \Omega \times (0, 1)\Big\}
    \end{aligned}
\end{equation*}
with $\textbf{n}$ being the unit outer normal of the domain $\Omega$.

This temporal discrete scheme~\eqref{JKO} is called the generalized JKO scheme, named after the seminal work \cite{Jordan1996} developed by Jordan, Kinderlehrer and Otto. In~\eqref{JKO}, energy dissipation is guaranteed by the minimization, mass conservation is satisfied by the no-flux boundary condition on $\textbf{m}$, and bound preservation is ensured by the feasibility of the action function $f$ \citep[Proposition 1]{Carrillo2023}. The existence of the time-discrete solution and its weak convergence as $\Delta t \to 0$ to the weak solution of~\eqref{1-1} has been established for linear mobilities \cite{Jordan1996, Blanchet2008}, and concave mobilities \cite{Lisini2012}.

Building upon the JKO formulation~\eqref{JKO}, many structure-preserving numerical methods are proposed for computing ~\eqref{1-1}. These methods typically begin by constructing fully discrete approximations of~\eqref{JKO} and then apply large-scale convex optimization solvers. For Wasserstein gradient flows, common approximation methods include finite difference methods \cite{Carrillo2019, Li2019}, finite volume methods \cite{Cancs2019}, and particle methods \cite{Westdickenberg2008, Lee2023}. For Wasserstein-like gradient flows with nonlinear mobilities, approximations include finite difference methods \cite{Carrillo2023} and finite element methods \cite{Li2023}. However, due to the non-smoothness of the action function $f$ defined in \eqref{1-4}, caused by the degeneracy of mobility, solving the optimization problem~\eqref{JKO} can be non-trivial. For this reason, some methods introduce artificial entropic regularization (such as Fisher information regularization \cite{Li2019}) to enhance smoothness and convexity, but this also introduces additional numerical diffusion at the PDE level. Alternatively, our previous works \cite{Carrillo2019, Carrillo2023} developed a primal-dual method based on proximal splitting techniques (see \cite{Papadakis2013} and the references therein) to solve \eqref{JKO} without the need for artificial terms. However, as we will explain below, its application is limited to gradient flows with simple mobilities such as linear mobilities \cite{Carrillo2019} or quadratic mobilities \cite{Carrillo2023}, due to the challenges involved in computing proximal operators. Developing a general framework that can handle a variety of nonlinear mobilities remains a significant challenge.

In this work, we introduce a new saddle-point formulation of~\eqref{JKO}, and develop a primal-dual forward-backward (PDFB) splitting method for computing~\eqref{1-1} with {\it general} nonlinear mobilities. It avoids the need to develop complex algorithms for computing the proximal operators with respect to different mobilities, and can be efficiently solved by primal-dual iterates with simple-to-evaluate proximal operators. The proposed algorithm can be regarded as an extension of the primal-dual three-operators (PD3O) splitting method \cite{Yan2016} adopted in the previous works \cite{Carrillo2019, Carrillo2023}, and it shares many similarities with primal-dual methods \cite{Hamedani2018, Malitsky2018, Zhu2022} for solving general convex-concave problems. We provide rigorous convergence analysis of the PDFB splitting method, and some remarks on its application to Wasserstein-like gradient flows. It is remarkable from convergence analysis that the convergence rate of the proposed PDFB splitting method is independent of the spatial grid size, which is important for solving large-scale problems \cite{Jacobs2018}. We also introduce a nonlinear preconditioner according to our convergence analysis, which further accelerates the convergence of our PDFB method.

\section{Core formulation}
In this section, we present the core formulation in our method, and compare it with the previous framework to highlight our motivation and methodological novelty.
\subsection{Previous approach}
In previous approach \cite{Carrillo2023}, they considered a single-step approximation of the dynamical JKO scheme~\eqref{JKO} in which the numerical solution $\rho^{n + 1}$ is obtained as the optimal $\rho$ that minimizes the following problem:
\begin{equation}\label{c_scheme-old}
    \begin{aligned}
        &\min_{\rho, \textbf{m}} \Delta t\mathcal{E}(\rho) + \int_{\Omega} f\big(M(\rho), \textbf{m}\big)\,\text{d}x\\
        &\text{ s.t. }\rho - \rho^{n} + \nabla\cdot \textbf{m} = 0\text{ in } \Omega\,,\quad \textbf{m} \cdot \textbf{n} = 0 \text{ on } \partial \Omega\,.
    \end{aligned}
\end{equation}
\begin{remark}
It has been shown in \cite{Li2019, Carrillo2023, Li2023} that, in general, single-step discretizations for the dynamical formulation~\eqref{1-5}, which discretizes the integral of action function with respect to time with a right-hand side rule and the transport constraint with an implicit Euler method, do not violate the first-order accuracy of JKO schemes~\eqref{JKO}.
\end{remark}

Upon this single-step approximation, the optimization problem~\eqref{c_scheme-old} is handled by a Lagrangian multiplier $\phi$, in which the problem converted to a min-max saddle point problem given as follows:
\begin{equation}\label{osad}
    \min_{\rho, \textbf{m}} \max_{\phi} \Delta t\mathcal{E}(\rho) + \int_{\Omega} f\big(M(\rho), \textbf{m}\big)\,\text{d}x + \int_{\Omega} (\rho - \rho^{n}) \phi - \textbf{m} \cdot \nabla \phi \,\text{d}x\,.
\end{equation}
Let $u = (\rho, \textbf{m})$ be the primal variable, this problem can be then abbreviated as the following sum-of-three-operators form optimization problem:
\begin{equation}\label{osad-a}
    \min_{u} \max_{\phi} E(u) + \langle \mathcal{K}\phi, u \rangle - F(\phi) + G(u),
\end{equation}
where $E(u) = \Delta t \mathcal{E}(\rho_{1})$, $\mathcal{K}\phi = (\phi, -\nabla \phi)$,  and
\begin{equation*}
    F(\phi) = \int_{\Omega} \rho^{n}(x) \phi(x) \,\mathrm{d}x, \quad G(u) = \int_{\Omega} f\big(M(\rho), \textbf{m}\big)\,\text{d}x\,.
\end{equation*}
Here $\langle\cdot, \cdot\rangle$ denotes the sum of pair-wise $L^{2}$ inner products.

Such an optimization problem was considered in the primal-dual three splitting method (PD3O) method proposed by Yan \cite{Yan2016}, an extended version of primal-dual hybrid gradient (PDHG) method \cite{Chambolle2011, Pock2011, Chambolle2016} and the Davis-Yin three-operator splitting algorithm \cite{Davis2015}. In the previous approach~\cite{Carrillo2023}, we applied the PD3O method after suitable finite-dimensional approximation, in which the primal and dual variables in~\eqref{osad-a} are updated through the following iteration scheme:
\begin{equation}\label{PD3O}
    \begin{aligned}
        v^{(\ell + 1)} &= \text{prox}_{\sigma F}\big(v^{(\ell)} + \sigma \mathcal{K}\overline{u}^{(\ell)}\big),\\
        u^{(\ell + 1)} &=\text{prox}_{\tau G}\big(u^{(\ell)} -\tau \nabla_{u} E(u^{(\ell)})-\tau \mathcal{K}^{T}v^{(\ell+1)}\big),\\
        \overline{u}^{(\ell+1)}&= 2u^{(\ell+1)}-u^{(\ell)} -\tau (\nabla_{u} E(u^{(\ell+1)})-\nabla_{u} E(u^{(\ell)})).
    \end{aligned}
\end{equation}
This method converges when the step sizes satisfy $\sigma \tau \|\mathcal{K}\mathcal{K}^{T}\|_{\ast}\leq 1$ and $\tau\leq 2/C_E$ \cite{Yan2016}, where $C_E$ is the Lipschitz constant of $\nabla E(u)$.

A key drawback is that the operator $\mathcal{K}$ becomes ill-conditioned after discretization, which forces step sizes $(\sigma, \tau)$ to be extremely small on fine grids to satisfy the convergence criteria $\sigma \tau \|\mathcal{K}\mathcal{K}^{T}\|_{\ast}\leq 1$. In the previous work \cite{Carrillo2023}, we introduced a preconditioned version of the PD3O method to alleviate the severe restrictions on step sizes, yet the convergence of the algorithm remains grid-resolution-dependent. 

Another main technical difficulty in the implementation of PD3O lies in computing the proximal operator of $G$, which reduces to the separate subproblem for given $(\rho_{0}, \textbf{m}_{0}) \in \mathbb{R} \times \mathbb{R}^{d}$ at each grid if considering mid-point integration rule:
\begin{equation}\label{of}
    \min_{(\rho, \textbf{m}) \in \mathbb{R} \times \mathbb{R}^{d}} f\big(M(\rho), \textbf{m}\big) + \frac{1}{2\tau}|\rho - \rho_{0}|^{2} + \frac{1}{2\tau}|\textbf{m} - \textbf{m}_{0}|^{2}.
\end{equation}
For linear mobilities, problem~\eqref{of} admits a closed-form solution \citep[Proposition 1]{Papadakis2013}; while for quadratic concave mobilities, one can use an efficient Newton-based method with a careful criterion for initial values \citep[Section 3.3]{Carrillo2023} to obtain the minimizer in \eqref{of}. However, solving such minimization problems for general mobilities is \textit{non-trivial} due to the non-smoothness of $f$ shown in~\eqref{1-4}, not to mention the potential computational inefficiencies. This difficulty limits the applicability of the framework \eqref{osad}--\eqref{PD3O} to gradient flows with general mobilities, as well as the flexibility of designing higher-order discretizations.

\subsection{New formulation}
With the above mentioned limitations, one can see that adopting the framework in~\eqref{osad} to compute~\eqref{1-1} with general mobilities is difficult, as it requires developing a new algorithm to solve~\eqref{of} for each different $M(\rho)$. In this regard, we consider an alternative approach to address the action function $f$.

Our main novelty is to use Legendre–Fenchel transform of the action function $f$, introduced by Benamou and Brenier \cite{Benamou2000, Papadakis2013}:
\begin{equation}\label{3-12}
    f(M, \textbf{m}) = \max_{\phi, \psi} \Big\{M \phi + \textbf{m} \cdot \psi : \text{ s.t. } \phi + \tfrac{1}{2}|\psi|^{2} \le 0\Big\}\,,
\end{equation}
to circumvent the non-smoothness issue of $f$ and decouple the nonlinear relationship between $M$ and $\textbf{m}$ from $f(M, \textbf{m})$. Substituting \eqref{3-12} into~\eqref{JKO}, we can obtain an alternative saddle point formulation where the optimal $\rho$ provides the numerical solution $\rho^{n + 1}$:
\begin{equation}\label{csad}
    \begin{cases}
        \displaystyle\min_{\rho, \textbf{m}} \max_{\phi, \psi} \Delta t \mathcal{E}(\rho) + \int_{\Omega} M\big(\tfrac{1}{2}(\rho^{n} + \rho)\big)  \phi + \textbf{m} \cdot \psi \,\text{d}x\\
        \text{ s.t. } \rho - \rho^{n} + \nabla\cdot \textbf{m} = 0\text{ in } \Omega, \text{ and } \beta_{0} \le \rho \le \beta_{1} \text{ in } \Omega\,,\\
        \text{ and } \phi + \tfrac{1}{2}|\psi|^{2} \le 0 \text{ in } \Omega \times (0, 1)\,.
    \end{cases}
\end{equation}
\begin{remark}
This saddle-point formulation corresponds to an alternative single-step approximation of the dynamical JKO scheme~\eqref{JKO}. Combined with the staggered grid approximation given in Section~\ref{sec2} below, it yields a space-time staggered grid approximation of the dynamical JKO formulation:
\begin{equation}\label{c_scheme}
    \begin{aligned}
        &\min_{\rho, \textbf{m}} \Delta t\mathcal{E}(\rho) + \int_{\Omega} f\Big(M\big(\tfrac{1}{2}(\rho^{n} + \rho)\big), \textbf{m}\Big)\,\text{d}x\\
        &\text{ s.t. }\rho - \rho^{n} + \nabla\cdot \textbf{m} = 0\text{ in } \Omega,\\
        &\qquad\textbf{m} \cdot \textbf{n} = 0 \text{ on } \partial \Omega \text{ and } \beta_{0} \le \rho \le \beta_{1} \text{ in } \Omega\,.
    \end{aligned}
\end{equation}
Here we impose the global bound $[\beta_{0}, \beta_{1}]$ for $\rho$ such that $M(\rho)\geq 0$ to ensure bound-preservation at the discrete level. Introducing the box constraint in this formulation is essential given that the feasibility of $f$ only guarantees the bound of $\frac{1}{2}(\rho^{n} + \rho)$, rather than directly enforcing the bound on $\rho$.
\end{remark}

Let $u = (\rho, \textbf{m})$ and $v = (\phi, \psi)$ be the primal and dual variables respectively, the new formulation~\eqref{csad} can be written in the following form of general convex-concave saddle point problems:
\begin{equation*}\label{saddle1}
    \min_{u} \max_{v}\, \Phi(u, v) - F(v) + G(u)
\end{equation*}
with a nonlinear primal-dual coupling function $\Phi$ and two non-smooth convex functions $F$ and $G$ that are independent of specific forms of mobilities $M(\rho)$:
\begin{align*}
    &\Phi(u, v) = \Delta t \mathcal{E}(\rho) + \int_{\Omega} M\big(\tfrac{1}{2}(\rho^{n} + \rho)\big)  \phi + \textbf{m} \cdot \psi \,\text{d}x, \\
    &F(v) = \int_{\Omega} \delta_{\mathfrak{K}}(\phi, \psi)\,\text{d}x\,,\quad \text{and} \quad G(u) = \delta_{\mathfrak{D}}(u)\,,
\end{align*}
where $\mathfrak{K}$ and $\mathfrak{D}$ denote the constraint sets of the multipliers and the primal variables in~\eqref{csad}, respectively, and $\delta_{A}$ denotes the indicator function of the set $A$. 

Motivated by the primal-dual approaches for solving general convex-concave saddle point problems \cite{Malitsky2018, Hamedani2018, Zhu2022}, we propose a primal-dual forward-backward (PDFB) splitting method (Algorithm~\ref{algorithm1}), as a generalization of the PD3O method \cite{Yan2016}, for the new saddle-point formulation~\eqref{csad} with nonlinear primal-dual coupling. 

The new formulation, together with the proposed PDFB splitting method, accommodates general nonlinear mobilities by converting the challenging calculation of the proximal operator of non-smooth nonlinear function $f$ associated with general mobility $M(\rho)$ into the computation of proximal operators of $F$ and $G$, which are easily calculated projections onto convex sets independent of $M(\rho)$, and the explicit evaluation of the gradient of $\Phi$ that is straightforward for any smooth mobility. Furthermore, we establish the convergence analysis of the proposed PDFB splitting method, which is non-trivial due to the nonlinear dependence of the primal-dual coupling $\Phi(u,v)$ on $u$. In order to overcome this difficulty, we introduce the linearized gradient $\widehat{\nabla}_{v} \Phi(\overline{u}^{(\ell)}, v^{(\ell)}; u^{(\ell)})$ in the update of $v$ in Algorithm~\ref{algorithm1} to address the challenge. Remarkably, the proposed PDFB splitting method can achieve grid-resolution-independent convergence for one-step JKO problem, as demonstrated in the numerical experiments. 


\begin{algorithm}[t!]
\caption{PDFB splitting method}\label{algorithm1}
    \begin{algorithmic}[1]
        \State{\textbf{input} $\text{iter}_{\max}$, $\tau$, $\sigma$, $\text{tolerance}$}
        \State{\textbf{initialize} $u^{(0)} = 0$, $v^{(0)} = 0$, $\overline{u}^{(0)} = u^{(0)}$}
        \For{$\ell = 0 : \text{iter}_{\max}$}
            \State{Update dual variable
            \begin{equation*}
                \begin{aligned}
                    v^{(\ell + 1)} = \text{prox}_{\sigma F}\big(v^{(\ell)} &+ \sigma\widehat{\nabla}_{v}\Phi(\overline{u}^{(\ell)}, v^{(\ell)}; u^{(\ell)})\big)
                \end{aligned}
            \end{equation*}}
            \State{Update primal variable
            \begin{equation*}
                u^{(\ell + 1)} = \text{prox}_{\tau G}\big(u^{(\ell)} - \tau \nabla_{u}\Phi(u^{(\ell)}, v^{(\ell + 1)})\big)
            \end{equation*}}
            \State{Gradient-informed reflection
            \begin{equation*}
                \overline{u}^{(\ell + 1)} = 2u^{(\ell + 1)} - u^{(\ell)} - \tau \big(\nabla_{u}\Phi(u^{(\ell + 1)}, v^{(\ell + 1)}) - \nabla_{u}\Phi(u^{(\ell)}, v^{(\ell + 1)})\big)
            \end{equation*}
            }
            
            \State{Convergence criteria
            \begin{equation*}
                \frac{\|u^{(\ell + 1)} - u^{(\ell)}\|}{\|u^{(\ell + 1)}\|} \le \text{tolerance}.
            \end{equation*}}
        \EndFor
    \end{algorithmic}
\end{algorithm}
The rest of the paper is organized as follows. In Section~\ref{sec2}, we build a staggered grid approximation of the saddle point formulation~\eqref{csad} of one-step JKO problem. In Section~\ref{sec3}, we present and analyze the PDFB splitting method for solving a general saddle point problem corresponding to~\eqref{csad}. In Section~\ref{sec:app}, we discuss the implementation of the PDFB splitting method. In Section~\ref{sec:4}, we present the numerical results to demonstrate the performance of the proposed method. Finally, Section~\ref{sec5} provides our conclusion.

\section{Approximation on staggered grids}\label{sec2}
\subsection{Conforming approximation on staggered grids}
Abundant fully discrete approximations for min-max saddle problems have been considered in mixed or hybrid finite element methods~\cite{Boffi2013}. In this part, we construct a staggered grid approximation for~\eqref{csad}. It yields a lowest order Raviart-Thomas $H(\Omega \times (0, 1), \text{div})$-conforming element approximation of the dynamical JKO scheme~\eqref{JKO} on cubes. This conforming element discretization has been considered in \cite{Papadakis2013, Chambolle2019} for its desirable stability and conservative properties. Here to avoid complicated notations, we derived this approximation based on a finite volume description.

\begin{figure}[h]
    \centering
    \includegraphics[width = \textwidth]{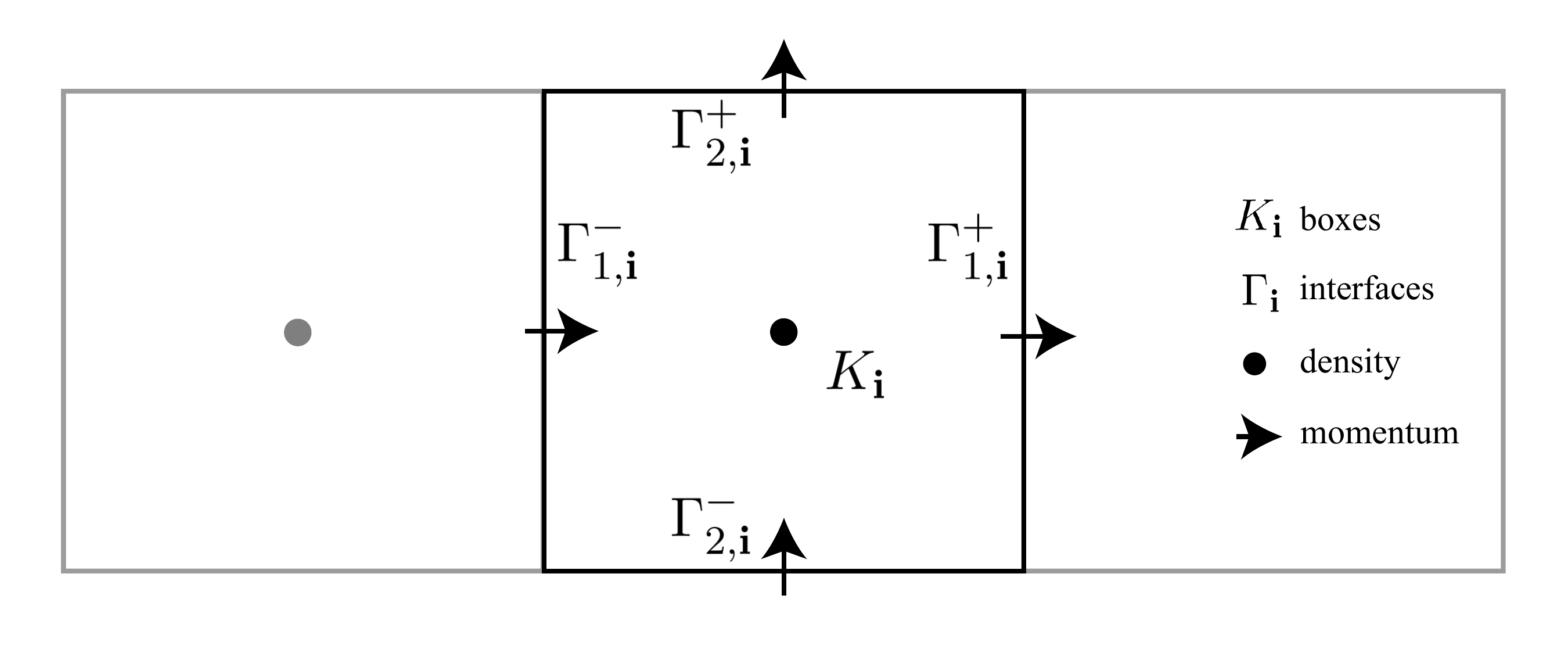}
    \caption{Illustration of the domain decomposition and the discrete variables. Here the dot denotes the location of discrete density in a box, and the arrows denotes the location of normal fluxes of discrete momentum on the bounding faces.}
    \label{fig:volume}
\end{figure}

For spatial discretization, as illustrated by Figure~\ref{fig:volume}, we cover the domain $\Omega\subset \mathbb{R}^{d}$ with non-overlapping uniform boxes $K_{\textbf{i}}$ defined by
\begin{equation*}
    K_{\textbf{i}} = \big\{x : \big|x - x_{\textbf{i}}\big|_{\infty} \le \tfrac{1}{2}h \big\},
\end{equation*}
where $h$ is the grid spacing and $x_{\textbf{i}}$ is the center of the box. The boundary of the box consists of two $(d-1)$-dimensional rectangles in each direction:
\begin{equation*}
    \partial K_{\textbf{i}} = \bigcup_{q = 1}^{d}\Gamma^{\pm}_{q, \textbf{i}}
\end{equation*}
In order to write the scheme in a finite-difference stencil, we adopt a global indices
\begin{equation*}
    \mathfrak{I} = \big\{\textbf{i} = \big(i_{1},..., i_{d}\big) \in \mathbb{N}^{d}\big\}
\end{equation*}
corresponding to the location of the volume $K_{\textbf{i}}$. Also, we introduce the fractional index
\begin{equation*}
    \textbf{i} \pm \tfrac{1}{2}e_{q} = \big(i_{1},..., i_{q} \pm \tfrac{1}{2}, ... i_{d}\big),
\end{equation*}
which corresponds to the location of the boundary interfaces $\Gamma_{q, \textbf{i}}^{\pm} \in \partial K_{\textbf{i}}$ orthogonal to the standard basis vector $e_{q}$ pointing in the $q$-th direction.

By integrating the discrete continuity equation constraint~\eqref{c_scheme} over an element $K_{\textbf{i}}$ and applying the divergence theorem, one has
\begin{equation}\label{dcont}
    \begin{aligned}
        \int_{K_{\boldsymbol{\mathrm{i}}}}\rho(x) - \rho^{n}(x) &+ \nabla \cdot \textbf{m}(x)\,\text{d}x\\
        &= \int_{K_{\boldsymbol{\mathrm{i}}}}\rho(x) - \rho^{n}(x) \,\text{d}x + \int_{\partial K_{\boldsymbol{\mathrm{i}}}} \textbf{n}_{K_{\boldsymbol{\mathrm{i}}}}(x) \cdot \textbf{m}(x) \,\text{d}s,
    \end{aligned}
\end{equation}
where $\textbf{n}_{K_{\boldsymbol{\mathrm{i}}}}$ is the outer unit normal to the volume $K_{\boldsymbol{\mathrm{i}}}$, and $\text{d}s$ is the surface measure on $\partial K_{\boldsymbol{\mathrm{i}}}$. Here, we compute the average over the volumes $K_{\boldsymbol{\mathrm{i}}}$ and interfaces $\partial K_{\boldsymbol{\mathrm{i}}}$ as piece-wise constant approximations of the variables. In particular, we let $\rho_{\textbf{i}}$ be the average over the $\textbf{i}$-th volume such that
\begin{equation*}
    \rho_{\textbf{i}} = \frac{1}{h^d}\int_{K_{\textbf{i}}} \rho(x) \text{d}x,\quad \rho_{\textbf{i}}^{n} = \frac{1}{h^d}\int_{K_{\textbf{i}}} \rho^{n}(x) \text{d}x,
\end{equation*}
and $m_{\textbf{i} \pm \frac{1}{2}e_{q}}^{q}$ be the average normal fluxes for the corresponding faces:
\begin{equation*}
    m_{\textbf{i} + \frac{1}{2}e_{q}}^{q}(x) = \frac{1}{h^{d- 1}}\int_{\Gamma_{q, \textbf{i}}^{+}} e_{q} \cdot \textbf{m}(x) \text{d}s,\quad m_{\textbf{i} - \frac{1}{2}e_{q}}^{q}(x) = \frac{1}{h^{d-1}}\int_{\Gamma_{q, \textbf{i}}^{-}} e_{q} \cdot \textbf{m}(x) \text{d}s.
\end{equation*}
By decomposing the boundary integral term in~\eqref{dcont} into the sum of average normal fluxes on interfaces $\partial K_{\textbf{i}} = \bigcup_{q = 1}^{d}\Gamma^{\pm}_{q, \textbf{i}}$, the constraints at the fully discrete level can be then written as follows,
\begin{equation}\label{constraint}
    \begin{aligned}
        &\rho_{\textbf{i}} - \rho_{\textbf{i}}^{n} + (\mathcal{A}\textbf{m})_{\textbf{i}} = 0,\quad \beta_{0} \le \rho_{\textbf{i}} \le \beta_{1},\text{ for all } \textbf{i} \in \mathfrak{I},
    \end{aligned}
\end{equation}
where $\mathcal{A}$ is the discrete divergence operator defined by
\begin{equation*}
    (\mathcal{A}\textbf{m})_{\textbf{i}} = \frac{1}{h}\sum_{q = 1}^{d} \Big(m_{\textbf{i} + \frac{1}{2}e_{q}}^{q} - m_{\textbf{i} - \frac{1}{2}e_{q}}^{q}\Big).
\end{equation*}
The no-flux boundary condition is satisfied by setting the average normal fluxes of the momentum field $\textbf{m}$ on $\partial \Omega$ to zero, which can be embedded into the design of the discrete matrix $\mathcal{A}$ (see Appendix~\ref{api}).
\subsection{Energy discretization}
We approximate the action functional at the discrete level as follows:
\begin{equation}\label{2-3}
    \int_{\Omega} M\big(\tfrac{1}{2}(\rho^{n} + \rho)\big) \phi + \textbf{m} \cdot \psi \,\text{d}x \approx h^{d}\sum_{\textbf{i} \in \mathfrak{I}} M\big(\tfrac{1}{2}(\rho_{\textbf{i}}^{n} + \rho_{\textbf{i}})\big) \phi_{\textbf{i}} + (\mathcal{I}\textbf{m})_{\textbf{i}} \cdot \psi_{\textbf{i}},
\end{equation}
where $\mathcal{I}$ is the averaging operator in space:
\begin{equation*}
     (\mathcal{I}\textbf{m})_{\textbf{i}} = \left(\frac{m_{\textbf{i} + \frac{1}{2}e_{1}}^{1} + m_{\textbf{i} - \frac{1}{2}e_{1}}^{1}}{2}, \cdots, \frac{m_{\textbf{i} + \frac{1}{2}e_{d}}^{d} + m_{\textbf{i} - \frac{1}{2}e_{d}}^{d}}{2}\right).
\end{equation*}
and subject the dual variables to the constraint $\phi_{\textbf{i}} + \frac{1}{2}|\psi_{\textbf{i}}|^{2} \le 0$ on each grid $\textbf{i}$.

On the other hand, the discretization of the free energy functional $\mathcal{E}(\rho)$ given in~\eqref{1-12} on uniform grids was considered in \cite{Carrillo2019, Carrillo2023}: For the internal and potential energies, a piecewise constant approximation using the midpoint rule can be used, and for the Dirichlet energy, we apply the central difference on staggered grids for $\nabla\rho$ and a trapezoidal rule for the integral. This trapezoidal rule leads to a centered difference approximation of the Laplacian that appears in the corresponding chemical potential. In addition, if Dirichlet energy is considered, a Neumann boundary condition on $\rho$ will be handled as described in \cite{Carrillo2023}. Further details will be provided in the numerical examples as needed. 

We let $\mathcal{E}_{h}$ to be the discrete free energy functional. Combining the discretizations above, we are able to get the following fully discrete min-max saddle point problem corresponding to:
\begin{equation}\label{2-8}
    \begin{aligned}
        \min_{\rho, \textbf{m}} \max_{\phi, \psi} \Delta t \widehat{\mathcal{E}}_{h}(\rho) + \sum_{\textbf{i} \in \mathfrak{I}} M\big(\tfrac{1}{2}(\rho^{n}_{\textbf{i}} + \rho_{\textbf{i}})\big) \phi_{\textbf{i}} &+ (\mathcal{I}\textbf{m})_{\textbf{i}} \cdot \psi_{\textbf{i}}\\
        &- \sum_{\textbf{i} \in \mathfrak{I}}\delta_{\mathfrak{K}}(\phi_{\textbf{i}}, \psi_{\textbf{i}}) + \delta_{\mathfrak{D}[h]}(\rho, \textbf{m}),
    \end{aligned}
\end{equation}
where $\widehat{\mathcal{E}}_{h} := \mathcal{E}_{h}/h^{d}$ is the rescaled discrete energy,  and $\delta_{\mathfrak{J}}$ is the indicator function of a convex set $\mathfrak{J}$ defined as follows
\begin{equation*}
    \begin{cases}
        \delta_{\mathfrak{J}}(v) = 0 &\text{ if } v \in \mathfrak{J},\\
        \delta_{\mathfrak{J}}(v) = +\infty &\text{ otherwise. }
    \end{cases}
\end{equation*}
We have introduced the penalty terms of $\delta_{\mathfrak{D}[h]}, \delta_{\mathfrak{K}}$ to remove the constraints on the primal and dual variables, respectively, where $\mathfrak{D}[h]$ and $\mathfrak{K}$ are the discrete admissible sets defined by
\begin{equation*}
    \begin{aligned}
    \mathfrak{D}[h] &= \Big\{(\rho, \textbf{m}) : \rho_{\textbf{i}} - \rho^{n}_{\textbf{i}} + (\mathcal{A}\textbf{m})_{\textbf{i}} = 0,\quad \beta_{0} \le \rho_{\textbf{i}} \le \beta_{1},\quad \textbf{i} \in \mathfrak{I}\Big\},\\
    \mathfrak{K} &= \Big\{(\phi, \psi) \in \mathbb{R} \times \mathbb{R}^{d} : \phi + \tfrac{1}{2}|\psi|^{2} \le 0\Big\}.
    \end{aligned}
\end{equation*}
This constructed staggered grid approximation has the ability to preserve the structural properties of the Wasserstein-like gradient flow~\eqref{1-1} at the discrete level given as follows.

\begin{theorem}\label{ts}Given a suitable initial value $\rho_{0}$, if there exists a time-sequence of solutions $(\rho^{n})_{n \ge 1}$ generated by the fully discrete scheme~\eqref{2-8}, we have
\begin{enumerate}
    \item energy dissipation $\mathcal{E}_{h}(\rho^{n + 1}) \le \mathcal{E}_{h}(\rho^{n})$.
    \item mass conservation $\sum_{\textnormal{\textbf{i}} \in \mathfrak{I}} \rho_{\textnormal{\textbf{i}}}^{n + 1} = \sum_{\textnormal{\textbf{i}} \in \mathfrak{I}} \rho_{\textnormal{\textbf{i}}}^{n}$.
    \item bound preservation $M(\rho_{\textnormal{\textbf{i}}}^{n + 1}) \ge 0,~ \text{ for all } \textnormal{\textbf{i}} \in \mathfrak{I}$.
\end{enumerate}
\end{theorem}
\begin{proof}
The proof is given in Appendix~\ref{ast}.
\end{proof}
\begin{remark}
The existence of minimizers for continuous JKO schemes with nonlinear concave mobilities and appropriate free energies has been established in \cite{Lisini2012}. Such existence result can be extended to the fully discrete scheme by following similar proof routine for linear mobility case in \cite{Carrillo2022}. Since the main focus of this paper is on the efficient numerical realization of JKO schemes, we leave the rigorous analysis of well-posedness of JKO scheme itself for more general cases in future work.
\end{remark}

\section{PDFB splitting method}\label{sec3}
In this section, we provide a detailed derivation of the PDFB splitting method to solve the discrete saddle point problem~\eqref{2-8}. For notational simplicity, we denote the primal variables by $u$ and the dual variables by $v$, and write the discrete saddle point~\eqref{2-8} as follows:
\begin{equation}\label{abbsad}
    \min_{u = (\rho, \textbf{m})} \max_{v = (\phi, \psi)} \Phi(u, v) - \sum_{\textbf{i} \in \mathfrak{I}} \delta_{\mathfrak{K}}(v_{\textbf{i}}) + \delta_{\mathfrak{D}[h]}(u),
\end{equation}
where $v_{\textbf{i}} = (\phi_{\textbf{i}}, \psi_{\textbf{i}})$ denotes the dual variables on the grid $K_{\textbf{i}}$, and 
\begin{equation}\label{Phi}
    \Phi(u, v) = \Delta t \widehat{\mathcal{E}}_{h}(\rho) + \sum_{\textbf{i} \in \mathfrak{I}} M\big(\tfrac{1}{2}(\rho^{n}_{\textbf{i}} + \rho_{\textbf{i}})\big) \phi_{\textbf{i}} + (\mathcal{I}\textbf{m})_{\textbf{i}} \cdot \psi_{\textbf{i}}
\end{equation}
is the discrete objective functional.

Consider the following general form of problem~\eqref{abbsad} as a convex-concave saddle point problem \cite{Hamedani2018}:
\begin{equation}\label{3-20}
    \min_{u}\max_{v} \mathscr{L}(u, v) := \Phi(u, v) - F(v) + G(u),
\end{equation}
where $F$ and $G$ are convex (but possibly non-smooth), and $\Phi(u, v)$ is continuously differentiable. We assume that $\Phi(u, v)$ is affine with respect to the dual variable $v$, which is satisfied by our problem~\eqref{Phi}. Additionally, we assume that problem~\eqref{3-20} admits a solution $(u^{\ast}, v^{\ast})$ characterized by the following variational inequality
\begin{equation*}
    \begin{aligned}
        0 &\in \nabla_{v}\Phi(u^{\ast}, v^{\ast}) - \partial F(v^{\ast}) \,;\\
        0 &\in \nabla_{u}\Phi(u^{\ast}, v^{\ast}) + \partial G(u^{\ast})\,,
    \end{aligned}
\end{equation*}
where $\nabla_{u}$ and $\nabla_{v}$ are the gradient operators with respect to $u$ and $v$ respectively, and $\partial F$ denotes the subgradient of $F$. We also assume that $F$ and $G$ are ``simple" in the sense that their proximal operators can be computed either in closed form or with high accuracy by a few iterations of some subroutine \cite{Chambolle2011}. In the following, we will use $\langle u, v\rangle = u^{T}v$ to denote the inner product of two vectors of the same dimension, and we will frequently use $\text{prox}_{F}$, the proximal operator of $F$, defined as
\begin{equation*}
    \text{prox}_{F}(v_{0}) = \arg\min_{v} F(v) + \frac{1}{2}\|v - v_{0}\|^{2},
\end{equation*}
where $\|\cdot\| = \langle \cdot, \cdot\rangle^{\frac{1}{2}}$ is the Euclidean norm. We will also use $\|\cdot\|_{\ast}$ to denote the spectral norm of a matrix, defined as its largest eigenvalue.

\subsection{Algorithm}
As demonstrated in Algorithm~\ref{algorithm1}, our PDFB splitting method is a variant of the operator splitting optimization methods proposed in \cite{Malitsky2018, Hamedani2018, Zhu2022} and an extension of the PD3O method \cite{Yan2016}. Unlike PD3O which considers a functional $\Phi$ of the form~\eqref{osad-a}, we consider a general form as given in~\eqref{Phi}, where the linear term $\mathcal{K}u$ is replaced by a nonlinear function of $u$. Compared to standard forward-backward splitting \cite{Chen1997, Malitsky2018},
\begin{equation*}
    \begin{aligned}
        v^{(\ell + 1)} &= \text{prox}_{\sigma F}\big(v^{(\ell)} + \sigma\nabla_{v}\Phi(u^{(\ell)}, v^{(\ell)})\big)\,;\\
        u^{(\ell + 1)} &= \text{prox}_{\tau G}\big(u^{(\ell)} - \tau \nabla_{u}\Phi(u^{(\ell)}, v^{(\ell)})\big)\,;\\
    \end{aligned}
\end{equation*}
we incorporate elements from the primal-dual hybrid gradient (PDHG) method \cite{Chambolle2011, Pock2011, Chambolle2016} and PD3O, such as using different step sizes for handling $F$ and $G$ and introducing a reflection step, which improve convergence behavior. 

More specifically, following the update rule for the dual variable $v$, we have:
\begin{equation*}
    v^{(\ell + 1)} \approx \text{prox}_{\sigma F}\big(v^{(\ell)} + \sigma\nabla_{v}\Phi(\overline{u}^{(\ell)}, v^{(\ell)})\big),
\end{equation*}
where $\nabla_v \Phi(u,v):  \mathbb{R}^{\text{dim}(u)} \to \mathbb{R}^{\text{dim}(v)}$ is a vector-valued function that depends only on $u$, as $\Phi(u,v)$ is affine in $v$. However, $\overline{u}^{(\ell)}$ may not lie within the domain of $\nabla_{v}\Phi$. To facilitate a rigorous convergence analysis, we approximate the gradient using a first-order expansion around  $u^{(\ell)}$:
\begin{equation}\label{eq:Phi_v_expansion}
    \widehat{\nabla}_{v}\Phi(\overline{u}^{(\ell)}, v^{(\ell)}; u^{(\ell)}) := \mathcal{K}(u^{(\ell)}) (\overline{u}^{(\ell)} - u^{(\ell)}) + \nabla_{v} \Phi(u^{(\ell)}, v^{(\ell)}).
\end{equation}
Here, $\mathcal{K}(u)$ is the Jacobi matrix defined by
\begin{equation}\label{k}
    \begin{aligned}
        \mathcal{K}(u) &: \mathbb{R}^{\text{dim}(u)} \to \mathbb{R}^{\text{dim}(v) \times \text{dim}(u)}\\
        &\text{ such that } \mathcal{K}(u)w = \lim_{\varepsilon \to 0} \frac{\nabla_{v} \Phi(u + \varepsilon w, v) - \nabla_{v} \Phi(u, v)}{\varepsilon},\quad \forall w \in \mathbb{R}^{\text{dim}(u)}.
    \end{aligned}
\end{equation}
Consequently, we update the dual variable as follows:
\begin{equation*}
    v^{(\ell + 1)} = \text{prox}_{\sigma F}\big(v^{(\ell)} + \sigma\widehat{\nabla}_{v}\Phi(\overline{u}^{(\ell)}, v^{(\ell)}; u^{(\ell)})\big)\,.
\end{equation*}
The update for the primal variable follows closely from \eqref{PD3O}, and our final {\bf PDFB} splitting iterative procedure is: 
\begin{equation}\label{PDFB}
    \begin{aligned}
        v^{(\ell + 1)} &= \text{prox}_{\sigma F}\big(v^{(\ell)} + \sigma\widehat{\nabla}_{v}\Phi(\overline{u}^{(\ell)}, v^{(\ell)}; u^{(\ell)})\big)\,;\\
        u^{(\ell + 1)} &= \text{prox}_{\tau G}\big(u^{(\ell)} - \tau \nabla_{u}\Phi(u^{(\ell)}, v^{(\ell+1)})\big)\,;\\
        \overline{u}^{(\ell+1)} &= 2u^{(\ell+1)} - u^{(\ell)} - \tau \big(\nabla_{u}\Phi(u^{(\ell+1)}, v^{(\ell+1)}) - \nabla_{u}\Phi(u^{(\ell)}, v^{(\ell+1)})\big)\,.
    \end{aligned}
\end{equation}
We note that when $\Phi(u,v)=E(u)+\langle \mathcal{K}u, v\rangle$, our proposed procedure recovers the PD3O method. This is because, for this specific form of $\Phi$, we have:
\begin{equation*}
    \widehat{\nabla}_{v}\Phi(\overline{u}^{(\ell)}, v^{(\ell)}; u^{(\ell)})=\mathcal{K}\overline{u}^{(\ell)} \text{ and } \nabla_{u}\Phi(u^{(\ell)}, v^{(\ell + 1)})= \nabla_{u}E(u^{(\ell)}) + \mathcal{K}^{T}v^{(\ell + 1)}.
\end{equation*}
The last step is a “reflection” technique, providing an approximation $\overline{u}^{(\ell)} \approx u^{\ell + 1}$ of a fully implicit proximal gradient scheme, and is the critical part guaranteeing the convergence of the algorithm.

\subsection{Convergence analysis}
To analyze the convergence of the proposed PDFB splitting method, we rewrite the iteration~\eqref{PDFB} in the following form analogous to the Davis-Yin algorithm \cite{Davis2015}:
\begin{equation}\label{davis-yin}
    \begin{aligned}
        u^{(\ell)} &= \text{prox}_{\tau G}(q^{(\ell)})\\
        v^{(\ell + 1)} &= \text{prox}_{\sigma F}\big(v^{(\ell)} + \sigma \mathcal{K}(u^{(\ell)}) \overline{u}^{(\ell)} + \sigma L(u^{(\ell)})\big)\\
        q^{(\ell + 1)} &= u^{(\ell)} - \tau \nabla_{u}\Phi(u^{(\ell)}, v^{(\ell + 1)})
    \end{aligned}
\end{equation}
where $L(u^{(\ell)}) = \nabla_{v}\Phi(u^{(\ell)}, v^{(\ell)}) - \mathcal{K}(u^{(\ell)})u^{(\ell)}$. 

To facilitate the convergence analysis, we make the following assumptions: 
\begin{assumption}[Cocoercivity]\label{ccs}\rm
$G$ and $F$ are proper lower semi-continuous convex functions, and in addition $F$ is $\beta$-strongly convex, i.e., there exists $\beta>0$ such that
\begin{equation*}\label{4-14}
\begin{aligned}
    F(v) &\ge F(\overline{v}) + \langle p, v - \overline{v}\rangle + \frac{\beta}{2}\|v - \overline{v}\|^{2}\text{ for all $v, \overline{v} \in \textnormal{dom}(F)$ and $p \in \partial F(\overline{v})$.}
\end{aligned}
\end{equation*}
\end{assumption}
\begin{assumption}[Lipschitz continuity \cite{Hamedani2018, Zhu2022}]\label{ls}\rm
$\Phi(u, v)$ is continuously differentiable, convex in $u$, and affine in $v$. In addition, $\nabla_{u}\Phi(u, v)$ is Lipschitz continuous with respect to $u$. Specifically, for all $v \in \textnormal{dom}(F)$, we have
\begin{equation*}
    \big\|
    \nabla_{u}\Phi(u_{0}, v) - \nabla_{u}\Phi(u_{1}, v)\big\| \le C_{\Phi}\|u_{0} - u_{1}\|,\quad \forall u_{0}, u_{1} \in \textnormal{dom}(G),
\end{equation*}
with some constant $C_{\Phi}$.
\end{assumption}
\begin{assumption}[Continuity]\rm
$\mathcal{K}(u)$ is bounded such that, for all $u \in \textnormal{dom}(G)$,
\begin{equation*}
    \|\mathcal{K}(u)^{T}v\| \le C_{K}\|v\|,\quad \forall v \in \textnormal{dom}(F).
\end{equation*}
with some constant $C_{K}$.
\end{assumption}

\begin{theorem}\label{convergence}
Letting the sequence $\{(u^{(\ell)}, v^{(\ell)})\}$ be generated by Algorithm~\ref{algorithm1} and the above assumptions be satisfied, then the following ergodic sequence
\begin{equation*}
    u_{(N)} = \sum_{\ell = 0}^{N}u^{(\ell)}/N \quad  \text{and}\quad v_{(N)} = \sum_{\ell = 0}^{N}v^{(\ell + 1)}/N
\end{equation*}
converges to the saddle point $(u^{\ast}, v^{\ast})$ of the problem~\eqref{3-20},  if the step sizes $\tau$ and $\sigma$ satisfy
\begin{equation}\label{convergence_c}
    \tau \sigma \le \min_{u \in \textnormal{dom}(G)}\frac{1}{\|\mathcal{K}(u)\mathcal{K}(u)^{T}\|_{\ast}} ~~\text{ and }~~ \tau \le \frac{1}{C_{\Phi} + C_{0}},
\end{equation}
where $C_{0}$ is a positive constant that depends on $C_{\Phi}$, $C_{K}$, and $\beta$.

In addition, we have the following estimate of the primal-dual gap:
\begin{equation*}
    \mathscr{L}(u_{(N)}, v) - \mathscr{L}(u, v_{(N)}) \le \frac{\|q^{(0)} - q\|^{2}}{2N\tau} + \frac{\langle v^{(0)} - v, \mathcal{T}(u^{(0)})(v^{(0)} - v)\rangle}{2N\sigma}
\end{equation*}
where $q = u - \tau \nabla_{u}\Phi(u, v)$ and $\mathcal{T}(u) = I - \tau \sigma \mathcal{K}(u)\mathcal{K}(u)^{T}$.
\end{theorem}

\begin{remark}
    The proof is given in Section~\ref{appendix}, following the idea of Yan's proof of the convergence of PD3O method \cite{Yan2016}. Compared to his proof of convergence, the main technical difficulty in the analysis of the proposed PDFB splitting method comes from the nonlinearity of $\mathcal{K}(u)$. In this regard, we introduced the linearization of $\nabla_v \Phi(u,v)$ in \eqref{eq:Phi_v_expansion} to make analogy with the Davis-Yin iteration form \eqref{davis-yin}, and assume the strong convexity of $F$ and the boundedness of $\mathcal{K}(u)$ to derive the estimate of the primal-dual gap.
\end{remark}

\section{Application to Wasserstein-like gradient flows}\label{sec:app}
In this section, we provide implementation details and several remarks on the PDFB splitting method when applied to computing Wasserstein-like gradient flows formulated in~\eqref{2-8}. First, recall that:
\begin{equation*}
    \Phi(u, v) = \Delta t \widehat{\mathcal{E}}_{h}(\rho) + \sum_{\textbf{i} \in \mathfrak{I}} M\big(\tfrac{1}{2}(\rho^{n}_{\textbf{i}} + \rho_{\textbf{i}})\big) \phi_{\textbf{i}} + (\mathcal{I}\textbf{m})_{\textbf{i}} \cdot \psi_{\textbf{i}},
\end{equation*}
\begin{equation*}
    F(v) = \sum_{\textbf{i} \in \mathfrak{I}} \delta_{\mathfrak{K}}(\phi_{\textbf{i}}, \psi_{\textbf{i}}) \text{ and } G(u) = \delta_{\mathfrak{D}}(u),
\end{equation*}
where $u = (\rho, \textbf{m})$ and $v = (\phi, \psi)$.

\subsection{Implementation}
In this part, we sketch the implementation of the PDFB splitting method in Algorithm~\ref{algorithm1} by writing down the gradients $\nabla_{u}\Phi, \nabla_{v}\Phi$ and the Jacobi matrix $\mathcal{K}(u)$ explicitly.

In particular, given the $\ell$-th iterate $v^{(\ell)} = (\phi^{(\ell)}, \psi^{(\ell)})$ and $u^{(\ell)} = (\rho^{(\ell)}, \textbf{m}^{(\ell)})$ and $\overline{u}^{(\ell)} = (\overline{\rho}^{(\ell)}, \overline{\textbf{m}}^{(\ell)})$, we spell out  the Algorithm~\ref{algorithm1} as follows:
\begin{itemize}
    \item Conduct a gradient ascent step on the dual variables:
    \begin{equation*}
        \begin{aligned}
            \phi^{(\ell + \frac{1}{2})} = \phi^{(\ell)} + \sigma \Big[\tfrac{1}{2} M^{\prime}\big(\tfrac{1}{2}(\rho^{n} &+ \rho^{(\ell)})\big) (\overline{\rho}^{(\ell)} - \rho^{(\ell)}) + M\big(\tfrac{1}{2}(\rho^{n} + \rho^{(\ell)})\big)\Big],\\
            \psi^{(\ell + \frac{1}{2})} &= \psi^{(\ell)} + \sigma \mathcal{I}\overline{\textbf{m}}^{(\ell)}.
        \end{aligned}
    \end{equation*}
    \item Project the dual variables onto $\mathfrak{K}$ by a proximal gradient step:
    \begin{equation*}
        \big(\phi^{(\ell + 1)}, \psi^{(\ell + 1)}\big) = \text{prox}_{\sigma F}\big(\phi^{(\ell + \frac{1}{2})}, \psi^{(\ell + \frac{1}{2})}\big).
    \end{equation*}
    \item Update the primal variable by a gradient descent step: 
    \begin{equation*}
        \begin{aligned}
            \rho^{(\ell + \frac{1}{2})} = \rho^{(\ell)} - \tau \bigg[\Delta t \frac{\delta \widehat{\mathcal{E}}_{h}(\rho^{(\ell)})}{\delta \rho} &+ \tfrac{1}{2}M^{\prime}\big(\tfrac{1}{2}(\rho^{n} + \rho^{(\ell)})\big)\phi^{(\ell + 1)}\bigg],\\
            \textbf{m}^{(\ell + \frac{1}{2})} = \textbf{m}^{(\ell)} &- \tau \mathcal{I}^{T} \psi^{(\ell + 1)},
        \end{aligned}
    \end{equation*}
    where $\delta \widehat{\mathcal{E}}_{h}(\rho)/{\delta \rho}$ is the gradient of the discrete free energy functional $\widehat{\mathcal{E}}_{h}(\rho)$. 
    \item Project the primal variable onto $\mathfrak{D}$ by computing the proximal operator:
    \begin{equation*}
        \big(\rho^{(\ell + 1)}, \textbf{m}^{(\ell + 1)}\big) = \text{prox}_{\tau G}\big(\rho^{(\ell + \frac{1}{2})}, \textbf{m}^{(\ell + \frac{1}{2})}\big).
    \end{equation*}
    \item Update reflection estimate by:
    \begin{equation*}
        \begin{aligned}
            \overline{\rho}^{(\ell + 1)} &= 2\rho^{(\ell + 1)} - \rho^{(\ell)} - \tau \bigg[\Delta t \frac{\delta \widehat{\mathcal{E}}_{h}(\rho^{(\ell + 1)})}{\delta \rho} - \Delta t \frac{\delta \widehat{\mathcal{E}}_{h}(\rho^{(\ell)})}{\delta \rho}\\
            &\quad + \tfrac{1}{2}\Big[M^{\prime}\big(\tfrac{1}{2}(\rho^{n} + \rho^{(\ell + 1)})\big) - M^{\prime}\big(\tfrac{1}{2}(\rho^{n} + \rho^{(\ell)})\big)\Big]\phi^{(\ell + 1)}\bigg],
        \end{aligned}
    \end{equation*}
    \begin{equation*}
        \overline{\textbf{m}}^{(\ell + 1)} = 2\textbf{m}^{(\ell + 1)} - \textbf{m}^{(\ell)}.
    \end{equation*}
\end{itemize}
The success of the above procedure depends on the fact that proximal operators of $F$ and $G$ in the variational scheme~\eqref{2-8} can be computed either through closed-form expressions or with just a few iterations. Details will be provided in Section~\ref{sec:3-3}.

\begin{remark}\label{writeeng}
In this iterate, we only need to evaluate a discrete approximation of the ``effective" chemical potential $\delta \widehat{\mathcal{E}}_{h}(\rho)/{\delta \rho}$ in Algorithm~\ref{algorithm1}, whose explicit expression is always available as $\delta \mathcal{E}(\rho)/\delta \rho$ in~\eqref{1-1}. Hence, the PDFB splitting algorithm remains valid in practice even when the explicit expression of the free energy functional is unknown; see Example~\ref{non-e}.
\end{remark}

\subsection{Remarks on the convergence result}
In this part, we make several remarks regarding the convergence of the algorithm when applied to Wasserstein-like gradient flows.

\subsubsection{Strongly convex assumption}\label{remark7}

In the saddle point formulation~\eqref{2-8}, $F = \sum \delta_{\mathfrak{K}}$ is not strongly convex, thus failing to satisfy the coercivity condition in Assumption~\ref{ccs}. Recently, Zhu, Liu and Tran-Dinh~\cite{Zhu2022} proposed a primal-dual method that does not require the coercivity condition to prove convergence. However, this method involves a relatively complicated construction. 

Within the current PDFB splitting framework, one possible remedy is to apply certain regularization to $F$ to enforce strong convexity, such as the Moreau–Yosida regularization \cite{Clason2016}. Such regularization, however, inevitably introduces additional deviation in solutions. Alternatively, one might relax the coercivity condition by using a semi-implicit treatment of the mobility in the approximation of the induced transport distance; similar discretization have been directly applied to PDEs for Wasserstein gradient flows \cite{Chen2019CH-log,Hu2020SSIPNP,Shen2021PNP}. For example, for $M(\rho) = \rho(1 - \rho)$, one can approximate the mobility as $M(\rho) \approx \rho(1 - \rho^{n})$, which yields a linear function of current density. For such linear mobilities, the PDFB splitting method reduces to the PD3O method \cite{Yan2016}, and the resulting discrete optimization problem falls within the established convergence framework of the primal--dual splitting algorithm. However, this semi-implicit mobility scheme brings an extra approximation error and its convergence to the solution of continuous PDE requires a separate analysis.

We emphasize that the strong convexity assumptions in our convergence analysis are sufficient but may not be necessary for convergence. In practice, we find that the PDFB splitting method still converges reliably in the absence of strong convexity and produces more accurate solutions with the implicit mobility scheme. In contrast, while the Moreau–Yosida regularization and the semi-implicit treatment of the mobility guarantee the convergence of the primal-dual splitting algorithm, they also introduce extra regularization bias or approximation errors that reduce solution accuracy (see Examples~\ref{examp-2} and~\ref{examp-3} in Section~\ref{sec:4}).




\subsubsection{Condition on step size}
To verify the convergence requirement on the step sizes in~\eqref{convergence_c} when applied to~\eqref{2-8}, a direct computation based on~\eqref{Phi} yields: 
\begin{equation*}
    \begin{aligned}
        \big\|\mathcal{K}(u)\mathcal{K}(u)^{T}\big\|_{\ast} &= \max_{\textbf{i}}\Big\{\|\mathcal{I}\mathcal{I}^{T}\|_{\ast}, \tfrac{1}{4} M^{\prime}\big(\tfrac{1}{2}(\rho_{\textbf{i}} + \rho_{\textbf{i}}^{n})\big)^{2}\Big\}\\
        &\le \max_{\textbf{i}}\Big\{1, \tfrac{1}{4} M^{\prime}\big(\tfrac{1}{2}(\rho_{\textbf{i}} + \rho_{\textbf{i}}^{n})\big)^{2}\Big\}.
    \end{aligned}
\end{equation*}
where $u = (\rho, \textbf{m})$, which leads to the restriction on step sizes:
\begin{equation*}
    \tau \sigma \le \min_{\rho}\frac{1}{\max_{\textbf{i}} \Big\{1, \tfrac{1}{4} M^{\prime}\big(\tfrac{1}{2}(\rho_{\textbf{i}} + \rho_{\textbf{i}}^{n})\big)^{2}\Big\}}.
\end{equation*}
Here the derivative of the mobility function is typically bounded as long as $\rho_{\textbf{i}}$ remains bounded, which holds in all our numerical examples. More importantly, this criterion is independent of the grid size $h$, which suggests a $h$-independent convergence rate, provided that $\widehat{\mathcal{E}}_{h}$ has a Lipschitz gradient independent of $h$. A numerical demonstration of the rate of convergence with respect to $h$ is considered in Example~\ref{examp1} in Section~\ref{sec:4}.

\begin{remark}[Comparison with the previous framework]
For the previous framework \cite{Carrillo2019, Carrillo2023}, when applying the gradient descent algorithms (e.g., the PD3O algorithm \cite{Yan2016}) to~\eqref{osad}, the step sizes and convergence rates can be severely limited by the condition number of the discrete Laplace operator (i.e., $\|\mathcal{K}\mathcal{K}^{T}\|_{\ast}$), which increases significantly as the grid size reduces. In this regard, one needs preconditioning techniques \cite{Pock2011,Liu2021AccelerationPrimalDual} to accelerate the primal-dual methods \cite{Carrillo2023}. Implementing such a preconditioning technique for non-smooth problems can be non-trivial, which might limit its application to large-scale simulations \cite{Jacobs2018}.
\end{remark}

From the convergence condition in~\eqref{convergence_c}, the Lipschitz constant $C_{\Phi}$ of $\nabla_{u}\Phi(u, v)$ in Assumption~\ref{ls} also imposes a restriction on the step size $\tau$:
\begin{equation} \label{0402}
    \tau \le (C_{\Phi} + C_{0})^{-1} \le C_{\Phi}^{-1}.
\end{equation}
In practice, this restriction on $\tau$ is mainly determined by the Lipschitz continuity of the gradient of $\widehat{\mathcal{E}}_{h}$. However, some energy functionals may be “ill-conditioned", meaning that their corresponding Lipschitz constant can be extremely large. Consequently, the step size $\tau$ can be severely limited when solving these models. For example, consider the Gibbs-Boltzmann entropy functional:
\begin{equation*}
    \mathcal{E}(\rho) = \int_{\Omega} \rho(x)\log(\rho(x)) - \rho(x)\text{d}x \text{ where } \frac{\delta \mathcal{E}(\rho)}{\delta \rho}(x) = \log(\rho(x)).
\end{equation*}
In this example, the step size $\tau$ needs to be sufficiently small to satisfy \eqref{0402} when $\rho \approx 0$ at some grid. Certain preconditioning techniques can be utilized to address ill-conditioned free energy functionals, which are detailed in Section~\ref{sec:3-4}.

\subsection{Convex splitting as a nonlinear preconditioning}\label{sec:3-4}
In this part, we introduce a convex splitting approach for the variational problem~\eqref{2-8} when energy functional $\widehat{\mathcal{E}}_{h}$ is ill-conditioned. This approach can be served as a nonlinear preconditioning to our PDFB splitting method introduced above. Specifically, we split $\widehat{\mathcal{E}}_{h}$ into two parts
\begin{equation}\label{4-10}
    \widehat{\mathcal{E}}_{h}(\rho) = \mathcal{U}(\rho) + \mathcal{V}(\rho),
\end{equation}
where $\mathcal{U}$ is ill-conditioned but convex, and the Lipschitz constant of the gradient of $\mathcal{V}$ is uniformly bounded. The main idea is to treat the ill-conditioned part $\mathcal{U}$ with a backward gradient step, while addressing $\mathcal{V}$ with the forward gradient step.

By applying the convex conjugate to $\mathcal{U}$, problem~\eqref{2-8} can be reformulated as:
\begin{equation}\label{3-15}
    \begin{aligned}
        \min_{\rho, \textbf{m}} \max_{\phi, \psi, \mu} \Delta t \mathcal{V}(\rho) + \sum_{\textbf{i} \in \mathfrak{I}} M\big(\tfrac{1}{2}(\rho_{\textbf{i}} &+ \rho^{n}_{\textbf{i}})\big) \phi_{\textbf{i}} + (\mathcal{I}\textbf{m})_{\textbf{i}} \cdot \psi_{\textbf{i}} + \rho_{\textbf{i}} \mu_{\textbf{i}}\\
        &- \delta_{\mathfrak{K}}(\phi_{\textbf{i}}, \psi_{\textbf{i}}) - \Delta t \mathcal{U}^{\ast}(\mu/\Delta t) + \delta_{\mathfrak{D}}(\rho, \textbf{m}),
    \end{aligned}
\end{equation}
where $\mathcal{U}^{\ast}$ is the convex conjugate of $\mathcal{U}$ and $\mu$ is the corresponding dual variable introduced by taking the convex conjugate. Abbreviating the primal variable and the dual variable as $u$ and $v$ respectively, we can rewrite~\eqref{3-15} in a compact form:
\begin{equation*}
    \max_{u = (\rho, \textbf{m})}\min_{v = (\phi, \psi, \mu)} \Phi(u, v) - F(v) + G(u),
\end{equation*}
where
\begin{equation*}
    \Phi(u, v) = \Delta t \mathcal{V}(\rho) + \sum_{\textbf{i} \in \mathfrak{I}} M\big(\tfrac{1}{2}(\rho^{n}_{\textbf{i}} + \rho_{\textbf{i}})\big) \phi_{\textbf{i}} + (\mathcal{I}\textbf{m})_{\textbf{i}} \cdot \psi_{\textbf{i}} + \rho_{\textbf{i}} \mu_{\textbf{i}},
\end{equation*}
\begin{equation}\label{3-25}
    F(v) = \sum_{\textbf{i} \in \mathfrak{I}} \delta_{\mathfrak{K}}(\phi_{\textbf{i}}, \psi_{\textbf{i}}) + \Delta t \mathcal{U}^{\ast}(\mu/\Delta t) \text{ and } G(u) = \delta_{\mathfrak{D}}(u),
\end{equation}
which can be solved by the PDFB splitting method (Algorithm~\ref{algorithm1}). By separating the ill-conditioned part $\mathcal{U}$ from the objective functional $\Phi$, the Lipschitz constant $C_{\Phi}$ can be made uniformly bounded, and therefore, the step size $\tau$ can be relaxed from the convergence criteria in~\eqref{convergence_c}.

\subsection{Computing proximal operators}\label{sec:3-3}
We hereby provide details for computing the proximal operators in Algorithm~\ref{algorithm1} for functionals $F$ and $G$ in \eqref{3-25}.

\subsubsection{Proximal operator of \textit{F} for dual variables}
It is important to note that the dual variables $(\phi, \psi)$ and $\mu$ are separable when computing the proximal operator of $F$. We will discuss them as two individual subproblems.

\vspace{6pt}\noindent 1. For dual variable $(\phi, \psi)$, we solve the following problem:
\begin{equation*}
    \min_{\phi, \psi} \sigma \sum_{\textbf{i} \in \mathfrak{I}} \delta_{\mathfrak{K}}(\phi_{\textbf{i}}, \psi_{\textbf{i}}) + \frac{1}{2}\|\phi - \phi_{0}\|^{2} + \frac{1}{2}\|\psi - \psi_{0}\|^{2},
\end{equation*}
which reduces to a pointwise projection of $(\phi_{0}, \psi_{0})$ onto the convex set $\mathfrak{K}$:
\begin{equation}\label{2-9}
    \min_{\phi_{\textbf{i}}, \psi_{\textbf{i}}} \frac{1}{2}\big|\phi_{\textbf{i}} - \phi_{0, \textbf{i}}\big|^{2} + \frac{1}{2}\big|\psi_{\textbf{i}} - \psi_{0, \textbf{i}}\big|^{2} ~~\text{ s.t. }~~\phi_{\textbf{i}} + \frac{1}{2}|\psi_{\textbf{i}}|^{2} \le 0.
\end{equation}
This projection has closed-form solutions as below:

\noindent (i) If the inequality constraint already holds for the initial points $(\phi_{0, \textbf{i}}, \psi_{0, \textbf{i}})$, then the solution is given by $$\phi_{\textbf{i}} = \phi_{0, \textbf{i}},~\psi_{\textbf{i}} = \psi_{0, \textbf{i}}.$$

\noindent (ii) Otherwise, consider the Karush–Kuhn–Tucker (KKT) optimality condition:
\begin{equation}\label{3-13}
    \begin{cases}
        \phi_{\textbf{i}} + \lambda = \phi_{0, \textbf{i}}\\
        (1 + \lambda)\psi_{\textbf{i}} = \psi_{0, \textbf{i}}\\
        \phi_{\textbf{i}} + \tfrac{1}{2}|\psi_{\textbf{i}}|^{2} = 0
    \end{cases}
\end{equation}
where $\lambda$ is the Lagrangian multiplier associated with the constraint. Then the solution is given by 
$$\phi_{\textbf{i}}=\phi_{0,\textbf{i}}-\lambda^{\ast},~ \psi_{\textbf{i}} =\frac{1}{1+\lambda^{\ast}}\psi_{0, \textbf{i}},$$
where $\lambda^{\ast}$ is the largest real root of a third order polynomial:
\begin{equation*}
    (1 + \lambda)^{2}(\phi_{0, \textbf{i}} - \lambda) + \frac{1}{2}|\psi_{0, \textbf{i}}|^{2} = 0,
\end{equation*}
which can be solved efficiently using Newton's method with a relatively large initial guess \cite{Papadakis2013}, or using the root formula for cubic polynomial \cite{Carrillo2019}.

\vspace{6pt}\noindent 2. For the dual variable $\mu$, we solve the minimization problem:
\begin{equation}\label{pu}
    \min_{\mu} F_{0}(\mu) = \sigma\Delta t \mathcal{U}^{\ast}(\mu/\Delta t) + \frac{1}{2}\|\mu - \mu_{0}\|^{2}.
\end{equation}
The general principle for solving this problem is as follows.

\noindent (i) When explicit expression of the convex conjugate $\mathcal{U}^{\ast}$ is available, one can solve~\eqref{pu} using Newton's method to compute its optimality condition  $\nabla_{\mu}F_0(\mu) = 0$. At each step $\mu_{k}$, we update
\begin{equation*}
    \mu_{k + 1} = \mu_{k} - \big(\nabla^{2}_{\mu} F_{0}(\mu_{k})\big)^{-1}\nabla_{\mu} F_{0}(\mu_{k}).
\end{equation*}
The Hessian $\nabla^{2}_{\mu} F_{0}(\mu_{k})$ is positive definite since $\mathcal{U}$ is convex, and is typically diagonal in our examples.

\noindent (ii) Otherwise, one can consider the Moreau’s identity (see for example \cite{Chambolle2011}):
\begin{equation}\label{3-23}
    \text{prox}_{\sigma\Delta t \mathcal{U}^{\ast}(\mu/\Delta t)}(\mu_{0}) = \mu_{0} - \sigma\text{prox}_{\sigma^{-1}\Delta t\mathcal{U}(\mu)}(\mu_{0}/\sigma).
\end{equation}
Then the proximal problem on the RHS
\begin{equation}\label{3-23-1}
    \min_{\mu} \sigma^{-1}\Delta t\mathcal{U}(\mu) + \frac{1}{2}\|\mu - \mu_{0}/\sigma\|^{2}
\end{equation} 
can be solved using Newton's method as above.

\begin{remark}\label{wr3}
Similar to Remark~\ref{writeeng}, when computing the proximal operator~\eqref{pu} based on~\eqref{3-23}, it suffices to solve the optimality condition of~\eqref{3-23-1}:
\begin{equation*}
    \Delta t \frac{\delta \mathcal{U}(\mu)}{\delta \mu} + \sigma \mu = \mu_{0}.
\end{equation*}
Hence, the explicit expression of $\mathcal{U}$ is unnecessary when computing this proximal operator. Instead, it suffices to know the form of $\delta \mathcal{U}(\mu)/\delta \mu$, which is always available from the discrete approximation of $\delta\mathcal{E}(\rho)/\delta \rho$. In Example~\ref{non-e}, we will see that even when the model does not possess an energy functional, one can still apply the PDFB splitting by evaluating the potential field through~\eqref{3-23}.
\end{remark}
\begin{remark}
We also point out that in our numerical examples, proximal problems can be solved efficiently with initial guess $\mu_{0}$ for~\eqref{pu} and initial guess $\mu_{0}/\sigma$ for~\eqref{3-23-1}.
\end{remark}

\subsubsection{Proximal of \textit{G} for primal variable}
For the primal variables, we need to solve the following optimization problem:
\begin{equation*}
    \min_{\rho, \textbf{m}} \tau \delta_{\mathfrak{D}}(\rho, \textbf{m}) + \frac{1}{2}\|\rho - \rho_{0}\|^{2} + \frac{1}{2}\|\textbf{m} - \textbf{m}_{0}\|^{2}.
\end{equation*}
which is a quadratic optimization problem with linear and box constraints:
\begin{equation}\label{3-18}
    \begin{aligned}
        &\min_{\rho, \textbf{m}} \frac{1}{2}\|\rho - \rho_{0}\|^{2} + \frac{1}{2}\|\textbf{m} - \textbf{m}_{0}\|^{2},\\
        &\text{ s.t. }\rho - \rho^{n} + \mathcal{A}\textbf{m} = 0,\quad\beta_{0} \le \rho_{\textbf{i}} \le \beta_{1},\text{ for all } \textbf{i} \in \mathfrak{I}.
    \end{aligned}
\end{equation}
The optimality conditions of the above problem are given by:
\begin{equation}\label{3-3}
    \begin{cases}
        \rho - \rho_{0} + \eta + \lambda = 0\\
        \textbf{m} - \textbf{m}_{0} + \mathcal{A}^{T}\eta = 0\\
        \rho - \rho^{n} + \mathcal{A}\textbf{m} = 0\\
        C(\rho, \lambda) = 0
    \end{cases}
\end{equation}
where $\eta$ and $\lambda$ are the Lagrangian multipliers associated with the linear constraint and box constraint respectively, and $C(\rho, \lambda)$ is the complementary function corresponding to the box constraint:
\begin{equation*}
    C(\rho, \lambda) = \lambda - \big(\lambda + c_{0}(\rho - \beta_{0})\big)_{-} - \big(\lambda + c_{0}(\rho - \beta_{1})\big)_{+}
\end{equation*}
with $c_{0} > 0$ where we fix $c_{0} = 1$ in our computation, where $(v)_{\pm}$ stands for the positive part and negative part of the vector $v$ respectively.

We use the classical primal-dual active set method \citep[Page 319]{Manzoni2021} to solve~\eqref{3-3} with iterations. In each iteration, we update the $(k+1)$-th iterates $(\rho_{k + 1}, \textbf{m}_{k + 1})$ and multiplier $(\eta_{k + 1}, \lambda_{k + 1})$, given the $k$-th active sets
\begin{equation*}
    A^{-}_{k} = \big\{\textbf{i} : \lambda_{k, \textbf{i}} + \rho_{k, \textbf{i}} - \beta_{0} < 0 \big\} \text{ and } A^{+}_{k} = \big\{\textbf{i} : \lambda_{k, \textbf{i}} + \rho_{k, \textbf{i}} - \beta_{1} > 0\big\},
\end{equation*}
and the inactive set $\mathfrak{I}_{k} = \mathfrak{I} - (A^{-}_{k}\cup A^{+}_{k})$. By denoting $P_{A}$ the selection matrix such that $\lambda_{A} = P_{A}\lambda$ is the vector of components of $\lambda$ associated with the index set $A$, we can conversely recover the full multiplier vector, given that $\lambda_{k+1,\textbf{i}}=0$ for $\textbf{i}\in \mathfrak{I}_{k}$, by
\begin{equation}\label{eq:lam_recover}
\lambda_{k + 1} = P^{T}_{A^{-}_{k}}\lambda_{k + 1, A^{-}_{k}} + P^{T}_{A^{+}_{k}}\lambda_{k + 1, A^{+}_{k}}.
\end{equation}
Then the system~\eqref{3-3} can be rewritten in the matrix-vector form
\begin{equation}\label{KKT}
    \left[
    \begin{array}{ccccc}
        I & 0 & I & P^{T}_{A^{-}_{k}} & P^{T}_{A^{+}_{k}}\\
        0 & I & \mathcal{A}^{T} & 0 & 0 \\
        I & \mathcal{A} & 0 & 0 & 0\\
        P_{A^{-}_{k}} & 0 & 0 & 0 & 0\\
        P_{A^{+}_{k}} & 0 & 0 & 0 & 0\\
    \end{array}\right]
    \left[\begin{array}{c}
        \rho_{k + 1}\\
        \textbf{m}_{k + 1}\\
        \eta_{k + 1}\\
        \lambda_{k + 1, A^{-}_{k}}\\
        \lambda_{k + 1, A^{+}_{k}}
    \end{array}\right] = 
    \left[\begin{array}{c}
        \rho_{0}\\
        \textbf{m}_{0}\\
        \rho^{n}\\
        \beta_{0} \textbf{e}_{A^{-}_{k}}\\
        \beta_{1} \textbf{e}_{A^{+}_{k}}
    \end{array}\right]
\end{equation}
where $I$ is the identity matrix, and $\textbf{e}_{A} = (1, \cdots, 1)^{T}$ is a vector with a length equals to the cardinality number $|A|$. By further introducing the matrix and vectors
\begin{equation*}
    \mathcal{J}_{k} =
    \left[\begin{array}{cc}
        I & \mathcal{A} \\
        P_{A^{-}_{k}} & 0 \\
        P_{A^{+}_{k}} & 0 \\
    \end{array}\right], \quad
    \begin{aligned}
        &\textbf{x}_{k + 1} = [\rho_{k + 1}; \textbf{m}_{k + 1}],  &&\textbf{y}_{k + 1} = \big[\eta_{k + 1}; \lambda_{k + 1, A^{+}_{k}}; \lambda_{k + 1, A^{-}_{k}}\big],\\
        &\textbf{x}_{0} = [\rho_{0}; \textbf{m}_{0}],  &&\textbf{y}_{0} = \big[\rho^{n}; \beta_{0} \textbf{e}_{A_{k}^{-}}; \beta_{1} \textbf{e}_{A_{k}^{+}}\big],
    \end{aligned}
\end{equation*}
the system~\eqref{KKT} can be written compactly in a saddle point form:
\begin{equation*}
    \left[\begin{array}{cc}
        I & \mathcal{J}^{T}_{k} \\
        \mathcal{J}_{k} & 0
    \end{array}\right] \left[\begin{array}{c}
        \textbf{x}_{k + 1}\\
        \textbf{y}_{k + 1}
    \end{array}\right] = \left[\begin{array}{c}
        \textbf{x}_{0}\\
        \textbf{y}_{0}
    \end{array}\right].
\end{equation*}
Using a simple computation, the above linear system can be decoupled as
\begin{equation}\label{3-17}
    \begin{aligned}
        \mathcal{J}_{k}\mathcal{J}^{T}_{k}\textbf{y}_{k + 1} &= \mathcal{J}_{k}\textbf{x}_{0} - \textbf{y}_{0},\\
        \textbf{x}_{k + 1} &= \textbf{x}_{0} - \mathcal{J}_{k}^{T}\textbf{y}_{k + 1}.
    \end{aligned}
\end{equation}
In the active set method, we repeatedly solve~\eqref{3-17} to obtain $\textbf{y}_{k + 1}$ and $\textbf{x}_{k + 1}$, and update the active sets $A^{\pm}_{k+1}$ and the matrix $\mathcal{J}_{k+1}$ by recovering the multiplier $\lambda_{k + 1}$ by \eqref{eq:lam_recover} until the active sets remain unchanged, i.e., $A^{\pm}_{k + 1} = A^{\pm}_{k}$. 

\begin{remark}[Sparse linear solver]
It is noted that in each iteration we need to solve a sparse linear system~\eqref{3-17} with a positive definite matrix $\mathcal{J}_{k}\mathcal{J}^{T}_{k}$. The system can be efficiently solved by using banded matrix solver, or applying the preconditioned conjugate gradient method with a block-Jacobi preconditioner:
\begin{equation*}
    \mathcal{M}_{k} = \left[\begin{array}{cc}
        I + \mathcal{A}\mathcal{A}^{T} & 0 \\
        0 & I_{k}
    \end{array}\right]
\end{equation*}
where $I_{k}$ is an identity matrix with length equal to the size of active set $|A_{k}^{-} \cup A_{k}^{+}|$. In particular, the inversion of the preconditioner $\mathcal{M}_{k}$ can be efficiently computed by using fast Fourier transform-based algorithms.
\end{remark}

\section{Numerical Examples}\label{sec:4}
In this section, we explore the capability of the proposed PDFB splitting method for computing Wasserstein-like gradient flows with general nonlinear mobilities~\eqref{1-1}. When needed, the convex splitting technique in Section~\ref{sec:3-4} will be applied. Furthermore, the semi-implicit mobility scheme in Section~\ref{remark7} will be used as comparison in selected examples. The code for generating the numerical results in this section is available at \url{https://github.com/Yunhong0335/PDFB-for-optimal-transport}.

\begin{figure}[t]
    \centering
    \includegraphics[width = \textwidth]{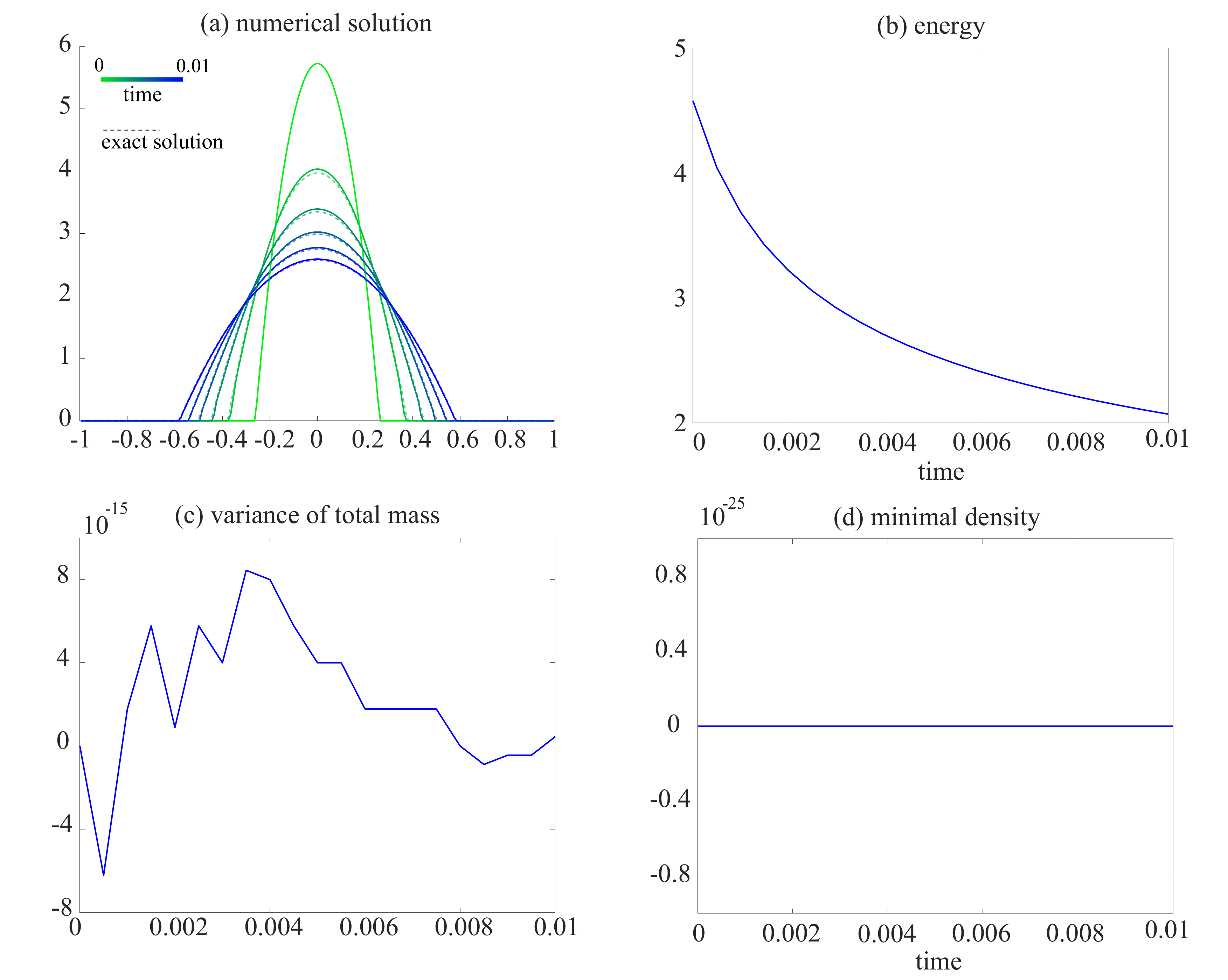}
    \caption{\small{Numerical results for porous media diffusion equation in \textit{Example 5.1}. (a) The numerical solution at several time steps computed using the PDFB splitting method, where the dash line represents the exact solution given in~\eqref{exam_11}. (b) The energy evolution of the numerical solution. (c) The error of the numerical total mass with respect to the initial total mass over time. (d) The minimal value of the numerical solution over time.}}
    \label{exam0}
\end{figure}

\subsection{1D equations}
\begin{example}[\cite{Otto2001}]\label{examp1}\rm
Consider first the porous media equation:
\begin{equation}\label{exam_1}
    \partial_{t}\rho - \nabla \cdot\big(\rho \nabla U^{\prime}(\rho)\big) = 0 \text{ where } U(\rho) = \rho^{2}\,,
\end{equation}
which admits a family of exact solutions of Barenblatt profiles \cite{Vzquez2006}. In particular, we consider the following exact solution as the reference of the numerical solution
\begin{equation}\label{exam_11}
    \rho(x, t) = (t + t_{0})^{-\frac{1}{3}} \bigg(\Big(\frac{3}{16}\Big)^{\frac{1}{3}} - (t + t_{0})^{-\frac{2}{3}}\frac{x^{2}}{12}\bigg)_{+} \text{ with } t_{0} = 10^{-3}.
\end{equation}

In this example, we solve~\eqref{exam_1} on the domain $\Omega = [-1, 1]$ with a grid spacing $h = 0.01$ and a time step size $\Delta t = 0.0005$. To ensure bound preservation $\rho \ge 0$, we set $\beta_{0} = 0$ and $\beta_{1} = \infty$ in the box constraint. The step sizes in the PDFB splitting method are chosen as $\tau = 1$ and $\sigma = 1$, and the tolerance for convergence is set to $10^{-5}$. We also use the midpoint rule for the discrete free energy:
\begin{equation*}
    \mathcal{E}_{h}(\rho) = \sum_{\textbf{i} \in \mathfrak{I}} U(\rho_{\textbf{i}}) h.
\end{equation*}

Numerical results are given in Figure~\ref{exam0}. It can be seen from Figure~\ref{exam0} (a) that the sharp fronts on the boundary of the compact supports are nicely captured in the numerical solution. We also demonstrate the structure-preserving features of the numerical solutions: Figure~\ref{exam0} (b) shows the dissipation of the discrete energy; Figure~\ref{exam0} (c) shows that the total mass can be conserved up to an error close to machine epsilon; Figure~\ref{exam0} (d) shows the positivity-preserving property. 
We also show the rate of convergence in Figure~\ref{exam01}. In particular, we plot the number of iterations needed for convergence with different spatial grid sizes in Figure~\ref{exam01} (b). It is worth noting that the rate of convergence is independent of the spatial grid spacing.

\begin{figure}[h]
    \centering
    \includegraphics[width = \textwidth]{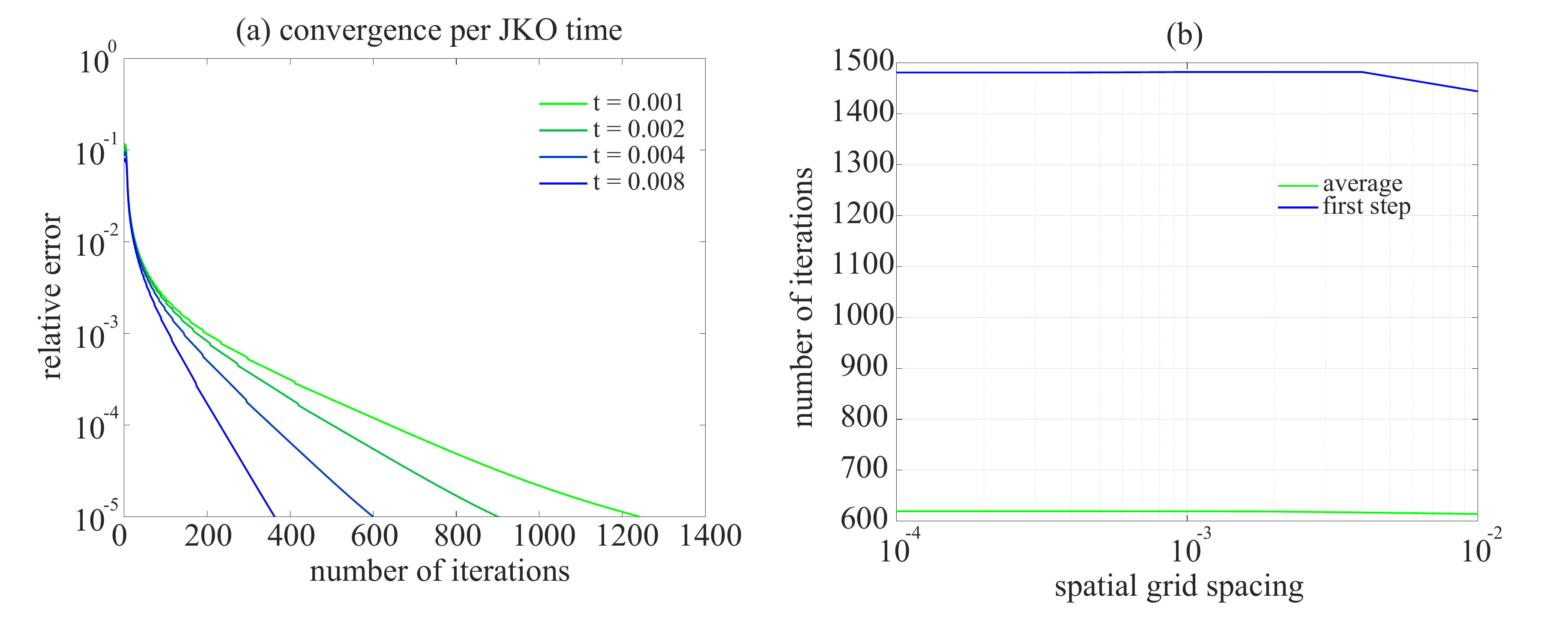}
    \caption{\small{The convergence behavior for \textit{Example 5.1}. (a) The convergence rate of the PDFB splitting method in each time step. (b) The averaged number of iterations over all time steps and the number of iterations in the first time step with different spatial grid spacing. }}
    \label{exam01}
\end{figure}
\end{example}

\begin{example}[\cite{Carrillo2022, Carrillo2023}]\label{examp-2}\rm
Consider the following Fokker-Planck equation with saturation effect:
\begin{equation}\label{example1}
    \partial_{t}\rho - \nabla \cdot\big(\rho(1 - \rho) \nabla (\ln(\rho) + V(x))\big) = 0 \text{ where } V(x) = \frac{|x|^{2}}{2}.
\end{equation}
We use an initial value considered in \cite{Carrillo2023}, which is a uniform density given by
\begin{equation*}
    \rho_{0}(x) = 0.415.
\end{equation*}

In this example, we solve~\eqref{example1} on the domain $\Omega = [-4, 4]$ with a grid spacing $h = 0.02$ and a time step size $\Delta t = 0.1$. To ensure bound preservation $0 \le \rho \le 1$, we set $\beta_{0} = 0$ and $\beta_{1} = 1$ in the box constraint. The step sizes in the PDFB splitting method are chosen as $\tau = 0.2$ and $\sigma = 1/0.2$, and the tolerance for convergence is set to $10^{-7}$. 

The discrete free energy of the model in this example computed as: 
\begin{equation*}
    \mathcal{E}_{h}(\rho) = \sum_{\textbf{i} \in \mathfrak{I}} \big(\rho_{\textbf{i}} \ln(\rho_{\textbf{i}}) - \rho_{\textbf{i}} + V(x_{\textbf{i}}) \rho_{\textbf{i}}\big) h.
\end{equation*}
As the Gibbs-Boltzmann entropy functional $\sum_{\textbf{i} \in \mathfrak{I}}\rho_{\textbf{i}} \ln(\rho_{\textbf{i}}) - \rho_{\textbf{i}}$ is ill-conditioned, we consider the convex splitting approach by taking
\begin{equation*}
    \mathcal{U}(\rho) = \sum_{\textbf{i} \in \mathfrak{I}} \rho_{\textbf{i}} \ln(\rho_{\textbf{i}}) - \rho_{\textbf{i}} \quad \text{and} \quad \mathcal{V}(\rho) = \sum_{\textbf{i} \in \mathfrak{I}} V_{\textbf{i}} \rho_{\textbf{i}}.
\end{equation*}
Here the convex conjugate to $\mathcal{U}$ can be explicitly written as follows:
\begin{equation*}
    \mathcal{U}^{\ast}(\mu) = \sum_{\textbf{i} \in \mathfrak{I}}e^{\mu_{\textbf{i}}}.
\end{equation*}

Additionally, we consider the semi-implicit mobility scheme mentioned in Section~\ref{remark7}. Specifically, we approximate the mobility $\rho(1 - \rho)$ using $\rho(1 - \rho^{n})$ where $\rho^{n}$ is the numerical solution at the previous step. In this scheme, we set $\Delta t = 0.01$ while keeping other parameters the same.

\begin{figure}[t]
    \centering
    \includegraphics[width = \textwidth]{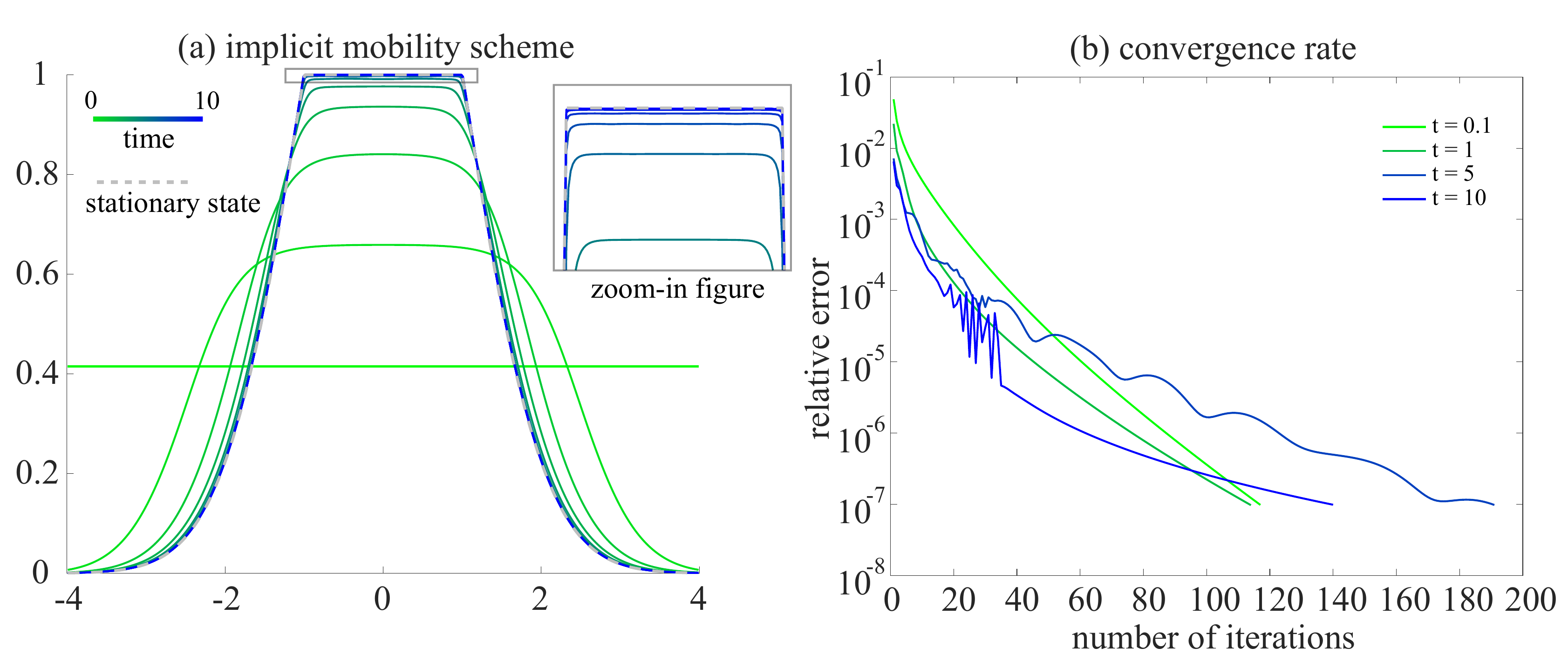}
    \caption{\small{Numerical results for \textit{Example 5.2} with full implicit mobility scheme. (a) The numerical solution at several time steps. Here the dashed gray line represents the stationary state of the equation with the considered initial value, which has been given in \cite{Carrillo2022, Carrillo2023}. (b) The convergence rate of the PDFB splitting method in several time steps.}}
    \label{exam1}
\end{figure}
\begin{figure}[t]
    \centering
    \includegraphics[width = \textwidth]{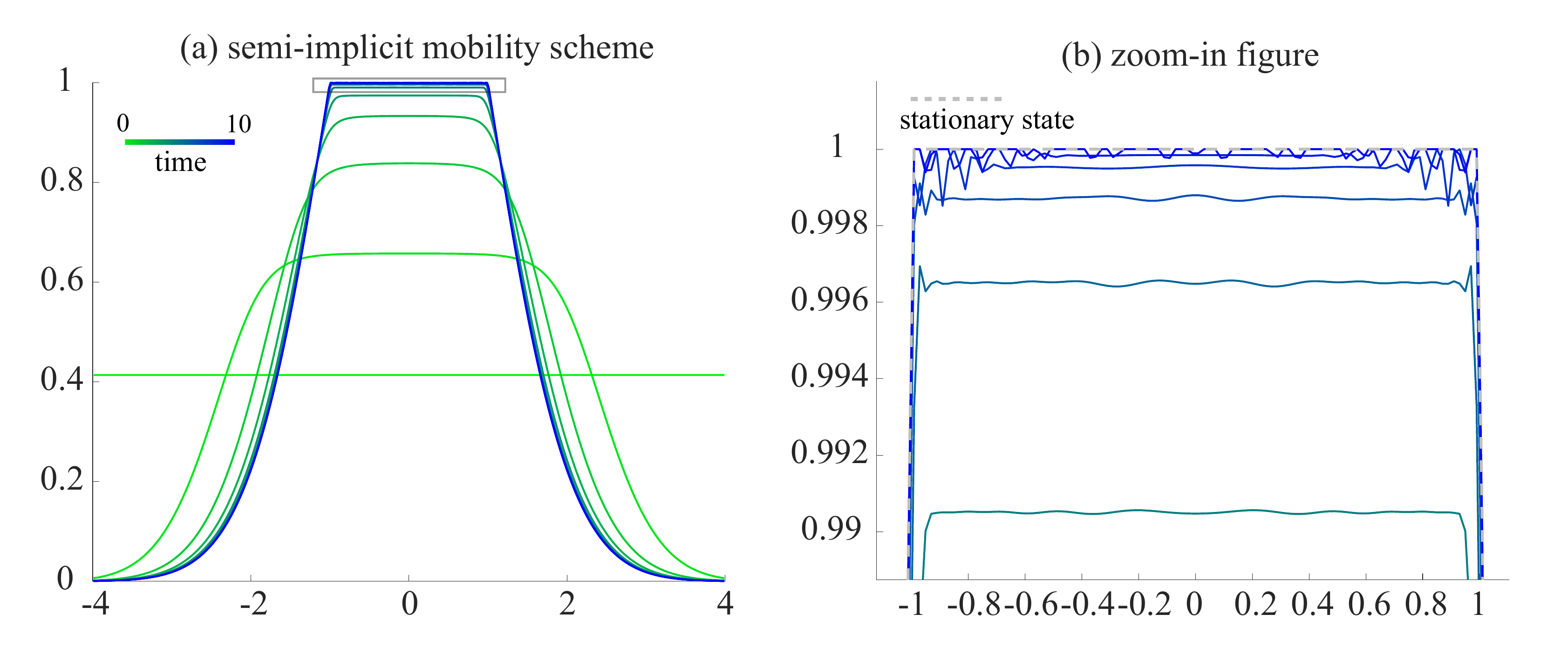}
    \caption{\small{Numerical results for \textit{Example 5.2} with semi-implicit mobility scheme. (a) The numerical solution at several time steps. (b) The zoom-in figure of the numerical solution singled out by the frame.}}
    \label{exam1-2}
\end{figure}
\begin{figure}[t]
    \centering
    \includegraphics[width = \textwidth]{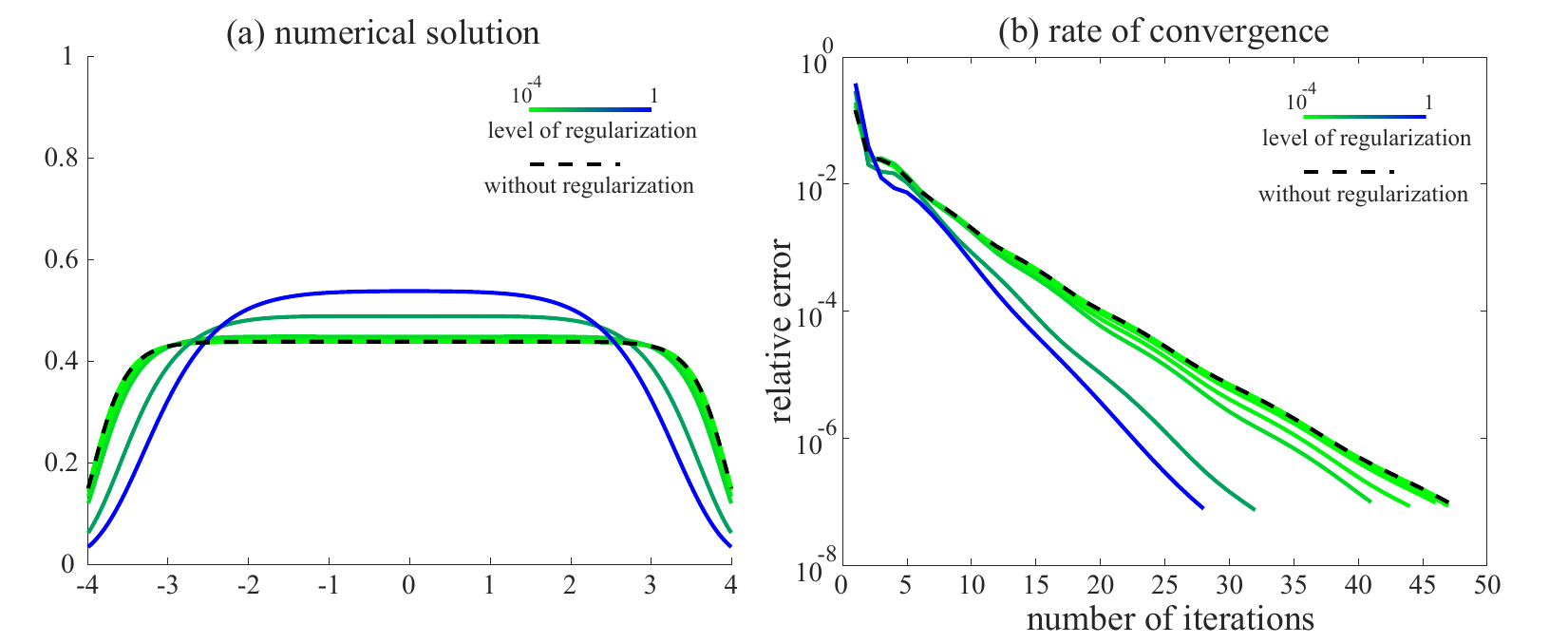}
    \caption{\small{Numerical results for \textit{Example 5.2} with different levels of Moreau-Yosida regularization. In particular, we run a one-step generalized JKO scheme with parameters $\gamma = 5\times 10^{-4}, 10^{-3}, 5\times 10^{-3}, 10^{-2}, 5\times 10^{-2}, 0.1, 0.5, 1$. As the comparison, the black dash-lines represent numerical results of non-regularized generalized JKO scheme with $\tau = 0.6$ and $\sigma = 1/\tau$. The stepsizes are optimized in this example, in particular, we take $\tau = 0.6, 0.6, 0.6, 0.6, 0.7, 0.8, 1.2, 1.6$ with $\sigma = 1/\tau$. (a) The numerical solution at $t = 0.1$ with different levels of regularization. (b) The rate of convergence in the PDFB splitting iteration.}}
    \label{exam1-2-regularization}
\end{figure}

The results of the PDFB splitting method with full implicit mobility and semi-implicit mobility schemes are given in Figure~\ref{exam1} and Figure~\ref{exam1-2}, respectively. Both schemes ensure the bound-preserving property of the numerical solution (i.e., $0\leq \rho \leq 1$). However, no noticeable oscillations are observed in Figure~\ref{exam1} (a) when the solution $\rho$ approaches $1$, using the implicit mobility scheme. In contrast, the results in Figure~\ref{exam1-2} uses the semi-implicit mobility scheme with a relatively small time step and still shows oscillation behavior. Comparing the two, it is evident that the full implicit mobility scheme offers a clear advantage over the semi-implicit scheme, as it allows for larger time step sizes while still providing accurate (non-oscillatory) numerical solutions and reducing computational cost. We also note that, compared to the results in \citep[Figure~1]{Carrillo2023}, which employs the scheme~\eqref{osad} with a finer spatial grid, the method proposed here demonstrates superior performance in reducing oscillations. Besides, to reduce the oscillation in the semi-implicit mobility scheme shown Figure~\ref{exam1-2}, one can use the generalized Schr\"odinger bridge problem considered in \citep[Remark 4]{Carrillo2023}.

As we discussed in Section~\ref{remark7}, the $F(v) = \sum_{\textbf{i}} \delta(\phi_{\textbf{i}}, \psi_{\textbf{i}})$ is not strongly convex to satisfy Assumption~\ref{ccs} for the convergence analysis. In this example, we apply different levels of Moreau–Yosida regularization to the action function $f$ to ensure the strong convexity condition, where the associated convex conjugate is $F_{\gamma}(v) = F(v) + \gamma\|v\|^{2}/2$ with regularization parameter $\gamma>0$. We present in Figure~\ref{exam1-2-regularization} the effects of this regularization on the numerical solution and the convergence rate. It demonstrates that our PDFB splitting method without regularization remains stable and converges reliably in practice even in the absence of strong-convexity. Moreover, it produces more accurate solutions compared to the regularized ones. Specifically, although the regularization enforces strong convexity and hence accelerates the convergence of iterations, it introduces regularization bias in solutions, and a sufficiently small regularization parameter ($\gamma<0.1$) is required to obtain reasonably accurate solutions.
\end{example}

\begin{example}[\cite{Bertozzi1998, Bertozzi2002, Zhornitskaya1999}]\label{examp-3}\rm
Consider the following thin film lubrication-type equation:
\begin{equation}\label{example2}
    \partial_{t}\rho - \nabla\cdot\big(\rho^{3} \nabla(P(\rho) - \Delta\rho)\big) = 0\,,
\end{equation}
where $\rho$ represents the thickness of the thin film, and $P(\rho)$ represents the intermolecular force exerted on the film when the dewetting phenomenon is considered. It is noticeable that the mobility function is non-concave in this example.

First, we set $P(\rho) = 0$ and consider an initial condition given by:
\begin{equation*}
    \rho_{0}(x) = 0.8 - \cos(\pi x) + 0.25 \cos(2\pi x).
\end{equation*}

In this case, we solve~\eqref{example2} on the domain $\Omega = [-1, 1]$ with a grid spacing $h = 0.01$ and a time step size $\Delta t = 0.001$. To ensure bound preservation $\rho \ge 0$, we set $\beta_{0} = 0$ and $\beta_{1} = \infty$ in the box constraint. The step sizes in the PDFB splitting method are chosen as $\tau = 0.03$ and $\sigma = 1/0.03$, and the tolerance for convergence is set to $10^{-7}$.

The gradient of the discrete free energy of the model in this case corresponds to a centered approximation of the Laplacian (only gradient of the discrete free energy is needed to evaluate a forward step, see Remark~\ref{writeeng}). 

Numerical results of the PDFB splitting method are given in Figure~\ref{exam2-1}. We also compare the full implicit mobility scheme with the semi-implicit mobility scheme, which approximates the mobility $\rho^{3}$ by $(\rho^{n})^{2}\rho$. By comparing Figure~\ref{exam2-1} (a) with Figure~\ref{exam2-1} (b), we can observe a noticeable error in the numerical solution obtained using the semi-implicit scheme. This error is likely due to the drastic changes in the mobility $\rho^{3}$ as the solution evolves, further highlighting the advantage of the implicit mobility scheme.

\begin{figure}[h]
    \centering
    \includegraphics[width = \linewidth]{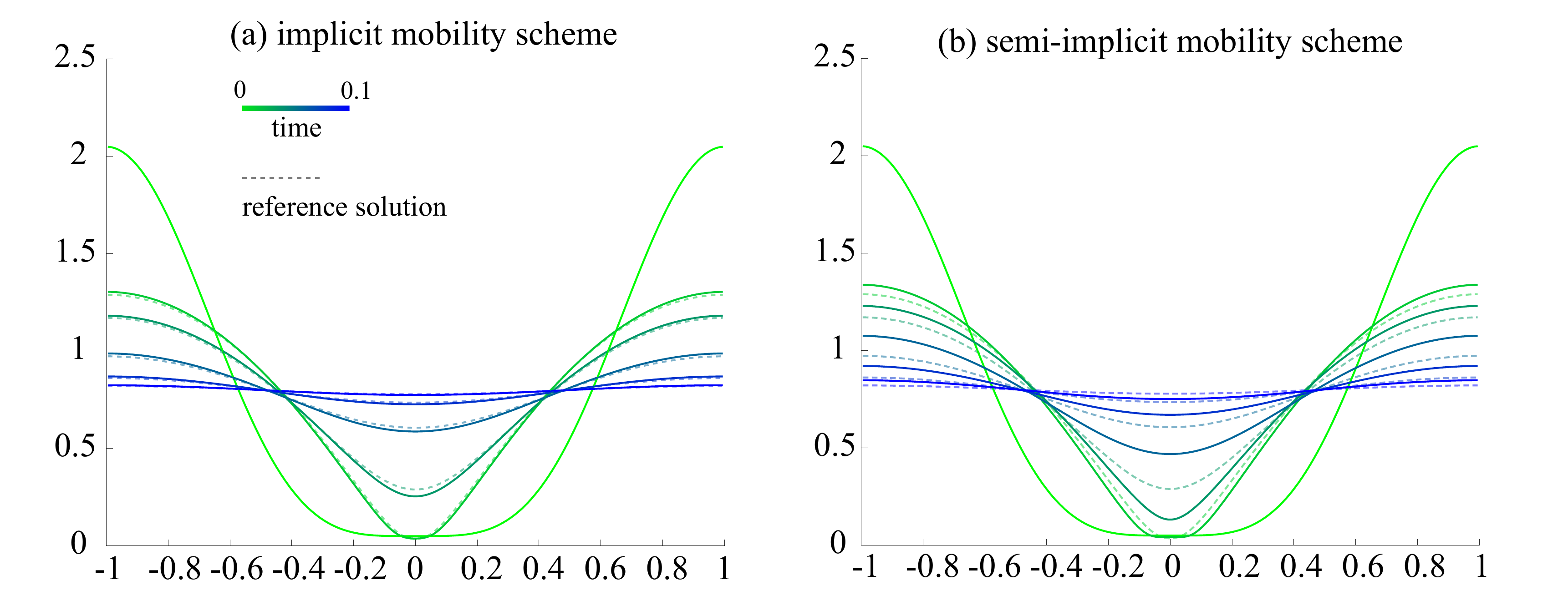}
    \caption{\small{Numerical results for thin film lubrication model in \textit{Example 5.3}. (a) Numerical solution of full implicit mobility scheme at several time steps. The dash line is a reference solution obtained using the typical backward Euler scheme with a step size $\Delta t = 0.0001$. (b) Numerical solution of the semi-implicit mobility scheme at several time steps.}}
    \label{exam2-1}
\end{figure}

Next, we consider two types of intermolecular forces  \cite{Oron1997} given by
\begin{equation}\label{tpr}
    P(\rho) = \rho^{-3} - \varepsilon \rho^{-4} \text{ and } P(\rho) = \rho^{-3} - \varepsilon^{6} \rho^{-9},
\end{equation}
where the first type represents the simultaneous action of the attractive long-range and repulsive short-range intermolecular forces, while the second corresponds to the Lennard-Jones potential for solid-liquid interaction. We fix the parameter $\varepsilon = 0.1$ and consider an initial value given by
\begin{equation*}
    \rho_{0}(x) = (1 - \varepsilon)\frac{\pi}{2} \cos\Big(\frac{\pi x}{2}\Big) + \varepsilon.
\end{equation*}

In this case, we solve~\eqref{example2} on the domain $\Omega = [0, 1]$ with a grid spacing $h = 0.005$ and a time step size $\Delta t = 0.001$. To ensure bound preservation $\rho \ge 0$, we set $\beta_{0} = 0$ and $\beta_{1} = \infty$ in the box constraint. The step sizes in the PDFB splitting method are chosen as $\tau = 0.01$ and $\sigma = 1/0.01$, and the tolerance for convergence is set to $10^{-7}$. 

In both cases, the discrete energy functionals corresponding to the intermolecular forces $P(\rho)$ in~\eqref{tpr} are extremely ill-conditioned as the thickness $\rho$ approaches 0. Therefore, we consider the convex splitting approach by taking
\begin{equation*}
    \frac{\delta \mathcal{U}(\rho)}{\delta \rho} = P(\rho) ~~\text{ and }~~ \frac{\delta \mathcal{V}(\rho)}{\delta \rho} = -\Delta_{h}\rho.
\end{equation*}
Here we treat $\delta \mathcal{U}/\delta \rho$ implicitly and $\delta \mathcal{V}/\delta \rho$ explicitly. As noted in Remark~\ref{writeeng} and Remark~\ref{wr3}, the PDFB splitting method can still be implemented in this example without the need for explicit expressions of $\mathcal{U}$ and $\mathcal{V}$.

Numerical results of the PDFB splitting method are given in Figure~\ref{exam2}. In Figure~\ref{exam2} (a), the numerical solution closely matches with the theoretical approximation of the steady solution \cite{Bertozzi2002}. 

\begin{figure}[h]
    \centering
    \includegraphics[width = \textwidth]{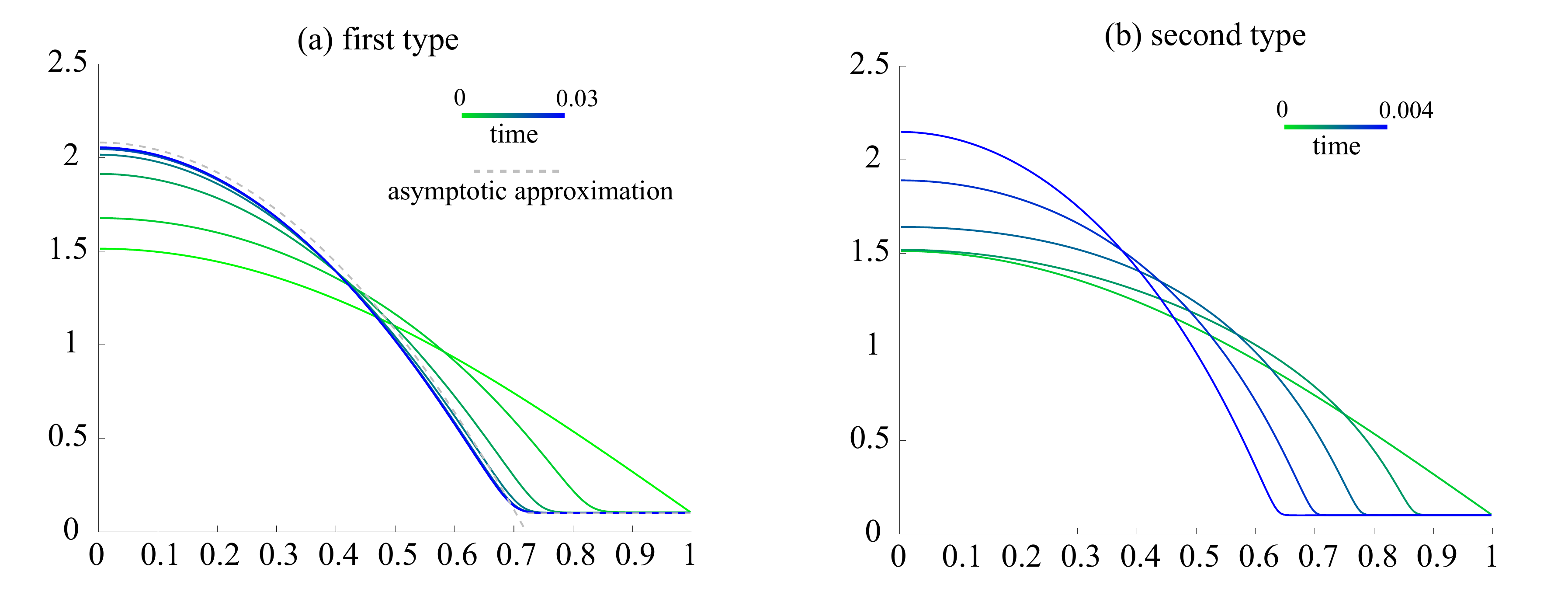}
    \caption{\small{Numerical results for thin film lubrication model in \textit{Example 5.3} with intermolecular force. (a) Numerical solution at several time steps, which corresponds to the first type of intermolecular force $P(\rho)$ in~\eqref{tpr}. The dash gray line represents the asymptotic approximation of the steady solution given in \cite{Bertozzi2002} (b) Numerical solution at several time steps, which corresponds to second type of intermolecular force $P(\rho)$ in~\eqref{tpr}.}}
    \label{exam2}
\end{figure}
\end{example}

\subsection{2D equations}
\begin{example}[\cite{Elliott1996, Cahn1996}]\label{examp4}\rm
Consider the following Cahn-Hilliard equation with degenerate mobility:
\begin{equation}\label{example3}
    \partial_{t}\rho - \nabla\cdot\big((1 - \rho^{2})\nabla(\rho(\rho^{2} - 1) - \varepsilon^{2} \Delta\rho)\big) = 0.
\end{equation}
The free energy functional associated with this equation is
\begin{equation*}
    \mathcal{E}(\rho) = \int_{\Omega} \frac{1}{4}\big(\rho(x)^{2} - 1\big)^{2} + \frac{\varepsilon^{2}}{2}|\nabla \rho(x)|^{2} \,\text{d}x.
\end{equation*}

In this example, we solve~\eqref{example3} on the domain $\Omega = [0, 1] \times [0, 1]$ with a grid spacing $h = 1/64$ and a time step size $\Delta t = 0.001$. We initialize the solution with a randomized value, drawn from a uniform distribution on the discrete domain. To ensure bound preservation $-1 \le \rho \le 1$, we set $\beta_{0} = -1$ and $\beta_{1} = 1$ in the box constraint. The step sizes in the PDFB splitting method are chosen as $\tau = 20$ and $\sigma = 1/20$, and the tolerance for convergence is set to $10^{-7}$. In each iteration of the PDFB splitting method, we do not consider the convex splitting, and only evaluate the explicit gradient of the discrete free energy as follows,
\begin{equation*}
    \bigg(\frac{\delta \widehat{\mathcal{E}}_{h}(\rho)}{\delta \rho}\bigg)_{\textbf{i}} = \rho_{\textbf{i}}(\rho^{2}_{\textbf{i}} - 1) - \varepsilon^{2} (\Delta_{h}\rho)_{\textbf{i}} \text{ for all } \textbf{i} \in \mathfrak{I}.
\end{equation*}

Numerical results of the PDFB splitting method are given in Figure~\ref{exam30}. In addition, we compare the rate of convergence between the proposed PDFB splitting method and the previous preconditioned PD3O method \cite{Carrillo2023} in Figure~\ref{exam3}. We can see that the PDFB splitting method is approximately twice faster than the preconditioned PD3O method. 

\begin{figure}[t]
    \centering
    \includegraphics[width = \textwidth]{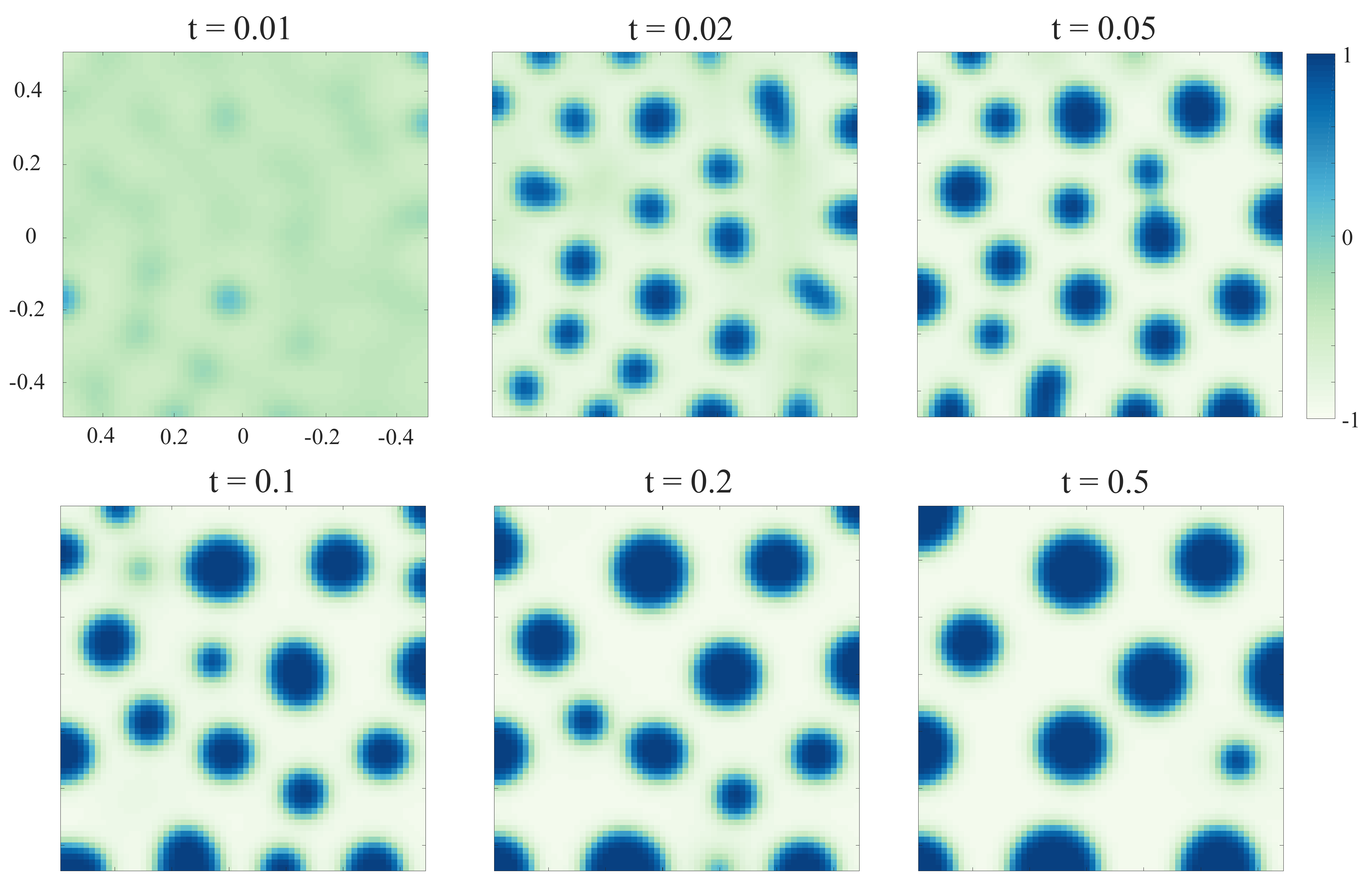}
    \caption{\small{Evolution of numerical solutions in \textit{Example 5.4} by using the PDFB splitting method.}}
    \label{exam30}
\end{figure}
\begin{figure}[h!]
    \centering
    \includegraphics[width = \textwidth]{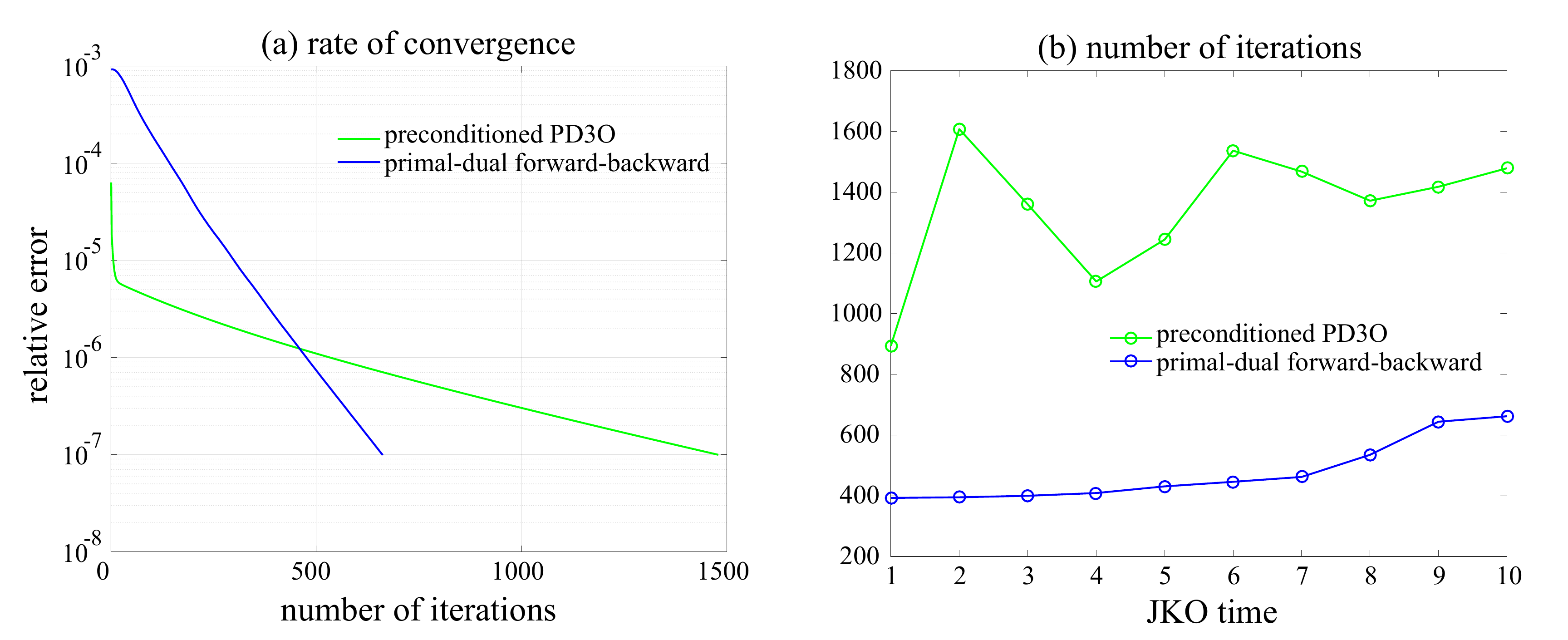}
    \caption{\small{The convergence of the PDFB splitting method and preconditioning PD3O method for \textit{Example 5.4}. (a) The comparison between the two methods on the rate of convergence at the last JKO time step. (b) The number of iterations needed for the relative error to achieve the tolerance for convergence in each JKO time step. Here, the parameters for preconditioned PD3O method is the same as that given in \citep[Figure 9]{Carrillo2023}.}}
    \label{exam3}
\end{figure}
\end{example}

\begin{example}[\cite{Cheng2020, Zhou2023}]\rm
Consider the following Cahn-Hilliard equation with degenerate mobility and anisotropic surface energy:
\begin{equation}\label{4-9}
    \partial_{t}\rho - \nabla\cdot \bigg((1 - \rho^{2})\nabla \frac{\delta \mathcal{E}(\rho)}{\delta \rho}\bigg) = 0,
\end{equation}
where $\mathcal{E}$ is the biharmonic regularized anisotropic Kobayashi-type free energy \cite{Cheng2020}:
\begin{equation*}
    \mathcal{E}(\rho) = \int_{\Omega} \frac{1}{4}\big(\rho(x)^{2} - 1\big)^{2} + \frac{\varepsilon^{2}}{2} \gamma(\boldsymbol{q}(x))^{2}|\nabla \rho(x)|^{2} + \frac{\beta\varepsilon^{2}}{2}(\Delta \rho(x))^{2}\,\text{d}x\,.
\end{equation*}
In this model, $\gamma(\boldsymbol{q})$ accounts for the surface energy anisotropy that depends on the local normal
\begin{equation*}
    \boldsymbol{q}(x) = 
    \begin{cases}
        \displaystyle\frac{\nabla\rho(x)}{|\nabla\rho(x)|} &\text{ if }|\nabla\rho(x)| > 0\\
        0 &\text{ otherwise, }
    \end{cases}
\end{equation*}
and $\beta>0$ is the regularization parameter.
We consider an initial value given by
\begin{equation*}
    \rho_{0}(x) = -\tanh\bigg(\frac{|x| - 0.3}{\varepsilon_{0}}\bigg) ~\text{ with } ~\varepsilon_{0} = 0.025.
\end{equation*}
In the following, we fix the parameters $\beta = 10^{-4}$ and $\varepsilon = 0.01$.

In this example, we solve~\eqref{4-9} on the domain $\Omega = [-\frac{1}{2}, \frac{1}{2}] \times [-\frac{1}{2}, \frac{1}{2}]$ with a grid spacing $h = 1/128$ and a time step size $\Delta t = 0.001$. To ensure bound preservation $-1 \le \rho \le 1$, we set $\beta_{0} = -1$ and $\beta_{1} = 1$ in the box constraint. The step sizes in the PDFB splitting method are chosen as $\tau = 2$ and $\sigma = 1/2$, and the tolerance for convergence is set to $10^{-5}$. The discrete free energy is computed as:
\begin{equation*}
    \mathcal{E}_{h}(\rho) =  \sum_{\textbf{i} \in \mathfrak{I}} \bigg(\frac{1}{4}(\rho_{\textbf{i}}^{2} - 1)^{2} + \frac{\varepsilon^{2}}{2} \gamma(\boldsymbol{q}_{\textbf{i}})^{2}|(\nabla_{h} \rho)_{\textbf{i}}|^{2} + \frac{\beta\varepsilon^{2}}{2}(\Delta_{h}\rho)_{\textbf{i}}^{2}\bigg)h^{2}\,,
\end{equation*}
where $\nabla_{h}$ is the central difference operator defined by:
\begin{equation*}
    (\nabla_{h}\rho)_{\textbf{i}} = \bigg(\frac{\rho_{\mathtt{i} + 1, \mathtt{j}} - \rho_{\mathtt{i} - 1, \mathtt{j}}}{2h}, \frac{\rho_{\mathtt{i}, \mathtt{j} + 1} - \rho_{\mathtt{i}, \mathtt{j} - 1}}{2h}\bigg) \text{ for all } \textbf{i} = (\mathtt{i}, \mathtt{j}) \in \mathfrak{I}\,.
\end{equation*}
We impose a periodic boundary condition on $\rho$, which makes the central difference operator well-defined on the boundary grids.

Since the biharmonic regularization term is ill-conditioned, we consider the convex splitting approach by taking
\begin{equation*}
    \mathcal{U}(\rho) = \sum_{\textbf{i}\in \mathfrak{I}} \frac{\beta\varepsilon^{2}}{2}(\Delta_{h}\rho)_{\textbf{i}}^{2} 
    \quad \text{and} \quad 
    \mathcal{V}(\rho) = \sum_{\textbf{i} \in \mathfrak{I}} \frac{1}{4}(\rho_{\textbf{i}}^{2} - 1)^{2} + \frac{\varepsilon^{2}}{2} \gamma(\boldsymbol{q}_{\textbf{i}})^{2}|(\nabla_{h} \rho)_{\textbf{i}}|^{2}.
\end{equation*}
The proximal operator of $\mathcal{U}^{\ast}$ can be evaluated using the Moreau’s identity~\eqref{3-23}, and the gradient of $\mathcal{V}$ is given in \citep[Lemma 3.4]{Zhou2023}.

\begin{figure}[h!]
    \centering
    \includegraphics[width = \textwidth]{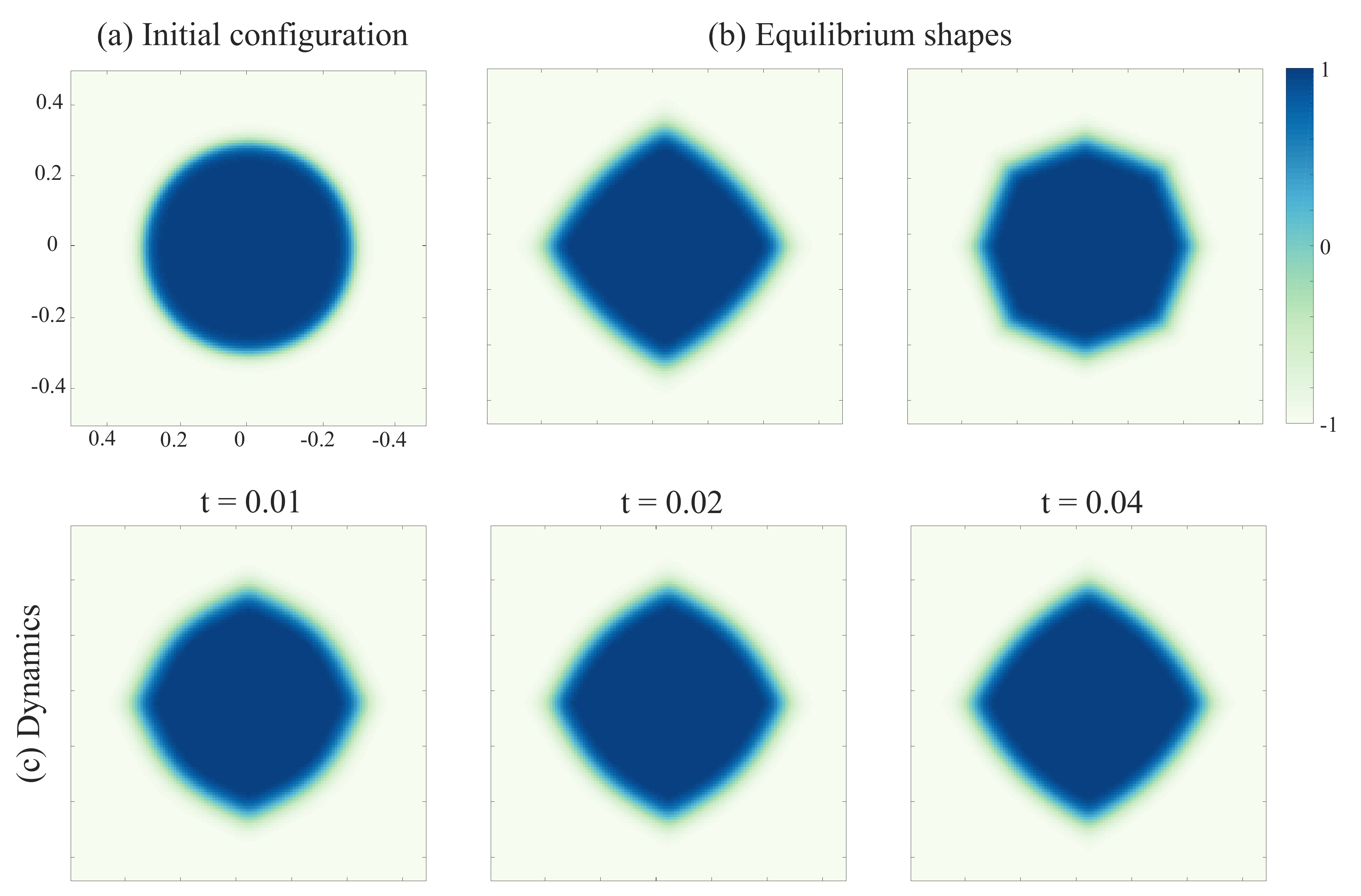}
    \caption{\small{Numerical results in \textit{Example 5} with surface densities of the type~\eqref{4-12} and~\eqref{4-13}. (a) Initial configuration. (b) Equilibrium shapes for surface energies in~\eqref{4-12} and~\eqref{4-13} from left to right, respectively, computed with relatively large time horizon. (c) The dynamics of numerical solution governed by the surface density~\eqref{4-12}.}}
    \label{fig:3-1}
\end{figure}

\begin{figure}[h]
    \centering
    \includegraphics[width = \textwidth]{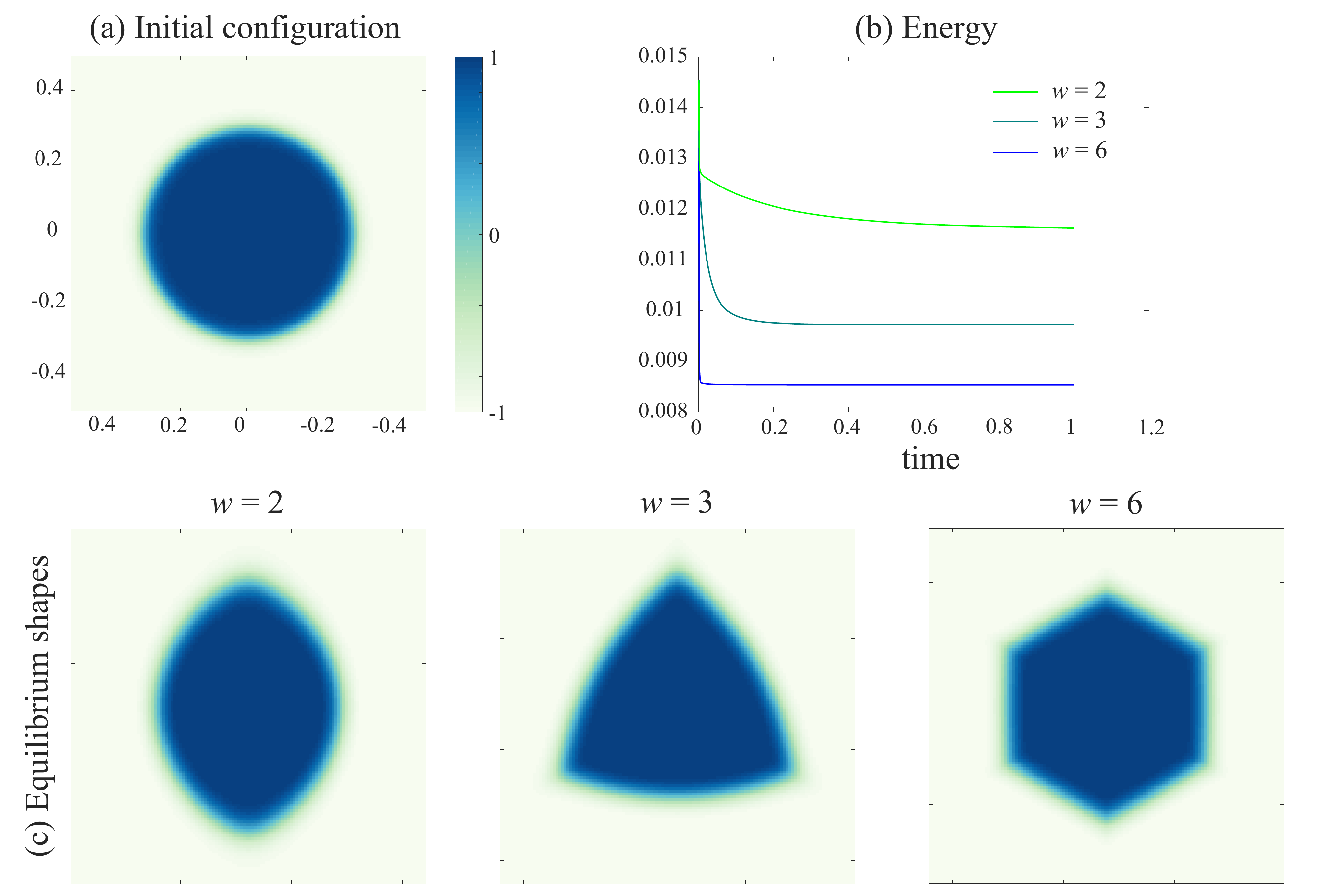}
    \caption{\small{Numerical results in \textit{Example 5} with surface densities of the type~\eqref{4-15}. (a) Initial configuration. (b) The energy evolution of the numerical solutions with different $\omega$. (c) Equilibrium shapes computed at the time $t = 1$ with different $\omega$.}}
    \label{fig:3-2}
\end{figure}

First, we solve~\eqref{4-9} with strongly anisotropic surface energy densities with different symmetries \cite{Cheng2020, Zhou2023} given by
\begin{equation}\label{4-12}
    \gamma(\boldsymbol{q}) = 1 + \alpha \big(4q_{1}^{4} + 4 q_{2}^{4} - 3\big),
\end{equation}
\begin{equation}\label{4-13}
    \gamma(\boldsymbol{q}) = 1 + \alpha \big(8(8q_{1}^{8} - 10q_{1}^{6} + q_{1}^{4}) + 8(8q_{2}^{8} - 10q_{2}^{6} + q_{2}^{4}) + 9\big),
\end{equation}
with $\alpha = 0.2$, where $\boldsymbol{q} = (q_{1}, q_{2})$. The equilibrium shapes of these two energies correspond to a four-fold and an eight-fold shape, respectively (see \cite{Cheng2020}). Numerical results of the equilibrium shapes and the dynamics of the four-fold shape are shown in Figure~\ref{fig:3-1}.

Next, we solve~\eqref{4-9} with the following surface densities \cite{Zhou2023}:
\begin{equation}\label{4-15}
    \gamma(\boldsymbol{q}) = 1 + \alpha \cos(\omega\theta) \text{ with } \tan(\theta) = \frac{q_{1}}{q_{2}}.
\end{equation}
with $\alpha = 0.4$, where $\boldsymbol{q} = (q_{1}, q_{2})$. The equilibrium shapes of this type of energies are shown in \cite{Zhou2023}. Numerical results of the equilibrium shapes and the energy evolution are shown in Figure~\ref{fig:3-2}.
\end{example}

\begin{example}[\cite{Rtz2006}]\label{non-e}\rm
Consider the following doubly degenerate diffuse interface model:
\begin{equation}\label{example5}
    \partial_{t}\rho - \nabla \cdot\bigg((1 - \rho^{2})^{2} \nabla \frac{\rho(\rho^{2} - 1) - \varepsilon^{2} \Delta \rho}{(1 - \rho^{2})^{2}}\bigg) = 0,
\end{equation}
where the denominator $(1 - \rho^{2})^{2}$ in the potential field acts as a diffusion-preventing term \cite{Bretin2020}, eliminating the bulk diffusion effect.  

In this example, the computational domain is $\Omega = [-\frac{1}{2},\frac{1}{2}]\times[-\frac{1}{2},\frac{1}{2}]$, with a grid spacing $h = 1/64$ and a time step size $\Delta t = 0.1$. To ensure bound preservation $-1 \le \rho \le 1$, we set $\beta_{0} = -1$ and $\beta_{1} = 1$ in the box constraint. The step sizes in the PDFB splitting method are chosen as $\tau = 0.2$ and $\sigma = 1/0.2$, and the tolerance for convergence is set to $10^{-5}$. In the following, we fix the parameter $\varepsilon = 0.02$. 

Although equation~\eqref{example5} is not a variational model (or gradient flow), the PDFB splitting method can still be implemented since the computation only requires the information of the ``potential field" (see Remark~\ref{wr3}). Therefore, from a computational perspective, we define a modified “free energy" $\mathcal{U}$ whose first variation is the potential field of the model such that
\begin{equation*}
    \bigg(\frac{\delta \mathcal{U}(\rho)}{\delta \rho}\bigg)_{\textbf{i}} = \frac{\rho_{\textbf{i}}(\rho^{2}_{\textbf{i}} - 1) -\varepsilon^{2}(\Delta_{h}\rho)_{\textbf{i}}}{(1 - \rho^{2}_{\textbf{i}})^{2}} \text{ for all } \textbf{i} \in \mathfrak{I}.
\end{equation*}
According to~\eqref{3-23}, in each step of the iteration, we use Newton method to solve:
\begin{equation*}
    \Delta t \frac{\delta \mathcal{U}(\mu)}{\delta \mu} + \sigma \mu - \mu_{0} = \Delta t \frac{\mu(\mu^{2} - 1) -\varepsilon^{2}\Delta_{h}\mu}{(1 - \mu^{2})^{2}} + \sigma \mu - \mu_{0} = 0,
\end{equation*}
which can be solved efficiently by using FFT-based fast matrix inversion after several iterations with initial guess $\mu = \mu_{0}$.

In Figure~\ref{fig:3-3}, we simulate the dynamics of phase field for the doubly degenerate diffuse interface model~\eqref{example5} and the degenerate Cahn-Hilliard equation~\eqref{example3} with the same initial configuration (Figure~\ref{fig:3-3} (a)) and same parameters. To highlight the differences between the two models, we also track the evolution of the free energy functional: 
\begin{equation}\label{GH}
    \mathcal{E}(\rho) = \int_{\Omega} \frac{1}{4}\big(\rho(x)^{2} - 1\big)^{2} + \frac{\varepsilon^{2}}{2} |\nabla \rho(x)|^{2}\,\text{d}x\,.
\end{equation}
We observe two different dynamics driven by two distinct diffusion regimes. As shown in Figure~\ref{fig:3-3} (c), for the degenerate Cahn-Hilliard equation with double-well potential, the smaller object is absorbed by the larger one due to the additional bulk diffusion, which is reflected in the drastic change in the free energy at $t = 1.5$, corresponding to a change in the topological structure. However, for the doubly degenerate diffuse interface model in Figure~\ref{fig:3-3} (d), the two objects remain separate, with no merging, demonstrating the interface dynamics driven by surface diffusion without bulk diffusion.
\begin{figure}[h!]
    \centering
    \includegraphics[width = \textwidth]{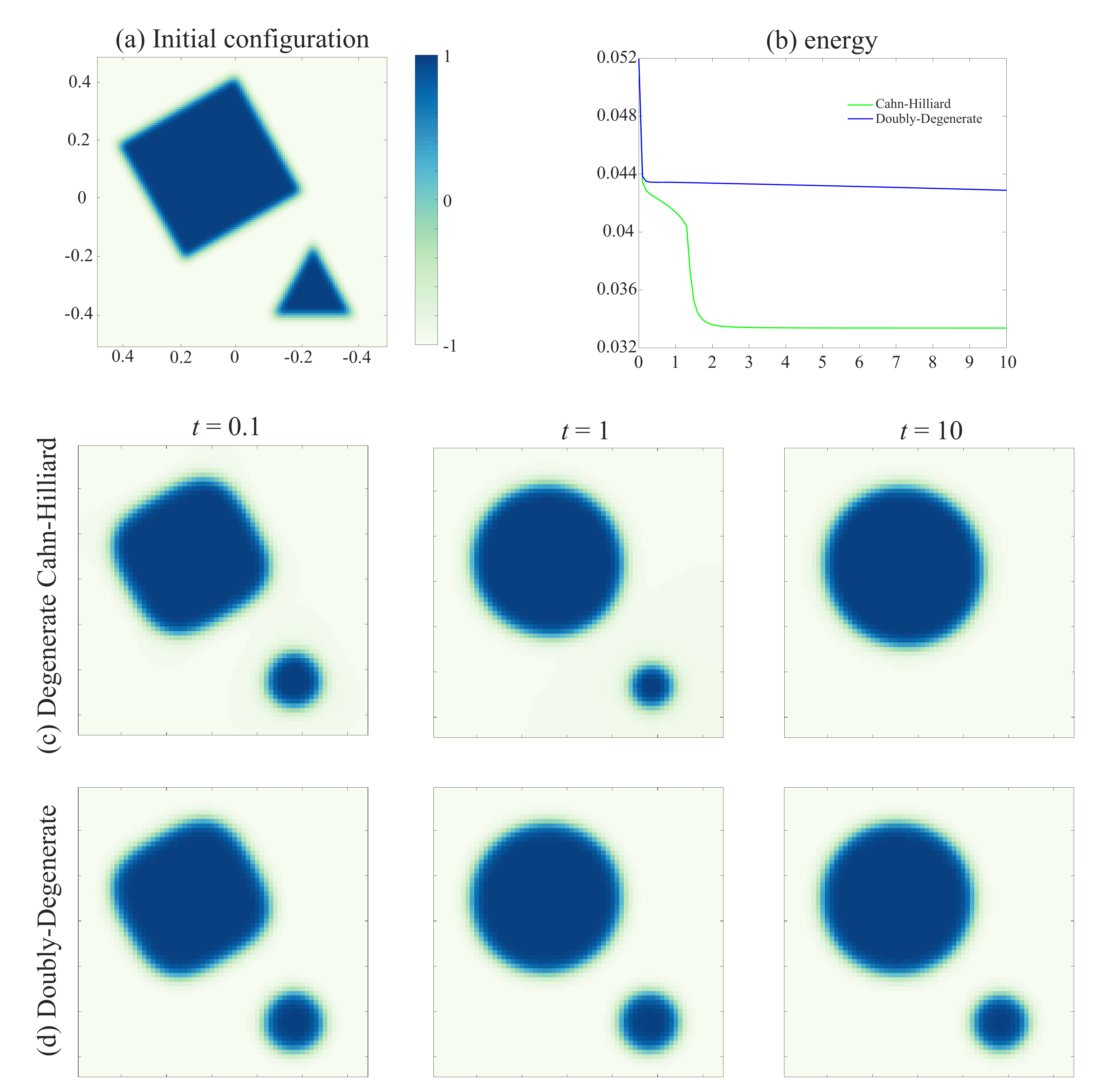}
    \caption{\small{The numerical results for the diffusion-driven interface dynamics in \textit{Example 6}. (a) Initial configuration. (b) The energy evolution for the numerical solutions of the two different equations. (c) The dynamics of numerical solution of the degenerate Cahn-Hilliard equation~\eqref{example3}. (d) The dynamics of numerical solution of the doubly degenerate diffuse interface model~\eqref{example5}.}}
    \label{fig:3-3}
\end{figure}
\end{example}

\section{Conclusion}\label{sec5}
In this work, we construct a new saddle point formulation for the generalized JKO scheme of Wasserstein-like gradient flows, and propose a PDFB splitting method for solving the resulting saddle point problem. We provide a convergence analysis of the proposed PDFB splitting method and discuss several key remarks regarding its application to Wasserstein-like gradient flows. Moreover, we offer implementation details on computing the proximal operators and utilizing a convex splitting nonlinear preconditioning technique. Compared to the previous work \cite{Carrillo2023}, the current framework can be easily adapted to gradient flows with general nonlinear mobilities and the algorithm exhibits a convergence rate independent of the grid size. The efficiency of the PDFB splitting method is demonstrated through several challenging models in the literature.

Concerning the future direction of the PDFB splitting method, it is natural to explore preconditioning techniques to further enhance convergence (e.g., \cite{Jacobs2018, chen2023TPD}). Since the JKO scheme is only first order accurate in time, developing high-order variational schemes that preserve the desired structure while maintaining computational efficiency remains both important and challenging. Additionally, the current PDFB splitting method is flexible and can be adapted to more complex spatial approximations such as high-order finite volume and finite element methods, to improve the spatial accuracy. These, along with other potential avenues, will be explored in our future research.

\section*{Acknowledgments}
The authors would like to thank Zhen Zhang and Wenxing Zhang for their useful observations and suggestions, which greatly improved this work.

\begin{appendix}
\section{Boundary condition}\label{api}
The no-flux boundary condition is given by setting the momentum variables to zero on the faces of boundaries. For $d = 2$, one has
\begin{equation*}
    \begin{aligned}
        m_{\frac{1}{2}, \mathtt{j}} &= m_{\mathtt{n} + \frac{1}{2}, \mathtt{j}} = 0,~\text{ for all } 1 \le \mathtt{j} \le \mathtt{n}\,;\\
        m_{\mathtt{i}, \frac{1}{2}} &= m_{\mathtt{i}, \mathtt{n} + \frac{1}{2}} = 0,~\text{ for all } 1 \le \mathtt{i} \le \mathtt{n}\,,
    \end{aligned}
\end{equation*}
where $\mathtt{n}$ is the number of total volumes in each direction. In practice, we directly embed this condition in constructing the matrices $\mathcal{A}$ and $\mathcal{I}$ appeared in \eqref{constraint} and {2-3}. More particularly, for $d = 1$, we omit the zero momentum variables $m_{\frac{1}{2}}$ and $m_{\mathtt{n} + \frac{1}{2}}$ at two boundary points and only consider the momentum on inner faces in the computation, that is: 
\begin{equation*}
    \textbf{m} = \big(m_{\mathtt{i} + \frac{1}{2}}\big)_{\mathtt{i} = 1}^{\mathtt{n} - 1}.
\end{equation*}
Correspondingly, the matrices are:
\begin{equation*}
    \mathcal{A} = \frac{1}{h}\left[
    \begin{array}{cccc}
        1 & ~ & ~ & ~ \\
        -1 & 1 & ~ & ~ \\
        ~ & \cdots &  & ~ \\
        ~ & ~ & -1 & 1 \\
        ~ & ~ & ~ & -1 \\
    \end{array}\right] \in \mathbb{R}^{\mathtt{n} \times \mathtt{n} - 1}\quad \mathcal{I} = \frac{1}{2}\left[
    \begin{array}{ccccc}
        1 & ~ & ~ & ~ \\
        1 & 1 & ~ & ~ \\
        ~ & \cdots & ~ \\
        ~ & ~ & 1 & 1 \\
        ~ & ~ & ~ & 1 \\
    \end{array}\right] \in \mathbb{R}^{\mathtt{n} \times \mathtt{n} - 1}.
\end{equation*}
High dimensional matrices can be constructed using the Kronecker tensor product.

\section{Proof of the structure-preserving theorem}\label{ast}
In this part, we prove Theorem~\ref{ts}. By taking the maximum of the dual variable $(\phi, \psi)$ in~\eqref{2-8}, the discrete saddle point form is equivalent to the following discrete JKO problem
\begin{equation}\label{2-11}
    \min_{\rho} \mathcal{E}_{h}(\rho) + \frac{1}{2 \Delta t}W_{M, h}(\rho^{n}, \rho)^{2},
\end{equation}
where $\frac{1}{2}W_{M, h}(\rho^{n}, \rho)^{2}$ is the minimum value of the following constrained optimization problem:
\begin{equation*}
    \begin{aligned}
         &\min_{\textbf{m}} \sum_{\textbf{i} \in \mathfrak{I}} f\Big(M\big(\tfrac{1}{2}(\rho^{n}_{\textbf{i}} + \rho_{\textbf{i}})\big), (\mathcal{I}\textbf{m})_{\textbf{i}}\Big)h^{d}\\
        &\text{ s.t. } \rho_{\textbf{i}} - \rho^{n}_{\textbf{i}} + (\mathcal{A}\textbf{m})_{\textbf{i}} = 0 \text{ and } \beta_{0} \le \rho_{\textbf{i}} \le \beta_{1} \text{ for all } \textbf{i} \in \mathfrak{I}.
    \end{aligned}
\end{equation*}

\subsection{Energy dissipation}
Since $\rho^{n + 1}$ is the minimizer~\eqref{2-11}, one has
\begin{equation}\label{th1-1}
    \mathcal{E}_{h}(\rho^{n + 1}) + \frac{1}{2\Delta t}W_{M, h}(\rho^{n}, \rho^{n + 1})^{2} \le \mathcal{E}_{h}(\rho^{n}) + \frac{1}{2\Delta t}W_{M, h}(\rho^{n}, \rho^{n})^{2}.
\end{equation}
Since
\begin{equation*}
    W_{M, h}(\rho^{n}, \rho^{n + 1}) \ge 0 \text{ and } W_{M, h}(\rho^{n}, \rho^{n}) = 0,
\end{equation*}
it is clearly from \eqref{th1-1} that we have
\begin{equation}\label{3-10}
    \mathcal{E}_{h}(\rho^{n + 1}) \le \mathcal{E}_{h}(\rho^{n + 1}) + \frac{1}{2\Delta t}W_{M, h}(\rho^{n}, \rho^{n + 1})^{2} \le \mathcal{E}_{h}(\rho^{n}).
\end{equation}

\subsection{Mass conservation}
There exists a momentum field $\textbf{m}$ with vanishing normal fluxes such that the discrete continuity equation is satisfied:
\begin{equation}\label{3-5}
    \rho^{n + 1}_{\textbf{i}} - \rho^{n}_{\textbf{i}} + (\mathcal{A}\textbf{m})_{\textbf{i}} = 0 \text{ for all } \textbf{i} \in \mathfrak{I}.
\end{equation}
By summing~\eqref{3-5} over all indices $\textbf{i}$, one can show that
\begin{equation}\label{3-8}
    \sum_{\textbf{i} \in \mathfrak{I}} \big(\rho_{\textbf{i}}^{n + 1} - \rho_{\textbf{i}}^{n}\big) = -\sum_{\textbf{i} \in \mathfrak{I}} (\mathcal{A}\textbf{m})_{\textbf{i}}.
\end{equation}
Since the normal fluxes of the momentum field $\textbf{m}$ is vanishing on the boundary, one has
\begin{equation}\label{3-9}
    \sum_{\textbf{i} \in \mathfrak{I}} (\mathcal{A}\textbf{m})_{\textbf{i}}h^{d} = \int_{\partial \Omega} \textbf{n} \cdot \textbf{m} \text{d}s = 0,
\end{equation}
where $\textbf{n} \cdot \textbf{m}$ in the integral should be understood as a piece-wise constant function on the boundary. Combining~\eqref{3-8} and~\eqref{3-9}, we have $\sum_{\textbf{i} \in \mathfrak{I}} \big(\rho_{\textbf{i}}^{n + 1} - \rho_{\textbf{i}}^{n}\big) = 0$.

\subsection{Bound preservation}
The bound preservation is a straightforward result of the constraint $\beta_{0} \le \rho_{\textbf{i}} \le \beta_{1}$ for all $\textbf{i} \in \mathfrak{I}$ in~\eqref{2-11}.

\section{Proof of the convergence theorem}\label{appendix}
In this part, we prove Theorem~\ref{convergence}. Recall that our goal is to solve the following saddle point problem:
\begin{equation*}
    \min_{u}\max_{v} \mathscr{L}(u, v) = \Phi(u, v) - F(v) + G(u)
\end{equation*}
by the PDFB splitting method (Algorithm~\ref{algorithm1}), which can be rewritten in the following form in analogous to the Davis-Yin algorithm \cite{Davis2015, Yan2016}:
\begin{equation}\label{A1}
    \begin{aligned}
        u^{(\ell)} &= \text{prox}_{\tau G}(q^{(\ell)})\\
        v^{(\ell + 1)} &= \text{prox}_{\sigma F}\big(v^{(\ell)} + \sigma \mathcal{K}(u^{(\ell)}) \overline{u}^{(\ell)} + \sigma L(u^{(\ell)})\big)\\
        q^{(\ell + 1)} &= u^{(\ell)} - \tau \nabla_{u}\Phi(u^{(\ell)}, v^{(\ell + 1)})
    \end{aligned}
\end{equation}
where $\mathcal{K}(u^{(\ell)})$ is the Jacobi matrix in~\eqref{k}, $L(u^{(\ell)}) = \nabla_{v}\Phi(u^{(\ell)}, v^{(\ell)}) - \mathcal{K}(u^{(\ell)})u^{(\ell)}$, and
\begin{equation*}
    \overline{u}^{(\ell)} = 2u^{(\ell)} - q^{(\ell)} - \tau \nabla_{u}\Phi(u^{(\ell)}, v^{(\ell)})\\
\end{equation*}

Our goal is to estimate the upper bound of the primal-dual gap
\begin{equation*}
    \text{gap}_{(\ell)}(u, v) = \mathscr{L}(u^{(\ell)}, v) - \mathscr{L}(u, v^{(\ell + 1)}),
\end{equation*}
and prove that $\text{gap}_{(\ell)}(u, v)$ converges to $0$ as $\ell \to \infty$ for all $u$ and $v$. 
Our strategy is to split the estimation into two parts,
\begin{equation*}
    \text{gap}_{(\ell)}(u, v) = \underbrace{\big(\mathscr{L}(u^{(\ell)}, v) - \mathscr{L}(u^{(\ell)}, v^{(\ell + 1)})\big)}_{\text{dual variable part}} + \underbrace{\big(\mathscr{L}(u^{(\ell)}, v^{(\ell + 1)}) - \mathscr{L}(u, v^{(\ell + 1)})\big)}_{\text{primal variable part}},
\end{equation*}
where the first term refers to the gap of the dual variable, and the second refers to the gap of the primal variable. 

\begin{theorem}\label{primalest}
For the primal variable, we have
\begin{equation*}
    \begin{aligned}
        \mathscr{L}(u^{(\ell)}, v^{(\ell + 1)}) - \mathscr{L}(u, v^{(\ell + 1)}) &\le \frac{\big<q^{(\ell + 1)} - q^{(\ell)}, u - u^{(\ell)}\big>}{\tau}\\
        &\quad - \frac{1}{2C_{\Phi}} \big\|\nabla_{u} \Phi(u^{(\ell)}, v^{(\ell + 1)}) - \nabla_{u} \Phi(u, v^{(\ell + 1)})\big\|^{2}.
    \end{aligned}
\end{equation*}
\end{theorem}

\begin{proof}
In the first step of the iteration~\eqref{A1}, we have
\begin{equation*}
    \frac{u^{(\ell)} - q^{(\ell)}}{\tau} + \partial G(u^{(\ell)}) \ni 0
\end{equation*}
Since $G$ is convex, it satisfies:
\begin{equation}\label{A2}
    G(u) - G(u^{(\ell)}) \ge \frac{\big<q^{(\ell)} - u^{(\ell)}, u - u^{(\ell)}\big>}{\tau},
\end{equation}
Then using the convexity and Lipschitz continuity of $\nabla_{u}\Phi(u, v)$, we have
\begin{equation}\label{A3}
    \begin{aligned}
        \Phi(u^{(\ell)}, v^{(\ell + 1)}) - \Phi(u, v^{(\ell + 1)}) &\le \big<\nabla_{u} \Phi(u^{(\ell)}, v^{(\ell + 1)}), u^{(\ell)} - u\big>\\
        &\quad - \frac{1}{2C_{\Phi}} \big\|\nabla_{u} \Phi(u^{(\ell)}, v^{(\ell + 1)}) - \nabla_{u} \Phi(u, v^{(\ell + 1)})\big\|^{2}.
    \end{aligned}
\end{equation}
This estimation is given by Yan \citep[Lemma 1]{Yan2016}, where $C_{\Phi}$ is the Lipschitz constant of $\nabla_u \Phi(u, v)$ such that:
\begin{equation*}
    \big\|\nabla_{u} \Phi(u^{(\ell)}, v) - \nabla_{u} \Phi(u, v)\big\| \le C_{\Phi}\|u^{(\ell)} - u\|.
\end{equation*}

Furthermore, notice that
\begin{equation*}
    \mathscr{L}(u^{(\ell)}, v^{(\ell + 1)}) - \mathscr{L}(u, v^{(\ell + 1)}) = \Phi(u^{(\ell)}, v^{(\ell + 1)}) - \Phi(u, v^{(\ell + 1)}) + G(u^{(\ell)}) - G(u),
\end{equation*}
and therefore, by combining~\eqref{A2} and~\eqref{A3}, we obtain:
\begin{equation}\label{A4}
    \begin{aligned}
        \mathscr{L}(u^{(\ell)}, v^{(\ell + 1)}) - \mathscr{L}(u, v^{(\ell + 1)}) &\le \big<\nabla_{u} \Phi(u^{(\ell)}, v^{(\ell + 1)}), u^{(\ell)} - u\big>\\
        &\quad - \frac{1}{2C_{\Phi}} \big\|\nabla_{u} \Phi(u^{(\ell)}, v^{(\ell + 1)}) - \nabla_{u} \Phi(u, v^{(\ell + 1)})\big\|^{2}\\
        &\quad + \frac{\big<q^{(\ell)} - u^{(\ell)}, u^{(\ell)} - u\big>}{\tau}.
    \end{aligned}
\end{equation}

Recall the definition of $q$, given by
\begin{equation*}
    q^{(\ell + 1)} = u^{(\ell)} - \tau\nabla_{u} \Phi(u^{(\ell)}, v^{(\ell + 1)}).
\end{equation*}
By replacing $u^{(\ell)}$ in terms of $q^{(\ell + 1)}$ in the last term in~\eqref{A4}, we proved the theorem.
\end{proof}

\begin{theorem}\label{dualest}
For the dual variable, we have
\begin{equation*}
    \begin{aligned}
        \mathscr{L}(u^{(\ell)}, v) - \mathscr{L}(u^{(\ell)}, v^{(\ell + 1)}) &\le \frac{\big<v^{(\ell + 1)} - v^{(\ell)}, v - v^{(\ell + 1)}\big>_{T(u^{(\ell)})}}{\tau}\\
        &\quad - \big<q^{(\ell + 1)} - q^{(\ell)}, \mathcal{K}(u^{(\ell)})(v - v^{(\ell + 1)})\big>\\
        &\quad - \frac{\beta}{2}\|v^{(\ell + 1)} - v\|^{2},
    \end{aligned}
\end{equation*}
where $\big<v_{1}, v_{2}\big>_{T} = \big<v_{1}, Tv_{2}\big>$, and $T(u^{(\ell)}) = \frac{\tau}{\sigma}\big(I - \tau \sigma \mathcal{K}(u^{(\ell)})\mathcal{K}(u^{(\ell)})^{T}\big)$.
\end{theorem}

\begin{remark}\label{remark-app}
The semi-norm $\langle v, v\rangle_{T}^{\frac{1}{2}}$ is well-defined as long as the matrix $T(u^{(\ell)})$ is positive semi-definite, which gives us the following convergence criteria
\begin{equation*}
    \tau \sigma \mathcal{K}(u^{(\ell)})\mathcal{K}(u^{(\ell)})^{T} \preceq I \Longrightarrow \tau\sigma \le \min_{u \in \textnormal{dom}(G)}\frac{1}{\|\mathcal{K}(u)\mathcal{K}(u)^{T}\|_{\ast}}.
\end{equation*}
\end{remark}

\begin{proof}
In the second step of the iteration, we note that
\begin{equation*}
    \mathcal{K}(u^{(\ell)})\overline{u}^{(\ell)} + L(u^{(\ell)}) - \frac{(v^{(\ell + 1)} - v^{(\ell)})}{\sigma} \ni \partial F(v^{(\ell + 1)}).
\end{equation*}
Since $F$ is $\beta$-strongly convex, we have:
\begin{equation}\label{A5}
    \begin{aligned}
        F(v^{(\ell + 1)}) - F(v) &\le \frac{\big<v^{(\ell + 1)} - v^{(\ell)}, v - v^{(\ell + 1)}\big>}{\sigma}\\
        &\quad - \big<\mathcal{K}(u^{(\ell)})\overline{u}^{(\ell)} + L(u^{(\ell)}), v - v^{(\ell + 1)}\big>\\
        &\quad - \frac{\beta}{2}\|v^{(\ell + 1)} - v\|^{2}\\
        &= \frac{\big<v^{(\ell + 1)} - v^{(\ell)}, v - v^{(\ell + 1)}\big>_{T(u^{(\ell)})}}{\tau} - \big<L(u^{(\ell)}), v - v^{(\ell + 1)}\big>\\
        &\quad - \big<q^{(\ell + 1)} - q^{(\ell)}, \mathcal{K}(u^{(\ell)})^{T}(v - v^{(\ell + 1)})\big>\\
        &\quad - \big<\mathcal{K}(u^{(\ell)})u^{(\ell)}, v - v^{(\ell + 1)}\big>\\
        &\quad - \frac{\beta}{2}\|v^{(\ell + 1)} - v\|^{2}.
    \end{aligned}
\end{equation}
This estimation is given by Yan \citep[equation (32)]{Yan2016}.

Furthermore, notice that
\begin{equation*}
    \mathscr{L}(u^{(\ell)}, v) - \mathscr{L}(u^{(\ell)}, v^{(\ell + 1)}) = \Phi(u^{(\ell)}, v) - \Phi(u^{(\ell)}, v^{(\ell + 1)}) + F(v^{(\ell + 1)}) - F(v),
\end{equation*}
then from~\eqref{A5}, we obtain:
\begin{equation}\label{A6}
    \begin{aligned}
        \mathscr{L}(u^{(\ell)}, v) - \mathscr{L}(u^{(\ell)}, v^{(\ell + 1)}) &\le \Phi(u^{(\ell)}, v) - \Phi(u^{(\ell)}, v^{(\ell + 1)})\\
        &\quad + \frac{\big<v^{(\ell + 1)} - v^{(\ell)}, v - v^{(\ell + 1)}\big>_{T(u^{(\ell)})}}{\tau}\\
        &\quad - \big<L(u^{(\ell)}), v - v^{(\ell + 1)}\big>\\
        &\quad - \big<q^{(\ell + 1)} - q^{(\ell)}, \mathcal{K}(u^{(\ell)})^{T}(v - v^{(\ell + 1)})\big>\\
        &\quad - \big<\mathcal{K}(u^{(\ell)})u^{(\ell)}, v - v^{(\ell + 1)}\big>\\
        &\quad - \frac{\beta}{2}\|v^{(\ell + 1)} - v\|^{2}.
    \end{aligned}
\end{equation}
Since $L(u^{(\ell)}) = \nabla_{v}\Phi(u^{(\ell)}, v) - \mathcal{K}(u^{(\ell)})u^{(\ell)}$ and $\Phi$ is affine in $v$, we have
\begin{equation}\label{A7}
    \Phi(u^{(\ell)}, v) - \Phi(u^{(\ell)}, v^{(\ell + 1)}) = \big<\mathcal{K}(u^{(\ell)})u^{(\ell)}, v - v^{(\ell + 1)}\big> + \big<L(u^{(\ell)}), v - v^{(\ell + 1)}\big>.
\end{equation}
Replacing the first line in~\eqref{A6} with~\eqref{A7}, the result follows. 
\end{proof}

By combining the above two estimations (Theorem~\ref{primalest} and Theorem~\ref{dualest}), we can show that
\begin{equation}\label{A8}
    \begin{aligned}
        \mathscr{L}(u^{(\ell)}, v) - \mathscr{L}(u, v^{(\ell + 1)}) &\le \frac{\big<q^{(\ell + 1)} - q^{(\ell)}, u - u^{(\ell)}\big>}{\tau} - \frac{\beta}{2}\|v^{(\ell + 1)} - v\|^{2}\\
        &\quad - \frac{1}{2C_{\Phi}} \big\|\nabla_{u} \Phi(u^{(\ell)}, v^{(\ell + 1)}) - \nabla_{u} \Phi(u, v^{(\ell + 1)})\big\|^{2}\\
        &\quad + \frac{\big<v^{(\ell + 1)} - v^{(\ell)}, v - v^{(\ell + 1)}\big>_{T(u^{(\ell)})}}{\tau}\\
        &\quad - \big<q^{(\ell + 1)} - q^{(\ell)}, \mathcal{K}(u^{(\ell)})^{T}(v - v^{(\ell + 1)})\big>.
    \end{aligned}
\end{equation}

In order to estimate the primal-dual gap, we introduce the following Pythagoras identity (see \citep[Fact 3]{Drori2015}).
\begin{lemma}\label{theorem8} For any positive definite matrix $T$, we have
\begin{equation*}
    \big<w - v, T(u - v)\big> = \frac{1}{2}\|w - v\|^{2}_{T} - \frac{1}{2}\|w - u\|_{T}^{2} + \frac{1}{2}\|u - v\|_{T}^{2}
\end{equation*}
where $\|v\|_{T}$ is a semi-norm defined by $\|v\|_{T}^{2} = \langle v, Tv\rangle$.
\end{lemma}

\begin{theorem}\label{theorem7}
By applying the Pythagoras identity to~\eqref{A8}, we have the following estimation of the primal-dual gap:
\begin{equation*}
    \begin{aligned}
        \mathscr{L}(u^{(\ell)}, v) - \mathscr{L}(u, v^{(\ell + 1)})
        &\le \frac{\|q^{(\ell)} - q\|^{2}}{2\tau} - \frac{\|q^{(\ell + 1)} - q\|^{2}}{2\tau}\\
        &\quad + \frac{\|v^{(\ell)} - v\|^{2}_{T(u^{(\ell)})}}{2\tau} - \frac{\|v^{(\ell + 1)} - v\|^{2}_{T(u^{(\ell)})}}{2\tau}\\
        &\quad - \bigg(\frac{1}{2\tau} - \frac{C_{\Phi}}{2}\bigg) \|q^{(\ell)} - q^{(\ell + 1)}\|^{2}\\
        &\quad - \frac{1}{2\tau}\|v^{(\ell)} - v^{(\ell + 1)}\|^{2}_{T(u^{(\ell)})} - \frac{\beta}{2}\|v^{(\ell + 1)} - v\|^{2}\\
        &\quad -\big<q^{(\ell + 1)} - q^{(\ell)}, (\mathcal{K}(u^{(\ell)}) - \mathcal{K}(u))^{T}(v - v^{(\ell + 1)})\big>.
    \end{aligned}
\end{equation*}
where $q = u - \tau \nabla_{u} \Phi(u, v)$ depending on both $u$ and $v$.
\end{theorem}

\begin{proof}
Since
\begin{equation*}
    q^{(\ell + 1)} = u^{(\ell)} - \tau \nabla_{u}\Phi(u^{(\ell)}, v^{(\ell + 1)}) \text{ and } q = u - \tau \nabla_{u} \Phi(u, v),
\end{equation*}
we have
\begin{equation}\label{A10}
    \begin{aligned}
        - \big<q^{(\ell + 1)} - q^{(\ell)}, \mathcal{K}(u^{(\ell)})&(v - v^{(\ell + 1)})\big> = \frac{\big<q^{(\ell + 1)} - q^{(\ell)}, q - q^{(\ell + 1)}\big>}{\tau}\\
        &\quad - \frac{\big<q^{(\ell + 1)} - q^{(\ell)}, u - u^{(\ell)}\big>}{\tau}\\
        &\quad - \big<q^{(\ell + 1)} - q^{(\ell)}, \nabla_{u}\Phi(u^{(\ell)}, v^{(\ell + 1)}) - \nabla_{u}\Phi(u, v^{(\ell + 1)})\big>\\
        &\quad - \big<q^{(\ell + 1)} - q^{(\ell)}, (\mathcal{K}(u^{(\ell)}) - \mathcal{K}(u))^{T}(v - v^{(\ell + 1)})\big>.
    \end{aligned}
\end{equation}
This resembles \citep[the equation following equation (36)]{Yan2016}.

Replacing the last term in~\eqref{A8} with~\eqref{A10}, it becomes:
\begin{equation}\label{A13}
    \begin{aligned}
        \mathscr{L}(u^{(\ell)}, v) - \mathscr{L}(u, v^{(\ell + 1)}) &\le \frac{\big<q^{(\ell + 1)} - q^{(\ell)}, q - q^{(\ell + 1)}\big>}{\tau}\\
        &\quad + \frac{\big<v^{(\ell + 1)} - v^{(\ell)}, v - v^{(\ell + 1)}\big>_{T(u^{(\ell)})}}{\tau} - \frac{\beta}{2}\|v^{(\ell + 1)} - v\|^{2}\\
        &\quad - \big<(q^{(\ell + 1)} - q^{(\ell)}), \nabla_{u}\Phi(u^{(\ell)}, v^{(\ell + 1)}) - \nabla_{u}\Phi(u, v^{(\ell + 1)})\big>\\
        &\quad - \frac{1}{2C_{\Phi}}\big\|\nabla_{u}\Phi(u^{(\ell)}, v^{(\ell + 1)}) - \nabla_{u}\Phi(u, v^{(\ell + 1)})\big\|^{2}\\
        & \quad - \big<q^{(\ell + 1)} - q^{(\ell)}, (\mathcal{K}(u^{(\ell)}) - \mathcal{K}(u))^{T}(v - v^{(\ell + 1)})\big>.
    \end{aligned}
\end{equation}
Using the Young's inequality, we have
\begin{equation}\label{A12}
    \begin{aligned}
        -\big<q^{(\ell + 1)} - q^{(\ell)}&, \nabla_{u}\Phi(u^{(\ell)}, v^{(\ell + 1)}) - \nabla_{u}\Phi(u, v^{(\ell + 1)})\big>\\
        &\quad - \frac{1}{2C_{\Phi}}\big\|\nabla_{u}\Phi(u^{(\ell)}, v^{(\ell + 1)}) - \nabla_{u}\Phi(u, v^{(\ell + 1)})\big\|^{2}\\
        &\quad \le \frac{C_{\Phi}}{2}\|q^{(\ell)} - q^{(\ell + 1)}\|^{2}.
    \end{aligned}
\end{equation}

Combining~\eqref{A13},~\eqref{A12}, and~\eqref{A11}, we see that
\begin{equation*}
    \begin{aligned}
        \mathscr{L}(u^{(\ell)}, v) - \mathscr{L}(u, v^{(\ell + 1)}) &\le \frac{\big<q^{(\ell + 1)} - q^{(\ell)}, q - q^{(\ell + 1)}\big>}{\tau}\\
        &\quad + \frac{\big<v^{(\ell + 1)} - v^{(\ell)}, v - v^{(\ell + 1)}\big>_{T(u^{(\ell)})}}{\tau} \\
        &\quad - \frac{\beta}{2}\|v^{(\ell + 1)} - v\|^{2}+ \frac{C_{\Phi}}{2}\|q^{(\ell)} - q^{(\ell + 1)}\|^{2} \\
        &\quad -\big<q^{(\ell + 1)} - q^{(\ell)}, (\mathcal{K}(u^{(\ell)}) - \mathcal{K}(u))^{T}(v - v^{(\ell + 1)})\big>.
    \end{aligned}
\end{equation*}
Finally, using the Pythagoras identity in Lemma~\ref{theorem8}, we have
\begin{equation}\label{A14}
    \big<q^{(\ell + 1)} - q^{(\ell)}, q - q^{(\ell + 1)}\big> = \frac{\|q^{(\ell)} - q\|^{2}}{2\tau} - \frac{\|q^{(\ell + 1)} - q\|^{2}}{2\tau} -\frac{\|q^{(\ell + 1)} - q^{(\ell)}\|^{2}}{2\tau}.
\end{equation}
\begin{equation}\label{A15}
    \begin{aligned}
        \big<v^{(\ell + 1)} - v^{(\ell)}, v - v^{(\ell + 1)}\big>_{T(u^{(\ell)})} &= \frac{\|v^{(\ell)} - v\|^{2}_{T(u^{(\ell)})}}{2\tau} - \frac{\|v^{(\ell + 1)} - v\|^{2}_{T(u^{(\ell)})}}{2\tau}\\
        &\quad - \frac{\|v^{(\ell + 1)} - v^{(\ell)}\|^{2}_{T(u^{(\ell)})}}{2\tau}.
    \end{aligned}
\end{equation}
By replacing the first two terms with~\eqref{A14} and~\eqref{A15}, we proved theorem.
\end{proof}

Hereinafter, we denote $s = (q, v)$, and define the norm $\|s\|_{I, T}^{2} = \|q\|^{2} + \|v\|_{T}^{2}$. The upper bound of the primal-dual gap in Theorem~\ref{theorem7} can be compactly written as follows,
\begin{equation}\label{A17}
    \begin{aligned}
        \mathscr{L}(u^{(\ell)}, v) - \mathscr{L}(u, v^{(\ell + 1)}) &\le \frac{\|s^{(\ell)} - s\|^{2}_{I, T(u^{(\ell)})}}{2\tau} - \frac{\|s^{(\ell + 1)} - s\|^{2}_{I, T(u^{(\ell)})}}{2\tau}\\
        &\quad - \bigg(\frac{1}{2\tau} - \frac{C_{\Phi}}{2}\bigg) \|q^{(\ell)} - q^{(\ell + 1)}\|^{2} - \frac{\beta}{2}\|v^{(\ell + 1)} - v\|^{2}\\
        &\quad - \big<q^{(\ell + 1)} - q^{(\ell)}, (\mathcal{K}(u^{(\ell)}) - \mathcal{K}(u))^{T}(v - v^{(\ell + 1)})\big>,
    \end{aligned}
\end{equation}
where we have omitted the term $-(2\tau)^{-1}\|v^{(\ell)} - v^{(\ell + 1)}\|^{2}_{T(u^{(\ell)})}$ for further relaxation.

\begin{remark}[Comparison with PD3O] In the problem of PD3O \cite{Yan2016}, since $T(u)\equiv T$ (or equivalently $\mathcal{K}(u)\equiv \mathcal{K}$) is a fixed operator independent of $u$, the gap estimate~\eqref{A17} is sufficient for the proof of convergence if $\tau\le 1/C_{\Phi}$, dispensing with the need of assuming the strong convexity of $F$ (i.e., $\beta=0$). However, in our situation, we need to take into account the difference between $T(u^{(\ell)})$ and $T(u^{(\ell + 1)})$ (or $\mathcal{K}(u^{(\ell)}) - \mathcal{K}(u)$). Therefore, we have to assume strong convexity of $F$ and the boundedness of $\mathcal{K}$ to complete the estimate of the primal-dual gap. 
\end{remark}

\begin{theorem}
By using the bounded property of the operator $\mathcal{K}(u)$,
\begin{equation*}
    \big\|\mathcal{K}(u^{(\ell)})^{T}(v^{(\ell + 1)} - v)\big\| \le C_{K}\|v^{(\ell + 1)} - v\|,
\end{equation*}
we have the following primal-dual gap estimate for any constant $C_0>0$
\begin{equation}\label{A27}
    \begin{aligned}
        \mathscr{L}(u^{(\ell)}, v) - \mathscr{L}(u, v^{(\ell + 1)}) &\le \frac{\|s^{(\ell)} - s\|^{2}_{I, T(u^{(\ell)})}}{2\tau} - \frac{\|s^{(\ell + 1)} - s\|^{2}_{I, T(u^{(\ell + 1)})}}{2\tau}\\
        &\quad - \bigg(\frac{1}{2\tau} - \frac{C_{\Phi}}{2} - \frac{C_{0}}{2}\bigg) \|q^{(\ell)} - q^{(\ell + 1)}\|^{2}\\
        &\quad - \bigg(\frac{\beta}{2} - \frac{\tau C_{K}^{2}}{2} - \frac{2 C_{K}^{2}}{C_{0}} \bigg) \|v^{(\ell + 1)} - v\|^{2}
    \end{aligned}
\end{equation}
\end{theorem}
\begin{proof}
Recall that 
\begin{equation*}
\begin{aligned}
    \| s\|^{2}_{I, T(u)} = \|q\|^{2} + \|v\|_{T(u)}^{2} &= \big<q,q\big>+\big<v,\frac{\tau}{\sigma}\big(I - \tau \sigma \mathcal{K}(u)\mathcal{K}(u)^{T}\big)v\big> \\
    &=\big<q,q\big>+\frac{\tau}{\sigma}\big<v,v\big>-\tau^2 \big<v,\mathcal{K}(u)\mathcal{K}(u)^{T}v\big>
\end{aligned}
\end{equation*}
By using the boundedness of $\mathcal{K}(u)$, we first obtain the following estimate 
\begin{equation}\label{A19}
    \begin{aligned}
        &\|s^{(\ell + 1)} - s\|^{2}_{I, T(u^{(\ell+ 1)})} - \|s^{(\ell + 1)} - s\|^{2}_{I, T(u^{(\ell)})}\\
        &= \tau^{2}\big<v^{(\ell + 1)} - v, (\mathcal{K}(u^{(\ell)})\mathcal{K}(u^{(\ell)})^{T} - \mathcal{K}(u^{(\ell + 1)})\mathcal{K}(u^{(\ell + 1)})^{T})(v^{(\ell + 1)} - v)\big>\\
        &= \tau^{2}\big\|\mathcal{K}(u^{(\ell)})^{T}(v^{(\ell + 1)} - v)\big\|^{2} - \tau^{2}\big\|\mathcal{K}(u^{(\ell + 1)})^{T}(v^{(\ell + 1)} - v)\big\|^{2}\\
        &\le \tau^{2}C_{K}^{2}\|v^{(\ell + 1)} - v\|^{2}.
    \end{aligned}
\end{equation}
Furthermore, we can deal with the last term of~\eqref{A17} as follows
\begin{equation}\label{A11}
    \begin{aligned}
        -\big<q^{(\ell + 1)} - q^{(\ell)}&, (\mathcal{K}(u^{(\ell)}) - \mathcal{K}(u))^{T}(v - v^{(\ell + 1)})\big>\\
        &\le \|q^{(\ell + 1)} - q^{(\ell)}\| \|(\mathcal{K}(u^{(\ell)}) - \mathcal{K}(u))^{T}(v - v^{(\ell + 1)})\|\\
        &\le 2C_{K}\|q^{(\ell + 1)} - q^{(\ell)}\| \|v^{(\ell + 1)} - v\| \\
        &\le \frac{C_{0}}{2} \|q^{(\ell)} - q^{(\ell + 1)}\|^{2} + \frac{2 C_{K}^{2}}{C_{0}} \|v^{(\ell + 1)} - v\|^{2}
    \end{aligned}
\end{equation}
for any $C_0>0$ by using the Young's inequality. Inserting \eqref{A19} and \eqref{A11} into \eqref{A17} leads to the result.
\end{proof}

Finally, we prove the convergence theorem in the ergodic sense given in Theorem~\ref{convergence} by choosing a suitable $\tau$ and $C_0$ to eliminate the last two terms of \eqref{A27}.

\begin{proof}
Here, we consider appropriate values for $C_{0}$ and $\tau$ such that
\begin{equation}\label{A99}
    \frac{1}{2\tau} - \frac{C_{\Phi}}{2} - \frac{C_{0}}{2} \ge 0,\quad \frac{\beta}{2}-\frac{\tau C_{K}^{2}}{2} - \frac{2 C_{K}^{2}}{C_{0}} \ge 0.
\end{equation}
In particular, the first equation in~\eqref{A99} gives us a limitation on the step size such that
\begin{equation*}
    \tau \le \frac{1}{C_{\Phi} + C_{0}}
\end{equation*}
By taking this bound into the second inequality in~\eqref{A99}, we have
\begin{equation*}
    \frac{C_{K}^{2}}{2(C_{\Phi} + C_{0})} + \frac{2 C_{K}^{2}}{C_{0}} \le \frac{\beta}{2}.
\end{equation*}
Since the LHS of the above equation decreases to $0$ as $C_{0}$ increases, we can always find a large enough $C_{0}$ that satisfies the above inequality as long as $\beta > 0$. 

Hence, we can show that
\begin{equation}\label{A31}
    \mathscr{L}(u^{(\ell)}, v) - \mathscr{L}(u, v^{(\ell + 1)}) \le \frac{\|s^{(\ell)} - s\|^{2}_{I, T(u^{(\ell)})}}{2\tau} - \frac{\|s^{(\ell + 1)} - s\|^{2}_{I, T(u^{(\ell + 1)})}}{2\tau}.
\end{equation}
for some constant $C_{0}$. Summing~\eqref{A31} from $\ell = 0$ to $\ell = N$, we have
\begin{equation*}
    \begin{aligned}
        \frac{1}{N}\sum_{\ell = 0}^{N}\big(\mathscr{L}(u^{(\ell)}, v) - \mathscr{L}(u, v^{(\ell + 1)})\big) &\le \frac{\|s^{(0)} - s\|^{2}_{I, T(u^{(0)})}}{2N\tau} - \frac{\|s^{(N + 1)} - s\|^{2}_{I, T(u^{(N + 1)})}}{2N\tau}\\
        &\le \frac{\|s^{(0)} - s\|^{2}_{I, T(u^{(0)})}}{2N\tau}.
    \end{aligned}
\end{equation*}
Finally, since $\mathscr{L}(u, v)$ is convex with respect to $u$ and affine with respect to $v$, by the Jensen's inequality, we have
\begin{equation*}
    \begin{aligned}
        \big(\mathscr{L}(u_{(N)}, v) - \mathscr{L}(u, v_{(N)})\big) &\le \frac{1}{N}\sum_{\ell = 0}^{N}\big(\mathscr{L}(u^{(\ell)}, v) - \mathscr{L}(u, v^{(\ell + 1)})\big)\\
        &\le \frac{\|s^{(0)} - s\|^{2}_{I, T(u^{(0)})}}{2N\tau}
    \end{aligned}
\end{equation*}
for all $s = (q, v)$, where $u_{(N)} = \sum_{\ell = 0}^{N}u^{(\ell)}/N$ and $v_{(N)} = \sum_{\ell = 0}^{N}v^{(\ell + 1)}/N$. Passing the limit $N \to \infty$, we can see that the primal-dual gap converges to $0$, and hence the ergodic sequence $(u_{(N)}, v_{(N)})$ converges to the saddle point.
\end{proof}
\end{appendix}

\bibliographystyle{amsplain}
\bibliography{bibliography.bib}

@article{Lisini2012,
    title = {{Cahn-Hilliard} and thin film equations with nonlinear mobility as gradient flows in weighted-{Wasserstein} metrics},
    author={Lisini, Stefano and Matthes, Daniel and Savar{\'e}, Giuseppe},
    journal={J. Differential Equations},
    volume={253},
    number={2},
    pages={814-850},
    year={2012}
}

@article{Carrillo2023,
    author = {Carrillo, Jos\'{e} A. and Wang, Li and Wei, Chaozhen},
    title = {Structure Preserving Primal Dual Methods for Gradient Flows with Nonlinear Mobility Transport Distances},
    journal = {SIAM J. Numer. Anal.},
    volume = {62},
    number = {1},
    pages = {376-399},
    year = {2024}
}

@article{Jordan1996,
    title={The variational formulation of the {Fokker-Planck} equation},
    author={Jordan, Richard and Kinderlehrer, David and Otto, Felix},
    journal={SIAM J. Math. Anal.},
    volume={29},
    number={1},
    pages={1-17},
    year={1998}
}

@article{Carrillo2003,
    title={Kinetic equilibration rates for granular media and related equations: entropy dissipation and mass transportation estimates},
    author={Jos{\'e} Antonio Carrillo and Robert J. McCann and C{\'e}dric Villani},
    journal={Rev. Mat. Iberoam.},
    volume={19},
    number={3},
    pages={971-1018},
    year={2003}
}

@book{Ambrosio2005,
    title={Gradient Flows: In Metric Spaces and in the Space of Probability Measures},
    author={Luigi Ambrosio and Nicola Gigli and Giuseppe Savar{\'e}},
    year={2005},
    publisher = {Birkh\"auser},
    address = "Basel"
}

@article{Burger2006,
    title={The {Keller-Segel} Model for Chemotaxis with Prevention of Overcrowding: Linear vs. Nonlinear Diffusion},
    author={Martin Burger and Marco Di Francesco and Yasmin Dolak-Struss},
    journal={SIAM J. Math. Anal.},
    year={2006},
    volume={38},
    number={4},
    pages={1288-1315}
}

@article{Cahn1996,
    title={The {Cahn–Hilliard} equation with a concentration dependent mobility: motion by minus the {Laplacian} of the mean curvature},
    author={John W. Cahn and Charles M. Elliott and Amy Novick-Cohen},
    journal={European J. Appl. Math.},
    year={1996},
    volume={7},
    number={3},
    pages={287 - 301}
}

@article{Elliott1996,
    title={On the {Cahn-Hilliard} equation with degenerate mobility},
    author={Charles M. Elliott and Harald Garcke},
    journal={SIAM J. Math. Anal.},
    year={1996},
    volume={27},
    number={2},
    pages={404-423}
}

@article{Bertozzi1998,
    title={The mathematics of moving contact lines in thin liquid films},
    author={Bertozzi, Andrea L},
    journal={Notices Amer. Math. Soc.},
    volume={45},
    number={6},
    pages={689-697},
    year={1998}
}

@incollection{Bertozzi2002,
    title={Thin film dynamics: theory and applications},
    author={Bertozzi, Andrea L and Bowen, Mark},
    booktitle={Modern Methods in Scientific Computing and Applications},
    pages={31-79},
    year={2002},
    publisher= {Springer Netherlands},
    address= {Dordrecht},
}

@article{Dolbeault2008,
    title={A new class of transport distances between measures},
    author={Jean Dolbeault and Bruno Nazaret and Giuseppe Savar{\'e}},
    journal={Calc. Var. Partial Differential Equations},
    year={2008},
    volume={34},
    number={2},
    pages={193-231}
}

@article{Benamou2000,
    title={A computational fluid mechanics solution to the {Monge-Kantorovich} mass transfer problem},
    author={Jean-David Benamou and Yann Brenier},
    journal={Numer. Math.},
    year={2000},
    volume={84},
    number={3},
    pages={375-393}
}

@article{Zhu2022,
    title={New Primal-Dual Algorithms for a Class of Nonsmooth and Nonlinear Convex-Concave Minimax Problems},
    author={Yuzixuan Zhu and Deyi Liu and Quoc Tran-Dinh},
    journal={SIAM J. Optim.},
    year={2022},
    volume={32},
    number={4},
    pages={2580-2611}
}

@book{Manzoni2021,
    author = {Andrea Manzoni and Alfio Maria Quarteroni and Sandro Salsa},
    title = {Optimal Control of Partial Differential Equations: Analysis, Approximation, and Applications},
    address = "Cham",
    publisher = {Springer},
    year = {2021},
}

@book{Vzquez2006,
    author = {Juan Luis V{\'a}zquez},
    title = {The Porous Medium Equation},
    publisher = {Oxford University Press},
    address = "Northants",
    year = {2006}
}

@article{Oron1997,
    title={Long-scale evolution of thin liquid films},
    author={Alexander Oron and S. H. Davis and S. George Bankoff},
    journal={Rev. Modern Phys.},
    year={1997},
    volume={69},
    pages={931-980}
}

@article{Carrillo2019,
    title={Primal Dual Methods for {Wasserstein} Gradient Flows},
    author={Jos{\'e} Antonio Carrillo and Katy Craig and Li Wang and Chaozhen Wei},
    journal={Found. Comput. Math.},
    year={2021},
    volume={22},
    pages={389 - 443}
}

@book{Boffi2013,
    author = {Boffi, Daniele and Brezzi, Franco and Fortin, Michel},
    title = {Mixed Finite Element Methods and Applications},
    publisher = {Springer},
    address = "Berlin",
    year = {2013}
}

@article{Eyre1998,
    title={Unconditionally gradient stable time marching the {Cahn--Hilliard} equation},
    author={David J. Eyre},
    journal={MRS Online Proceedings Library (OPL)},
    volume={529},
    number={39},
    pages={1676409},
    year={1998}
}

@article{beretta1995nonnegative,
    title={Nonnegative solutions of a fourth-order nonlinear degenerate parabolic equation},
    author={Beretta, Elena and Bertsch, Michiel and Dal Passo, Roberta},
    journal={Archive for rational mechanics and analysis},
    volume={129},
    number={2},
    pages={175--200},
    year={1995},
    publisher={Springer-Verlag}
}

@article{Lee2023,
    title = {Deep {JKO}: Time-implicit particle methods for general nonlinear gradient flows},
    journal = {J. Comput. Phys.},
    volume = {514},
    pages = {113187},
    year = {2024},
    issn = {0021-9991},
    author = {Wonjun Lee and Li Wang and Wuchen Li}
}

@article{Drori2015,
    title={A simple algorithm for a class of nonsmooth convex-concave saddle-point problems},
    author={Yoel Drori and Shoham Sabach and Marc Teboulle},
    journal={Oper. Res. Lett.},
    year={2015},
    volume={43},
    number={2},
    pages={209-214}
}

@article{Clason2016,
    title={Primal-Dual Extragradient Methods for Nonlinear Nonsmooth {PDE}-Constrained Optimization},
    author={Christian Clason and Tuomo Valkonen},
    journal={SIAM J. Optim.},
    year={2017},
    volume={27},
    number={3},
    pages={1314-1339}
}

@article{Westdickenberg2008,
    title={VARIATIONAL PARTICLE SCHEMES FOR THE POROUS MEDIUM EQUATION AND FOR THE SYSTEM OF ISENTROPIC {Euler} EQUATIONS},
    author={Michael Westdickenberg and Jon Wilkening},
    journal={ESAIM Math. Model. Numer. Anal.},
    year={2010},
    volume={44},
    number={1},
    pages={133-166}
}

@article{Li2019,
    title = {Fisher information regularization schemes for {Wasserstein} gradient flows},
    journal = {J. Comput. Phys.},
    volume = {416},
    pages = {109449},
    year = {2020},
    issn = {0021-9991},
    author = {Wuchen Li and Jianfeng Lu and Li Wang}
}

@article{Carrillo2009,
    title={Nonlinear mobility continuity equations and generalized displacement convexity},
    author={Jos{\'e} Antonio Carrillo and Stefano Lisini and Giuseppe Savar'e and Dejan Slepvcev},
    journal={J. Funct. Anal.},
    year={2010},
    volume={258},
    number={4},
    pages={1273-1309}
}

@article{Carrillo2022,
    author = {Bailo, Rafael and Carrillo, Jos\'{e} A. and Hu, Jingwei},
    title = {Bound-Preserving Finite-Volume Schemes for Systems of Continuity Equations with Saturation},
    journal = {SIAM J. Appl. Math.},
    volume = {83},
    number = {3},
    pages = {1315-1339},
    year = {2023}
}

@article{Carrillo2014,
    title={A Finite-Volume Method for Nonlinear Nonlocal Equations with a Gradient Flow Structure},
    author={Jos{\'e} Antonio Carrillo and Alina Chertock and Yanghong Huang},
    journal={Commun. Comput. Phys.},
    year={2015},
    volume={17},
    number={1},
    pages={233-258}
}

@article{Blanchet2008,
    title={Convergence of the Mass-Transport Steepest Descent Scheme for the Subcritical {Patlak-Keller-Segel} Model},
    author={Adrien Blanchet and Vincent Calvez and Jos{\'e} Antonio Carrillo},
    journal={SIAM J. Numer. Anal.},
    year={2008},
    volume={46},
    pages={691-721}
}

@article{Otto2001,
    title={THE GEOMETRY OF DISSIPATIVE EVOLUTION EQUATIONS: THE POROUS MEDIUM EQUATION},
    author={Felix Otto},
    journal={Commun. Partial Differ. Equ.},
    year={2001},
    volume={26},
    pages={101 - 174}
}

@article{Chambolle2019,
    title={Total roto-translational variation},
    author={A. Chambolle and Thomas Pock},
    journal={Numer. Math.},
    year={2019},
    volume={142},
    pages={611 - 666}
}

@article{Papadakis2013,
    title={Optimal Transport with Proximal Splitting},
    author={Nicolas Papadakis and Gabriel Peyr{\'e} and {\'E}douard Oudet},
    journal={SIAM J. Imaging Sci.},
    year={2014},
    volume={7},
    number={1},
    pages={212-238}
}

@article{Hamedani2018,
    author = {Hamedani, Erfan Yazdandoost and Aybat, Necdet Serhat},
    title = {A Primal-Dual Algorithm with Line Search for General Convex-Concave Saddle Point Problems},
    journal = {SIAM J. Optim.},
    volume = {31},
    number = {2},
    pages = {1299-1329},
    year = {2021}
}

@article{Rtz2006,
    title={Surface evolution of elastically stressed films under deposition by a diffuse interface model},
    author={Andreas R{\"a}tz and Angel Ribalta and Axel Voigt},
    journal={J. Comput. Phys.},
    year={2006},
    volume={214},
    pages={187-208}
}

@article{Li2023,
    title={High order spatial discretization for variational time implicit schemes: {Wasserstein} gradient flows and reaction-diffusion systems},
    author={Guosheng Fu and Stanley J. Osher and Wuchen Li},
    journal={J. Comput. Phys.},
    year={2023},
    volume={491},
    pages={112375}
}

@article{Chen1997,
    title={Convergence Rates in Forward-Backward Splitting},
    author={George H.-G. Chen and R. Tyrrell Rockafellar},
    journal={SIAM J. Optim.},
    year={1997},
    volume={7},
    number={2},
    pages={421-444}
}

@article{Sun2017,
    title={A discontinuous {Galerkin} method for nonlinear parabolic equations and gradient flow problems with interaction potentials},
    author={Zheng Sun and Jos{\'e} Antonio Carrillo and Chi-Wang Shu},
    journal={J. Comput. Phys.},
    year={2017},
    volume={352},
    pages={76-104}
}

@article{Shen2018,
    title={The scalar auxiliary variable ({SAV}) approach for gradient flows},
    author={Jie Shen and Jie Xu and Jiang Yang},
    journal={J. Comput. Phys.},
    year={2018},
    volume={353},
    pages={407-416}
}

@article{Xu2006,
    title={Stability Analysis of Large Time-Stepping Methods for Epitaxial Growth Models},
    author={Chuanju Xu and Tao Tang},
    journal={SIAM J. Numer. Anal.},
    year={2006},
    volume={44},
    number={4},
    pages={1759-1779}
}

@article{Yang2017,
    title={Linear, second order and unconditionally energy stable schemes for the viscous {Cahn-Hilliard} equation with hyperbolic relaxation using the invariant energy quadratization method},
    author={Xiaofeng Yang and Jia Zhao and Xiaoming He},
    journal={J. Comput. Appl. Math.},
    year={2017},
    volume={343},
    number={4},
    pages={80-97}
}

@article{Frank2020,
    title = {Bound-preserving flux limiting schemes for {DG}-discretizations of conservation laws with applications to the {Cahn–Hilliard} equation},
    journal = {Comput. Methods Appl. Mech. Eng.},
    volume = {359},
    pages = {112665},
    year = {2020},
    author = {Florian Frank and Andreas Rupp and Dmitri Kuzmin}
}

@article{Cheng20212,
    title={A new {Lagrange} multiplier approach for constructing structure-preserving schemes, {II}. {Bound} preserving},
    author={Qing Cheng and Jie Shen},
    journal={SIAM J. Numer. Anal.},
    year={2021},
    volume={60},
    number={3},
    pages={970-998}
}

@article{DeGiorgi,
    title={New problems on minimizing movements},
    author={De Giorgi, Ennio},
    journal={Ennio de Giorgi: selected papers},
    pages={699--713},
    year={1993},
    publisher={Springer}
}

@article{Lu2012,
    title={The cutoff method for the numerical computation of nonnegative solutions of parabolic {PDEs} with application to anisotropic diffusion and Lubrication-type equations},
    author={Changna Lu and Weizhang Huang and Erik S. Van Vleck},
    journal={J. Comput. Phys.},
    year={2012},
    volume={242},
    pages={24-36}
}

@article{Bretin2020,
    author = {Bretin, Elie and Masnou, Simon and Sengers, Arnaud and Terii, Garry},
    title = {Approximation of surface diffusion flow: A second-order variational {Cahn–Hilliard} model with degenerate mobilities},
    journal = {Math. Models Methods Appl. Sci.},
    volume = {32},
    number = {04},
    pages = {793-829},
    year = {2022}
}

@article{Malitsky2018,
    title={A Forward-Backward Splitting Method for Monotone Inclusions Without Cocoercivity},
    author={Yura Malitsky and Matthew K. Tam},
    journal={SIAM J. Optim.},
    year={2020},
    volume={30},
    number={2},
    pages={1451-1472}
}

@article{Zhornitskaya1999,
    title={Positivity-Preserving Numerical Schemes for Lubrication-Type Equations},
    author={L. Zhornitskaya and A. Bertozzi},
    journal={SIAM J. Numer. Anal.},
    year={1999},
    volume={37},
    number={2},
    pages={523-555}
}

@article{Cancs2019,
    title={A variational finite volume scheme for {Wasserstein} gradient flows},
    author={Cl{\'e}ment Canc{\`e}s and Thomas O. Gallou{\"e}t and Gabriele Todeschi},
    journal={Numer. Math.},
    year={2019},
    volume={146},
    pages={437 - 480}
}

@article{Jacobs2018,
    title={Solving Large-Scale Optimization Problems with a Convergence Rate Independent of Grid Size},
    author={Matt Jacobs and Flavien L{\'e}ger and Wuchen Li and Stanley J. Osher},
    journal={SIAM J. Numer. Anal.},
    year={2019},
    volume={57},
    number={3},
    pages={1100-1123}
}

@article{Pock2011,
    title={Diagonal preconditioning for first order primal-dual algorithms in convex optimization},
    author={Thomas Pock and A. Chambolle},
    journal={2011 International Conference on Computer Vision},
    year={2011},
    pages={1762-1769}
}

@article{Chambolle2011,
    title={A First-Order Primal-Dual Algorithm for Convex Problems with Applications to Imaging},
    author={A. Chambolle and Thomas Pock},
    journal={J. Math. Imaging Vision},
    year={2011},
    volume={40},
    pages={120-145}
}

@article{Chambolle2016,
    title={On the ergodic convergence rates of a first-order primal–dual algorithm},
    author={A. Chambolle and Thomas Pock},
    journal={Math. Program.},
    year={2016},
    volume={159},
    number={1},
    pages={253-287}
}

@article{Zhou2023,
    title={An Operator-Splitting Optimization Approach for Phase-Field Simulation of Equilibrium Shapes of Crystals},
    author={Zeyu Zhou and Wen Huang and Wei Jiang and Zhen Zhang},
    journal={SIAM J. Sci. Comput.},
    year={2024},
    volume={46},
    number={3},
    pages={331-353}
}

@article{Davis2015,
    title={A Three-Operator Splitting Scheme and its Optimization Applications},
    author={Damek Davis and Wotao Yin},
    journal={Set-Valued and Variational Analysis},
    year={2017},
    volume={25},
    pages={829-858}
}

@article{Bailo2021,
    title={Bound-Preserving Finite-Volume Schemes for Systems of Continuity Equations with Saturation},
    author={Rafael Bailo and Jos{\'e} Antonio Carrillo and Jingwei Hu},
    journal={SIAM J. Appl. Math.},
    year={2023},
    volume={83},
    number={3},
    pages={1315-1339}
}

@article{Sun2018,
    title={An entropy stable high-order discontinuous {Galerkin} method for cross-diffusion gradient flow systems},
    author={Zheng Sun and Jos{\'e} Antonio Carrillo and Chi-Wang Shu},
    journal={Kinetic \& Related Models},
    year={2018}
}

@article{Cheng2020,
    title={A weakly nonlinear, energy stable scheme for the strongly anisotropic {Cahn-Hilliard} equation and its convergence analysis},
    author={Kelong Cheng and Cheng Wang and Steven M. Wise},
    journal={J. Comput. Phys.},
    year={2020},
    volume={405},
    pages={109109}
}

@article{Yan2016,
    title={A New Primal–Dual Algorithm for Minimizing the Sum of Three Functions with a Linear Operator},
    author={Ming Yan},
    journal={J. Sci. Comput.},
    year={2018},
    volume={76},
    pages={1698 - 1717}
}

@article{Liu2021AccelerationPrimalDual,
	title={Acceleration of primal-dual methods by preconditioning and simple subproblem procedures},
	journal={J. Sci. Comput.},
	volume={86},
	number={2},
	pages={1-34},
	year={2021},
	author={Yanli Liu and Yunbei Xu and Wotao Yin}
}

@article{chen2023TPD,
    title={Transformed primal-dual methods with variable-preconditioners},
    author={Chen, Long and Guo, Ruchi and Wei, Jingrong},
    journal={arXiv preprint arXiv:2312.12355},
    year={2023}
}

@article{Du2021Siamreview,
author = {Du, Qiang and Ju, Lili and Li, Xiao and Qiao, Zhonghua},
title = {Maximum Bound Principles for a Class of Semilinear Parabolic Equations and Exponential Time-Differencing Schemes},
journal = {SIAM Review},
volume = {63},
number = {2},
pages = {317-359},
year = {2021},
}

@article{Chen2019CH-log,
title = {Positivity-preserving, energy stable numerical schemes for the Cahn-Hilliard equation with logarithmic potential},
author = {Wenbin Chen and Cheng Wang and Xiaoming Wang and Steven M. Wise},
journal = {Journal of Computational Physics: X},
volume = {3},
pages = {100031},
year = {2019},
issn = {2590-0552},
doi = {https://doi.org/10.1016/j.jcpx.2019.100031},
}

@article{Hu2020SSIPNP,
	title={{A fully discrete positivity-preserving and energy-dissipative finite difference scheme for Poisson-Nernst-Planck equations}},
	author={Hu, Jingwei and Huang, Xiaodong},
	journal={Numerische Mathematik},
	volume={145},
	number={1},
	pages={77-115},
	year={2020},
	publisher={Springer}
}

@article{Shen2021PNP,
	title={{Unconditionally positivity preserving and energy dissipative schemes for Poisson-Nernst-Planck equations}},
	author={Shen, Jie and Xu, Jie},
	journal={Numerische Mathematik},
	volume={148},
	number={3},
	pages={671-697},
	year={2021},
	publisher={Springer}
}

\end{document}